\documentclass[12pt]{article}
\usepackage{amsfonts}
\usepackage{amsmath}
\usepackage{pdfsync}

\usepackage[english]{babel}
\usepackage[latin1]{inputenc}
\usepackage[T1]{fontenc}

\newtheorem{theorem}{Theorem}
\newtheorem{definition}{Definition}
\newtheorem{lemma}{Lemma}

\newtheorem{corollary}{Corollary}
\newtheorem{proposition}{Proposition}

\numberwithin{equation}{section}

\begin{document}

\centerline{\large{\textbf{Hidden processes and hidden Markov processes: classical and quantum}}}

\bigskip\bigskip

\centerline{\textbf{Luigi Accardi}}

\centerline{Centro Vito Volterra, Universit\`a  di Roma ''Tor Vergata'',}
\centerline{Roma I-00133, Italy, accardi@volterra.uniroma2.it}

\centerline{\textbf{Soueidy El Gheteb}}
\centerline{Department of Mathematics and Informatics, Faculty of Sciences and Technologies,}
\centerline{University of Nouakchott Al Aasriya,}
\centerline{Nouakchott, Mauritania, elkotobmedsalem@gmail.com}

\centerline{\textbf{Yun Gang Lu}}
\centerline{Dipartimento di Matematica, Univerversit\`{a} di Bari ``Aldo Moro'',}
\centerline{Via E. Orabona 4, 70125, Bari, Italy,  E-mail:yungang.lu@uniba.it}

\centerline{\textbf{Abdessattar Souissi }}
\centerline{$^1$ Department of Accounting, College of Business Management}
\centerline{Qassim University, Ar Rass, Saudi Arabia, a.souaissi@qu.edu.sa }
\centerline{$^2$ Preparatory Institute for Scientific and Technical Studies La Marsa,}
\centerline{ Carthage University, Tunisia, abdessattar.souissi@ipest.rnu.tn}

\bigskip\bigskip

\begin{abstract}
This paper consists of $3$ parts. The first part only considers classical processes and
introduces two different extensions of the notion of hidden Markov process.
In the second part, the notion of quantum hidden process is introduced.
In the third part it is proven that, by restricting various types of quantum Markov chains
to appropriate commutative sub--algebras (diagonal sub--algebras) one recovers all
the classical hidden process and, in addition, one obtains families of processes which are
not usual hidden Markov process, but are included in the above mentioned extensions of
these processes.
In this paper we only deal with processes with an at most countable state space.
\end{abstract}

\tableofcontents
\eject

\section{ Introduction}

Classical hidden Markov processes (HMP) were introduced in 1966 by Baum and Petrie \cite{BaumPetr66} in order to extend \textit{the standard statistical estimation theory for independent sampling or Markov chains} (in particular, the maximum likelihood estimate and the $(\chi^2)$--theory of power, estimation and testing) to a larger class of functions of Markov processes, now called hidden
Markov processes.
They define this larger class of functions of Markov processes considering two correlated
processes, an \textit{observable} one and a Markov process (now called a \textit{hidden} (or \textit{underlying}) process). Their idea was to acquire information on the hidden process
through measurements of the observable one.\\
In their original construction, the processes in this new class differ from the usual Markov
processes by the presence of the so--called \textit{emission operator}, a Markov operator
describing the conditional probabilities of the observable process given the underlying one
which is described by its own Markov operator.
These two Markov operators, together with the initial distribution, uniquely define the joint
probabilities of the process through an assumption of conditional independence (see Definition \ref{df:class-HP} below).\\
The availability of additional parameters, compared to the usual Markov chains,
makes this new class of processes natural candidates for modelling complex systems
and in fact in the past years these processes have found a multiplicity of applications,
in speech recognition \cite{JelBahMer75}, \cite{RabLeSo83}, \cite{Rab89}, \cite{HuaYasMerv90},
security of cloud and credit cards \cite{RabJua86}, \cite{Algh2016},
computational molecular biology \cite{FelsChur92}, \cite{JasKell12},
pattern classification and recognition \cite{SatGuru93},
machine learning \cite{GhahrJord97},
computer vision \cite{GhahrJord97}, \cite{RosPent98},
models of human interactions \cite{RosPent98}, \cite{YamOhyIsh92},
bioinformatics  \cite{Eddy98}, genetics \cite{LiSte03},
musical performances \cite{PardBirm05}, economics and finance \cite{HasNat05}, \cite{Nguyet18},
\cite{RebSa17}.
A nice survey of applications of classical HMMs is \cite{MGK2021}.\\
We will use the term HMPs when dealing with the general probabilistic structure underlying all the
specific models and, when referring to the latter ones, we use the term hidden Markov models
(HMMs).\\

\noindent Several articles have been dedicated in recent years to the attempt to extend various
aspects of the definition of HMPs to the quantum case. In \cite{WiesnCrutc08} the evolution of
the system is governed by the application of quantum operations on a quantum state.
The sequence of emitted symbols defines the sequence of quantum operations being applied on
the initial state of the hidden quantum Markov model.\\
Monras et al. \cite{MonrWiesn11} propose to describe a class of Quantum HMM (HQMMs), related to
problems in many-body systems, in terms of a set of quantum operations associated with emission
symbols.\\
In \cite{SatGuru93} a maximum-likelihood-based unsupervised learning algorithm was proposed for HQMMs
in the framework of quantum information. \\
However these are only partial extensions in the sense that they deal only with some aspects
of HMP, like dynamics or different kinds of statistical algorithms, but
\textbf{no class of quantum processes has so far been produced} with the properties that:\\
(i) its restriction to some abelian sub-algebra of the algebra of the process gives a classical
HMP;\\
(ii) varying the quantum process and the abelian sub-algebra, one \textbf{can recover all
classical HMP}.\\

\noindent The aim of this article is to describe a solution of this problem which:\\
(i) not only \textbf{includes all the classical HMP} (by restriction to abelian sub-algebras),\\
\noindent (ii) but also suggests \textbf{several natural extensions of the classical HMP} in the
sense that we produce \textbf{easily constructable} new classes of classical processes
\textbf{not covered by the existing literature on classical HMP}.\\
Many of such new types of processes can be obtained by taking diagonal restrictions of
quantum Markov Chains and quantum HMP.\\
In Section \ref{sec:Alg-form-clas-stoc-proc}, we recall same basic ideas on the algebraic
formulation of classical stochastic processes. Classical hidden processes are introduced
in Section \ref{sec:Clas-HP} and they are shown to include classical hidden Markov processes.
In Section \ref{sec:Alg-clas-q-stoc-proc}, the second part of the paper, the notion of
algebraic stochastic process, including both the classical and the quantum case, is recalled
and the notion of quantum hidden process is introduced.
In the third part (Section \ref{Diag-QMC-Cl-HMP}) it is proven that, by restricting various
types of quantum Markov chains to appropriate commutative sub--algebras (diagonal sub--algebras)
one recovers all the classical hidden process and, in addition, one obtains families of processes which are not usual hidden Markov process, but are included in the extensions of these processes introduced in the previous sections.

\section{Algebraic formulation of classical stochastic processes}\label{sec:Alg-form-clas-stoc-proc}

The algebraic formulation of classical stochastic processes has two advantages: (1) it
better hilights the structure of the various processes; (2) it greatly simplifies the transition
from classical to quantum processes.\\
In the following, when no confusion is possible, we denote a measurable space
$(S,\mathcal{B})$ simply by the symbol $S$, \textbf{leaving implicit the} $\sigma$--algebra.
\begin{definition}\label{df:class-stoc-proc}{\rm
Let $T$ be a set and $(S_{n})_{n\in T}$ a family of measurable spaces.
A \textbf{classical stochastic process} with \textbf{state spaces} $S_{n}$ and \textbf{index set} $T$
is given by a probability space $(\Omega,\mathcal{F}, P)$ and family of \textbf{measurable functions} $X\equiv (X_{n})_{n\in T}$
$$
X_{n}:(\Omega,\mathcal{F})\to S_{n}
$$
Denote $\mathcal{F}_{fin}(T)$ the family of finite sub--sets of $T$ (if $F\in\mathcal{F}_{fin}(T)$,
we often simply write $F\subset_{fin} T$).
The \textbf{finite dimensional joint expectations} of the process $X$ (also called
\textbf{correlation kernels} or \textbf{correlation functions} in physics) are the expectation values
\begin{equation}\label{df-class-joint-exp}
P_{F}\left((f_{n})_{n\in F}\right)
:= P_{X}\left(\prod_{n\in F} f_{n}(X_{n})\right)
:= \int_{\Omega} \prod_{n\in F} f_{n}(X_{n}(\omega)) P_{X}(d\omega)
\end{equation}
for any $F\in\mathcal{F}_{fin}(T)$ and $f_{n}\in L^{\infty}_{\mathbb{C}}(S_{n})$
for each $n\in F$.\\
Two classical stochastic processes with state spaces $S_{n}$ and index set $T$ are called
\textbf{stochastically equivalent} if they have \textbf{the same} finite dimensional joint expectations. Any element in the stochastic equivalence class of a process $X$ is called a
\textbf{realization} of $X$.
}\end{definition}
\textbf{Remark}.
Choosing in \eqref{df-class-joint-exp} each $f_{n}$ to be the characteristic
function  $\chi_{I_{n}}$ of some measurable sub--set $I_{n}\subseteq S_{n}$, one obtains the
\textbf{finite dimensional joint probabilities} which also characterize the stochastic
equivalence class of a process $X$.\\
\noindent The finite dimensional joint expectations \eqref{df-class-joint-exp} are states, i.e.
positive, normalized linear functionals,
\begin{equation}\label{df-class-proj-exp}
P_{F} \colon  \mathcal{A}_{F}
:=\bigotimes_{n\in F} L^{\infty}_{\mathbb{C}}(S_{n})  \to \mathbb{C}
\quad,\quad F\in\mathcal{F}_{fin}(T)
\end{equation}
with the following properties:\\
1) If $F\subset G\subset_{fin} T$, there is a natural immersion
$$
j_{F,G} \colon f_{F} \in \mathcal{A}_{F} \to j_{F,G}(f_{F})\in \mathcal{A}_{G}
$$
where by definition
$$
j_{F,G}(f_{F})(\{x_{n}\}_{n\in G}) := f_{F}(\{x_{n}\}_{n\in F})
$$
We consider this immersion as an identification, where a function of the variables
$\{x_{n}\}_{n\in F}$ is considered as a function in the variables $\{x_{n}\}_{n\in G}$ constant in the variables $\{x_{n}\}_{n\in G\setminus F}$. Therefore we simply write
$\mathcal{A}_{F} \subseteq \mathcal{A}_{G}$.\\
2) For any $F\subset_{fin} T$, $P_{F}$ is a state on
$\mathcal{A}_{F}=\prod_{n\in F} L^{\infty}_{\mathbb{C}}(S_{n})$.\\
3) If $F\subset G\subset_{fin} T$,
\begin{equation}\label{df-class-proj-exp-alg}
P_{G}\big|_{\mathcal{A}_{F}} = P_{F}
\end{equation}
\begin{definition}{\rm
A family $(P_{F})_{F\subset_{fin} T}$ of functions of the form \eqref{df-class-proj-exp-alg}
satisfying conditions 2) and 3) above, is called a \textbf{projective family of
distribution functions} with state spaces $S_{n}$ and index set $\mathcal{F}_{fin}(T)$.
}\end{definition}
A fundamental result of probability theory is the following.
\begin{theorem}\label{Kolm-consist-thm}{\rm
(Kolmogorov consistency theorem)
If $S_{n}=S=\mathbb{R}^d$, for some $d\in\mathbb{N}$ and for all $n\in\mathbb{N}$,
given any projective family of distribution functions $(P_{F})_{F\subset_{fin} T}$ with state spaces $S$, there exists a unique probability measure $P$ on the space
$\Omega:=\prod_{n\in T}S\equiv \mathcal{F}(T;S)$ with the $\sigma$--algebra $\mathcal{F}$
generated by the projection functions
$$
\pi_{n}: \omega:= (s_{n})_{n\in T}\in\Omega \to \pi_{n}(\omega) := s_{n}\in S
\quad,\quad n\in T
$$
such that the finite dimensional joint expectations of the stochastic process
$\pi_{n} \colon (\Omega, \mathcal{F})\to S $ coincide with $(P_{F})_{F\subset_{fin} T}$.
}\end{theorem}
\textbf{Remark}.
Kolmogorov consistency (or compatibility) theorem holds for spaces much more general than
$S=\mathbb{R}^d$. The subtle issue, in the transition from cylindrical measure to measure, is
the proof of countable additivity.
To our knowledge the sharpest result in this direction is \cite{RaoSaz93-proj} which also contains
a complete bibliography.\\

\noindent A corollary of Kolmogorov consistency theorem is that any classical stochastic
process with state space $S=\mathbb{R}^d$ has a \textbf{Kolmogorov realization} on the probability
space defined by Theorem \ref{Kolm-consist-thm}.\\

\subsection{Algebraic classical stochastic process}\label{sec:alg-class-stoch-proc}

Any classical stochastic process $(X_{n})_{n\in T}$ with index set $T$ and state spaces $S_n$
($n\in T$) uniquely defines:\\
-- the $*$--homomorphisms
\begin{equation}\label{df-jn-class}
j_{n} \colon f\in L^{\infty}_{\mathbb{C}}(S_n) \to j_{n}(f)
:= f(X_{n})\in L^{\infty}_{\mathbb{C}}(\Omega,\mathcal{F}, P) \quad,\quad n\in T
\end{equation}
-- the state $P_{X}$ on $L^{\infty}_{\mathbb{C}}(\Omega,\mathcal{F}, P)$ given by
$$
P_{X}(F) := \int_{\Omega} F(\omega) P(d\omega)
\quad,\quad F\in L^{\infty}_{\mathbb{C}}(\Omega,\mathcal{F}, P)
$$
\begin{definition}\label{df:alg-class-stoch-proc}{\rm
The quadruple
\begin{equation}\label{alg-class-stoch-proc}
\left(L^{\infty}_{\mathbb{C}}(\Omega,\mathcal{F}),P_{X},  (L^{\infty}_{\mathbb{C}}(S_n))_{n\in T},
(j_{n})_{n\in T}\right)
\end{equation}
is called an \textbf{classical algebraic stochastic process} with index set $T$, \textbf{sample
algebra} $L^{\infty}_{\mathbb{C}}(\Omega,\mathcal{F}, P)$ and \textbf{state algebra}
$L^{\infty}_{\mathbb{C}}(S_n)$.
}\end{definition}
\textbf{Remark}.
The most general notion of classical algebraic stochastic process is obtained replacing, in
Definition \ref{df:alg-class-stoch-proc}, $L^{\infty}_{\mathbb{C}}(\Omega,\mathcal{F}, P)$ and
$L^{\infty}_{\mathbb{C}}(S_n)$ by arbitrary abelian $*$--algebras.\\

\noindent Using the $*$--homomorphisms \eqref{df-jn-class} to identify, for each $F\subset_{fin} T$,
the\\ $*$--algebra $\prod_{n\in F} L^{\infty}_{\mathbb{C}}(S_{n})$ with the $*$--sub--algebra of
$L^{\infty}_{\mathbb{C}}(\Omega,\mathcal{F}, P)$ generated by
$\{j_{n}(L^{\infty}_{\mathbb{C}}(S_n)) \colon n\in F\}$, the algebraic expression of the
finite dimensional joint expectations becomes
\begin{equation}\label{df-class-joint-exp-alg}
P_{F}\left(\prod_{n\in F}f_{n}\right)
:=P_{X}\left(\prod_{n\in F}j_{n}(f_{n})\right) \quad,\quad f_{n}\in L^{\infty}_{\mathbb{C}}(S_n)
\ , \ n\in F
\end{equation}
and the stochastic equivalence relation is expressed by the coincidence of the
finite dimensional joint expectations \eqref{df-class-joint-exp-alg}.\\

\noindent An algebraic formulation of the Kolmogorov realization of the algebraic classical
stochastic process \eqref{alg-class-stoch-proc} is obtained replacing the $*$--algebra
$L^{\infty}_{\mathbb{C}}(\Omega,\mathcal{F}, P),P_{X})$ by the dense sub--algebra
$\bigotimes_{T} L^{\infty}_{\mathbb{C}}(S_n)$ consisting of the algebraic tensor product
of $|T|$ copies of $L^{\infty}_{\mathbb{C}}(S_n)$. This leads to the stochastically equivalent
algebraic classical process
\begin{equation}\label{alg-class-stoch-proc-KR}
\left(\bigotimes_{n\in T} L^{\infty}_{\mathbb{C}}(S_{n}),P_{X}), L^{\infty}_{\mathbb{C}}(S_n), (j_{n})_{n\in T}\right)
\end{equation}
where the embeddings are given by
$$
j_{n} \colon f_{n} \in L^{\infty}_{\mathbb{C}}(S_{n}) \to
f_{n}\otimes 1_{T\setminus \{n\}}
$$
where $1_{T\setminus \{n\}}$ is the identity of
$\bigotimes_{k\in T\setminus \{n\}}L^{\infty}_{\mathbb{C}}(S_k)$.

\subsection{Algebraic classical Markov chains}

In the notations of Section \ref{sec:alg-class-stoch-proc} \textbf{but with the index set}
$T:=\mathbb{N}$, the algebraic classical stochastic process $(X_{n})_{n\in \mathbb{N}}$
is called a \textbf{backward Markov chain}, for any $n\in\mathbb{N}$ and
$g_{m}\in L^{\infty}_{\mathbb{C}}(S_{m})$, $m\in\{1,\dots,n\}$, one has
\begin{equation}\label{df-class-MP}
E_{X_{0}, \dots,X_{n}}\left(g_{n+1}(X_{n+1})\right)
= E_{X_{n}}\left(g_{n+1}(X_{n+1})\right)
\in j_{n}\left( L^{\infty}_{\mathbb{C}}(S_{n})\right) =:\mathcal{A}_{X_{n}}
\end{equation}
where $E_{X_{0}, \dots,X_{n}}$ (resp. $E_{X_{n}}$) is the $P_{X}$--conditional expectation onto
the algebra
$$
\mathcal{A}_{X;[0,n]} := \hbox{ algebraic span of }
\{\mathcal{A}_{X_{0}}, \ldots,\mathcal{A}_{X_{n}}\}
 \  (\hbox{resp. } \mathcal{A}_{X_{n}})
$$
It is known that condition \eqref{df-class-MP} is equivalent to the \textbf{conditional
independence identity}
\begin{equation}\label{MP-equiv-cond-ind}
E_{X_{n}}\left(a_{n)}a_{(n}\right)
= E_{X_{n}}\left(a_{n)}\right)E_{X_{n}}\left(a_{(n}\right)
\end{equation}
for any
$$
\begin{cases}
a_{n)}\in \mathcal{A}_{X;n)}
:= \hbox{ algebraic span of }\{\mathcal{A}_{X_{k}} \colon k<n\}  \\
a_{(n}\in \mathcal{A}_{X;(n}
:= \hbox{ algebraic span of }\{\mathcal{A}_{X_{k}} \colon k>n\}
\end{cases}
\quad,\quad\forall n\in\mathbb{N}
$$
If any of the two equivalent conditions \eqref{df-class-MP} or \eqref{MP-equiv-cond-ind}
is satisfied, $P_{X}$ is called a Markov state.
\noindent The structure theorem for classical Markov processes is a corollary of
Theorem \ref{proj-fam-equiv-st} as shown by the following known result.
\begin{theorem}\label{thm:struct-Mark-st}{\rm
Let $(X_{n})$ be a backward Markov chain and define
$$
\mathcal{A}_{X;[0,n]} :=  \hbox{ algebraic span of }\{\mathcal{A}_{X_{k}} \colon k\in [0,n]\}
$$
$$
\mathcal{A} = \bigcup_{n\in\mathbb{N}}\mathcal{A}_{X;[0,n]}
$$
Then $(X_{n})$ uniquely defines a pair $(P_{0}, (P_{n})_{n\in\mathbb{N}})$ such that:\\
-- $P_{0}$ is a state on $\mathcal{A}_{X_{0}}$.\\
-- For each $n\ge 1$, $P_{n} \colon\mathcal{A}_{X_{n}}\to \mathcal{A}_{X_{n-1}}$ is a Markov operator (completely positive identity preserving).\\
-- The family of states defined by
\begin{equation}\label{fin-dim-clas-MC-jnt-exp}
P_{[0,n]}(a_{0}\cdot a_{1}\cdots a_{n})
:= P_{0}\left(a_{0}P_{1}(a_{1}\cdots P_{n-1}(a_{n-1}P_{n}(a_{n})))\right)
\end{equation}
for any $n\in\mathbb{N}$ and $a_{m}\in\mathcal{A}_{X_{m}}$, $m\in\{1,\dots,n\}$, is projective.
In particular the limit
\begin{equation}\label{fin-dim-clas-MC-lim}
\lim_{N\to\infty} P_{[0,N]}(a) =: P(a)
\end{equation}
exists for any $a\in\mathcal{A}$ in the strongly finite sense.\\
Conversely, any pair $(P_{0}, (P_{n})_{n\in\mathbb{N}})$ as above uniquely defines, through
\eqref{fin-dim-clas-MC-jnt-exp} the finite dimensional joint expectations of a unique classical
Markov Chain.
}\end{theorem}

\noindent\textbf{Remark}.
In what follows we will say that the state $P_{X}$ on $(\mathcal{A}_{X}$, defined by Corollary
\eqref{fin-dim-clas-MC-lim}, is defined by the pair $(P_{0}, (P_{n})_{n\in\mathbb{N}})$.\\
$P_{X}$ is uniquely determined by its joint expectations \eqref{fin-dim-clas-MC-jnt-exp} up to
stochastic equivalence and states whose joint expectations have the form \eqref{fin-dim-clas-MC-lim}
are called \textbf{backward Markov states}.

\section{Classical hidden processes }\label{sec:Clas-HP}

Recall that, in the present paper, \textbf{only the case} $T=\mathbb{N}$ is discussed.\\
The following definition extends the notion of hidden Markov processes to the case where
\textbf{the hidden process is arbitrary}.
\begin{definition}\label{df:class-HP}{\rm
Given two sequences of measurable spaces $(S_{H_{n}})_{n\in\mathbb{N}}$
$(S_{O_{n}})_{n\in\mathbb{N}}$ denote
\begin{equation}\label{df-st-alg-n-cl}
\mathcal{B}_{H_{n}} := L^{\infty}_{\mathbb{C}}(S_{H_{n}}) \quad;\quad
\mathcal{B}_{O_{n}} := L^{\infty}_{\mathbb{C}}(S_{O_{n}})
\end{equation}
the corresponding state algebras and define:\\
-- the $H$--sample algebra
$$
\mathcal{A}_{H}
:=\bigotimes_{n\in \mathbb{N}}\mathcal{B}_{H_{n}}
$$
-- the $O$--sample algebra
$$
\mathcal{A}_{O}
:=\bigotimes_{n\in \mathbb{N}} \mathcal{B}_{O_{n}}
$$
-- the corresponding tensor embeddings ($n\in \mathbb{N}$)
\begin{equation}\label{df-H-embeds-cl}
j_{H_{n}}(g_{n})\equiv g(H_{n})
:= g_{n}\otimes 1_{\bigotimes_{m\in \{n\}^c}\mathcal{B}_{H_{m}}}
\quad,\quad\forall g_{n}\in \mathcal{B}_{H_{n}}
\end{equation}
\begin{equation}\label{df-O-embeds-cl}
j_{O_{n}}(f) \equiv f(O_{n}) := f_{n}\otimes 1_{\bigotimes_{m\in \{n\}^c}\mathcal{B}_{O_{m}}}
\quad,\quad\forall f_{n}\in \mathcal{B}_{O_{n}}
\end{equation}
where, for a set $I\subseteq T$, $I^c$ denotes its complement;\\
-- the $(H,O)$--sample algebra
\begin{equation}\label{df-HO-sampl-alg-cl}
\mathcal{A}_{H,O}
:=\mathcal{A}_{H}\otimes \mathcal{A}_{O}
\equiv \bigotimes_{n\in \mathbb{N}}
(\mathcal{B}_{H_{n}}\otimes \mathcal{B}_{O_{n}})
\end{equation}
-- the corresponding tensor embeddings ($n\in \mathbb{N}$)
\begin{equation}\label{df-H-embeds-cl-alg}
j_{H_{n}}\otimes j_{O_{n}} := \mathcal{B}_{H_{n}}\otimes \mathcal{B}_{O_{n}}
\to \mathcal{A}_{H,O}
\end{equation}
Let $P_{H,O} $ be a state on $\mathcal{A}_{H,O}$.
The classical stochastic process\\ ($(H,O)$--process)
\begin{equation}\label{df-HO-proc-clas}
\left( \mathcal{A}_{H,O}, P_{H,O},
(\mathcal{B}_{H_{n}}\otimes \mathcal{B}_{O_{n}})_{n\in \mathbb{N}},
(j_{H_{n}}\otimes j_{O_{n}})_{n\in \mathbb{N}}  \right)
\end{equation}
is called a \textbf{hidden process} if, introducing the $H$--sub--process,
called the \textbf{hidden or underlying process}
\begin{equation}\label{df-H-proc-cl}
\left( \mathcal{A}_{H}\equiv\mathcal{A}_{H}\otimes 1_{\mathcal{A}_{O}} ,
P_{H}\equiv P_{H,O}\big|_{\mathcal{A}_{H}\otimes 1_{\mathcal{A}_{O}}}, (j_{H_{n}}),
(\mathcal{B}_{H_{n}})  \right)
\end{equation}
($j_{H_{n}}\equiv j_{H_{n}}\otimes 1_{\mathcal{A}_{O}}$,
$\mathcal{B}_{H_{n}}\equiv \mathcal{B}_{H_{n}}\otimes 1_{\mathcal{A}_{O}}$)
and the $O$--sub--process, also called the \textbf{observable process}
\begin{equation}\label{df-O-proc-cl}
\left( \mathcal{A}_{O}\equiv 1_{\mathcal{A}_{H}}\otimes \mathcal{A}_{O} ,
P_{O}\equiv P_{H,O}\big|_{1_{\mathcal{A}_{H}}\otimes \mathcal{A}_{O}},
(j_{O_{n}}), (L^{\infty}_{\mathbb{C}}(S_{O_{n}})  \right)
\end{equation}
($j_{O_{n}}\equiv 1_{\mathcal{A}_{H}}\otimes  j_{O_{n}}$,
$\mathcal{B}_{O_{n}}\equiv 1_{\mathcal{A}_{H}}\otimes \mathcal{B}_{O_{n}}$),
the $O$--process is \textbf{conditionally independent} of the $H$--process
in the following sense.\\
For all $n\in\mathbb{N}$, the $P_{H,O}$--conditional expectation $E_{H_{0}, \dots,H_{n}}$ onto
the algebra
$$
\mathcal{A}_{H_{0}, \dots,H_{n}}
:= \bigvee_{m\in\{0,1,\dots,n\}}j_{H_{m}}\left(\mathcal{B}_{H_{m}}\right)
=: \bigvee_{m\in\{0,1,\dots,n\}} \mathcal{A}_{H_{m}}
$$
satisfies
\begin{equation}\label{cond-exps-O|H}
E_{H_{0}, \dots,H_{n}}\left(\prod_{m=0}^{n}j_{O_{m}}(f_{m})\right)
=\prod_{m=0}^{n}E_{H_{m}}\left(j_{O_{m}}(f_{m})\right)
\end{equation}
for all $f_{m}\in \mathcal{B}_{O_{m}}$ ($m\in \{0, \dots, n\}$), where $E_{H_{m}}$ denotes the
$P_{H,O}$--conditional expectation onto $\mathcal{A}_{H_{m}}$.
}\end{definition}
\noindent\textbf{Remark}.
Notice that, interpreting the index set $\mathbb{N}$ as time, in \eqref{cond-exps-O|H}, for
each $m\in\mathbb{N}$, the observable process $O_{m}$ is conditioned on the hidden process
$H_{m}$ \textbf{at the same time}.
These processes differ from the \textbf{backward hidden Markov process} (see Definition \ref{df:class-backw-HMP}) in which the observable process $O_{m}$ is conditioned on the
hidden process $H_{m-1}$ \textbf{at the previous time}.
\begin{definition}\label{df:stoc-eq-class-HMP}{\rm
Two classical hidden processes are called \textbf{stochastically equivalent} if the
\textbf{joint expectations of the associated observable processes coincide}.
}\end{definition}

\noindent Notice that this equivalence relation is specific for hidden processes:
for usual stochastic processes stochastic equivalence is defined as coincidence of \textbf{all}
joint expectations.\\

\noindent\textbf{Remark}.
The restriction of $E_{H_{m}}$ to the algebra
$$
\mathcal{A}_{O_{m}} := j_{O_{m}}(L^{\infty}_{\mathbb{C}}(S_{O_{m}}))
:=\{j_{m}(f_{m}) \colon f_{m}\in L^{\infty}_{\mathbb{C}}(S_{O_{m}})\} \subseteq \mathcal{A}_{O}
$$
is a \textbf{Markov operator} $E_{H_{m}}\colon \mathcal{A}_{O_{m}}
\to \mathcal{A}_{H_{m}}:= j_{O_{m}}(L^{\infty}_{\mathbb{C}}(S_{O_{m}}))$.
Since $j_{H_{m}}$ has a left inverse, the  linear operator
\begin{equation}\label{df-B(O,H,m)2}
B_{O_{m},H_{m}}  = j_{H_{m}}^{-1}\circ E_{H_{m}} \circ j_{O_{m}}
\colon \mathcal{B}_{O_{m}} = L^{\infty}_{\mathbb{C}}(S_{O_{m}})
\to \mathcal{B}_{H_{m}} = L^{\infty}_{\mathbb{C}}(S_{H_{m}})
\end{equation}
is well defined because the range of $E_{H_{m}}$ is in the domain of $j_{H_{m}}^{-1}$.
By construction $B_{O_{m},H_{m}}$, satisfies
\begin{equation}\label{df-B(O,H,m)}
j_{H_{m}}(B_{O_{m},H_{m}}f_{m}) =  E_{H_{m}}(j_{O_{m}}(f_{m}))
\ , \ \forall f_{m}\in \mathcal{B}_{O_{m}}
\end{equation}
\eqref{df-B(O,H,m)2} implies that $B_{O_{m},H_{m}}$ is a Markov operator being a composition of
completely positive identity preserving maps.\\
The Markov operator $B_{O,H,n}$ is called the \textbf{$n$--th emission operator}.
The \textbf{global emission operator} is the Markov operator
\begin{equation}\label{basic-form-HMP}
B_{O,H} := \bigotimes_{n\in \mathbb{N}}B_{O,H,n} \ \colon \
\bigotimes_{n\in \mathbb{N}}L^{\infty}_{\mathbb{C}}(S_{O_{n}}) \to
\bigotimes_{n\in \mathbb{N}}L^{\infty}_{\mathbb{C}}(S_{H_{n}})
\end{equation}
\textbf{Remark}.
Notice that, in Definition \ref{df:class-HP}, the hidden process $H\equiv\{H_n\}_{n\in\mathbb{N}} $, \textbf{is arbitrary}. For example one can take $H\equiv\{H_n\}_{n\in\mathbb{N}} $ to be a
sequence of independent identically distributed \textbf{gaussian} random variables. In this case
one will have a \textbf{hidden gaussian process}. Or one can choose
$H\equiv\{H_n\}_{n\in\mathbb{N}} $ to be a backward $d$--step Markov process (i.e with finite
memory of length $d$). In this case
one will have a \textbf{hidden backward $d$--step Markov process}.
In general infinitely many other choices are possible
\begin{theorem}\label{th:joint-exps-same-time(O,H)}{\rm
Given a hidden process in the sense of Definition \ref{df:class-HP}, the
probability distribution  $P_{H,O}$ of the $(H,O)$--process are uniquely determined by the
probability distribution $P_{H}$ of the $H$--process, the emission probability operators
$(B_{O,H,n})_{n\in\mathbb{N}}$ and the property of \textbf{conditional independence} through
the identity:
\begin{equation}\label{joint-exps-(O,H)-B}
P_{H,O} = P_{H} \circ B_{O,H}
\end{equation}
or, more explicitly:
\begin{align}
&P_{H,O}\left(\prod_{m=0}^{n}f_{m}(O_{m})g_{m}(H_{m})\right)
=P_{H}\left(\prod_{m=0}^{n}j_{H_{m}}(B_{O_{m},H_{m}}f_{m})j_{H_{m}} (g_{m})\right)\notag\\
=: &P_{H_{0},\dots,H_{n}}\left(\prod_{m=0}^{n}j_{H_{m}}((B_{O_{m},H_{m}}f_{m})\cdot g_{m})\right) \label{joint-exps-(O,H)-B-expl}
\end{align}
for all $n\in\mathbb{N}$, $m\in \{0, \dots, n\}$, $f_{m}\in L^{\infty}_{\mathbb{C}}(S_{O_{m}})$ and
$g_{m}\in L^{\infty}_{\mathbb{C}}(S_{H_{m}})$.\\
Conversely, given:\\
(1) a classical stochastic process $H\equiv\{H_n\}_{n\in\mathbb{N}} $ with state spaces $S_{H_{n}}$ and probability distribution $P_{H}$;\\
(2) for each $n\in \mathbb{N}$, a measurable space $S_{O_{n}}$;\\
(3) for each $n\in \mathbb{N}$, a Markov operator
\begin{equation}\label{df-B(O,H,n)}
B_{O,H,n} \colon  L^{\infty}_{\mathbb{C}}(S_{O_{n}}) \to L^{\infty}_{\mathbb{C}}(S_{H_{n}})
\end{equation}
the right hand side of \eqref{joint-exps-(O,H)-B} defines a unique state $P_{H,O}$ on the
$*$--algebra
\begin{equation}\label{df-alg-HO}
\mathcal{A}_{H,O} := \bigotimes_{n\in \mathbb{N}}
\left( L^{\infty}_{\mathbb{C}}(S_{H_{n}})\otimes L^{\infty}_{\mathbb{C}}(S_{O_{n}})\right)
\end{equation}
and a unique classical process $(O_{n})$ such that:\\
(i) the processes $(H_{n})$ and $(O_{n})$ satisfy the conditions of Definition \ref{df:class-HP};\\
(ii) the joint expectations of the $(H,O)$--process are given by \eqref{joint-exps-(O,H)-B}.
In particular, the restriction of $P_{H,O}$ on the $*$--sub--algebra
\begin{equation}\label{df-AH-sub-alg-HO}
\mathcal{A}_{H} \equiv \bigotimes_{n\in \mathbb{N}}
\left( L^{\infty}_{\mathbb{C}}(S_{H_{n}})\otimes 1_{L^{\infty}_{\mathbb{C}}(S_{O_{n}})}\right)
\end{equation}
coincides with $P_{H}$.
}\end{theorem}
\textbf{Proof}.
If the conditions of Definition \ref{df:class-HP} are satisfied, one has:
\begin{align}
&P\left(\prod_{m=0}^{n}f_{m}(O_{m})g_{m}(H_{m})\right)
=P\left(E_{H_{0},\dots,H_{n}}\left(\prod_{m=0}^{n}f_{m} (O_{m})\right)\prod_{m=0}^{n}g_{m}(H_{m})\right)\notag\\
\overset{\eqref{cond-exps-O|H}, \eqref{joint-exps-(O,H)-B-expl}}{=}&
P_{H_{0}, \dots,H_{n}}\left(\prod_{m=0}^{n}E_{H_{m}} \left(f_{m}(O_{m})\right)g_{m}(H_{m})\right) \label{joint-exps-(O,H)-E}
\end{align}
and both $P_{H_{0}, \dots,H_{n}}$ and the $E_{H_{m}}$ are given. Therefore
\eqref{joint-exps-(O,H)-B} follows from \eqref{joint-exps-(O,H)-E} because
$E_{H_{m}}\left(f_{m}(O_{m})\right) \overset{\eqref{df-B(O,H,m)}}{=}  j_{H_{m}}(B_{O_{m},H_{m}}f_{m}) $.\\
Conversely, suppose that conditions (1), (2), (3) are satisfied and use the right hand side of \eqref{joint-exps-(O,H)-B} to define the joint expectations on the left hand side.
The fact that each $B_{O,H,n}$ is a Markov operator allows to verify that these expectations
satisfy the Kolmogorov compatibility conditions for finite dimensional joint expectations.
Hence they define a unique state $P_{H,O}$ on $\mathcal{A}_{H,O}$ whose restriction on
the $*$--sub--algebra \eqref{df-AH-sub-alg-HO} coincides with $P_{H}$.
It remains to verify conditional independence.
This follows from
$$
P_{H,O}\left(\prod_{m=0}^{n}f_{m}(O_{m})g_{m}(H_{m})\right)
=P_{H,O}\left(\prod_{m=0}^{n}f_{m}(O_{m})\prod_{m=0}^{n}g_{m}(H_{m})\right)
$$
$$
\overset{\eqref{joint-exps-(O,H)-B-expl}}{=}P_{H_{0}, \dots,H_{n}}\left(E_{H_{0}, \dots,H_{n}}
\left(\prod_{m=0}^{n}f_{m}(O_{m})\right)\prod_{m=0}^{n}g_{m}(H_{m})\right)
$$
On the other hand, the definition of $P_{H,O}$ implies that
$$
P_{H,O}\left(\prod_{m=0}^{n}f_{m}(O_{m})g_{m}(H_{m})\right)
=P_{H_{0}, \dots,H_{n}}\left(\prod_{m=0}^{n}j_{H_{m}}(B_{O_{m},H_{m}}f_{m})g_{m}(H_{m})\right)
$$
$$
=P_{H_{0}, \dots,H_{n}}
\left(\prod_{m=0}^{n}j_{H_{m}}(B_{O_{m},H_{m}}f_{m})\prod_{m=0}^{n}g_{m}(H_{m})\right)
$$
Therefore
$$
P_{H_{0}, \dots,H_{n}}\left(E_{H_{0}, \dots,H_{n}}
\left(\prod_{m=0}^{n}f_{m}(O_{m})\right)\prod_{m=0}^{n}g_{m}(H_{m})\right)
$$
$$
=P_{H_{0}, \dots,H_{n}}
\left(\prod_{m=0}^{n}j_{H_{m}}(B_{O_{m},H_{m}}f_{m})\prod_{m=0}^{n}g_{m}(H_{m})\right)
$$
Since $\prod_{m=0}^{n}j_{H_{m}}(B_{O_{m},H_{m}}f_{m})$ belongs to the algebra
$\mathcal{A}_{H_{0}, \dots,H_{n}}$ and the products $\prod_{m=0}^{n}g_{m}(H_{m})$ are total in this
algebra, it follows that
\begin{equation}\label{cond-ind-O|H2}
E_{H_{0}, \dots,H_{n}}
\left(\prod_{m=0}^{n}f_{m}(O_{m})\right)
=\prod_{m=0}^{n}j_{H_{m}}(B_{O_{m},H_{m}}f_{m})
\end{equation}
and this proves \eqref{cond-exps-O|H}.
It remains to prove that, for each $m\in\{1,\dots,n\}$
$$
E_{H_{m}}(f_{m}(O_{m})) = j_{H_{m}}(B_{O_{m},H_{m}}f_{m})
$$
This follows by taking all the $f_{k}(O_{k})=1$ for $k\in\{1,\dots,n\}\setminus\{m\}$ in \eqref{cond-ind-O|H2}.
$\qquad\square$

\subsection{Classical hidden Markov processes}

\noindent The following Corollary shows that, when the hidden process is (backward) Markov, the family of states defined by Theorem \ref{th:joint-exps-same-time(O,H)} coincides with the family of usual classical hidden Markov processes defined in the literature (see e.g. \cite{RabJua86}, \cite{Rab89}, \cite{{SatGuru93}}).
\begin{corollary}\label{cor:joint-exps-H-Mark}{\rm
In the assumptions and notations of Theorem \ref{th:joint-exps-same-time(O,H)},
if the\\ $H$--process is the backward Markov process with initial state $p_{H_{0}}$ (state on
$L^{\infty}_{\mathbb{C}}(S_{H_{0}})$) and transition operators
$P_{H_{n}}\colon L^{\infty}_{\mathbb{C}}(S_{H_{n}+1})\to L^{\infty}_{\mathbb{C}}(S_{H_{n}})$,
$n\in\mathbb{N}$, then
\begin{align}\label{joint-probs-O-H-Mark}
&P_{H_{0}}\left(\prod_{m=0}^{n}f_{m}(O_{m})g_{m}(H_{m})\right)\\
= &p_{H_{0}}\left(B_{O,H,0}(f_{0})g_{0}P_{H_{0}}\left(B_{O,H,1} (f_{1})g_{1}\right)\left(\cdots P_{H_{n-1}}\left(B_{O,H,n} (f_{n})g_{n}\right)\right)\right)\notag
\end{align}
}\end{corollary}
\textbf{Proof}. Using the Markov property of the $H$--process, \eqref{joint-exps-(O,H)-B} becomes
\begin{align*}
&P_{H,O}\left(\prod_{m=0}^{n}f_{m}(O_{m})g_{m}(H_{m})\right)
=P_{H}\left(\prod_{m=0}^{n}j_{H_{m}}(B_{O_{m},H_{m}}f_{m})g_{m}(H_{m}) \right)\\
= &P_{H}\left(j_{H_{0}}(B_{O,H,0}f_{0})g_{0}(H_{0})
\prod_{m=1}^{n}j_{H_{m}}(B_{O_{m},H_{m}}f_{m})g_{m}(H_{m})\right)
\end{align*}
and by applying \eqref{fin-dim-clas-MC-jnt-exp} to the Markov process $\{ p_{H_{0}}, P_{H_{n}}\}$, one knows
that the above expression is equal to the right hand side of \eqref{joint-probs-O-H-Mark}.
$\qquad\square$\\

\noindent\textbf{Example}. Suppose that the state spaces of all the observable (resp. hidden)
random variables are equal to a single space $S_{O}$ (resp. $S_{H}$) and that the cardinality of both $S_{H}$ and $S_{O}$ are at most countable.
Then the identity \eqref{joint-probs-O-H-Mark} implies that the joint probabilities of the
$(H,O)$--process are given by
\begin{equation}\label{joint-exp-HO-discr-proc}
P_{H,O} \left(H_0=j_0, O_0=k_0; H_1=j_1, O_1=k_1; \dots; H_n=j_n,  O_n=k_n\right)
\end{equation}
$$
=p_{H;j_0}^{(0)} p_{H;j_0j_1}\cdots  p_{H;j_{n-1}j_n}\cdot
p(O_0=k_0|H_0=j_0)\cdots p(O_n=k_n|H_n=j_n)
$$
for any choice of $j_0, j_1, \dots, j_n \in S_{H}$ and of $k_0; k_1, \dots, k_n\in S_{O}$.\\
As announced at the beginning of this section, the right hand side of the identity
\eqref{joint-exp-HO-discr-proc} is the usual expression for the joint probabilities of a hidden Markov process as can be found in the literature on these processes.

\subsection{Classical hidden Markov processes as restrictions of Classical
Markov processes}

In this section we suppose that:\\
1)  The assumptions of Corollary \ref{cor:joint-exps-H-Mark} are verified.\\
2) For all $n\in\mathbb{N}$, both Markov operators
\begin{equation}\label{df-P(Hn)-P(On)}
P_{H_{n}}\colon L^{\infty}_{\mathbb{C}}(S_{H_{n+1}})\to L^{\infty}_{\mathbb{C}}(S_{H_{n}})
\quad;\quad
B_{O,H,n}\colon L^{\infty}_{\mathbb{C}}(S_{O_{n}})\to L^{\infty}_{\mathbb{C}}(S_{H_{n}})
\end{equation}
are given by regular Markov kernels denoted respectively
$P_{H_{n}}(h_{n}, dh_{n+1})$ and $P_{O_{n}}(h_{n}, do_{n})$ ($h_{n}\in S_{H_{m}}$),
i.e. that for all
$ g_{n}\in L^{\infty}_{\mathbb{C}}(S_{H_{n}})$ and
$ f_{n}\in L^{\infty}_{\mathbb{C}}(S_{O_{n}})$, they have the form
\begin{equation}\label{trans-ops-class-HMP}
(P_{H_{n}}f_{n})(H_{n}) =  \int_{S_{H_{n+1}}}f_{n}(h_{n+1}) P_{H_{n}}(H_{n}, dh_{n+1})
\end{equation}
\begin{equation}\label{trans-ops-class-BOH}
(B_{O_{m},H_{m}}f)(H_{m}) =   \int_{S_{O;m}}f_{m}(o_{m}) B_{O_{m}}(H_{m}, do_{m})
\end{equation}
In what follows, we will often us the simplified notation:
\begin{equation}\label{B(OHn)-equiv-P(On)}
B_{O,H,n} \equiv B_{O_{n}} \, ,\quad  1_{L^\infty_{\mathbb{C}} (S_{O_n})} \equiv 1_{O_n}\ \text{ and }\ 1_{L^\infty_{\mathbb{C}} (S_{H_n})} \equiv 1_{H_n}\ ,\quad\forall n\in\mathbb{N}
\end{equation}
For each $n\ge 1$, define the Markov operator
\begin{equation}\label{act-H-tens-O}
P_{H_{n-1}\otimes O_{n-1};H_{n-1}}\colon
L^{\infty}_{\mathbb{C}}(S_{H_{n}})\otimes L^{\infty}_{\mathbb{C}}(S_{O_{n-1}})
\to L^{\infty}_{\mathbb{C}}(S_{H_{n-1}})
\end{equation}
by
\begin{align}\label{H-tens-O-ops}
&P_{H_{n-1}\otimes O_{n-1};H_{n-1}}(g_{n}\otimes  f_{n-1})(h_{n-1})\\
:=&\left(P_{H_{n-1}}\circ B_{O_{n-1}}\right)(g_{n}\otimes f_{n-1})(h_{n-1})  = (P_{H_{n}}g_{n}\cdot  B_{O_{n-1}}f_{n-1})(h_{n-1})\notag\\
=&\int_{S_{H_{n}}\times S_{O;n-1}}
(P_{H_{n-1}}(h_{n-1}, dh_{n})B_{O_{n-1}}(h_{n-1}, do_{n-1}))g_{n}(h_{n})f_{n-1}(o_{n-1})\notag
\end{align}
With the identification
\begin{equation}\label{identif-m-ge-2}
L^{\infty}_{\mathbb{C}}(S_{H_{n}})
\equiv L^{\infty}_{\mathbb{C}}(S_{H_{n}})\otimes 1_{L^{\infty}_{\mathbb{C}}(S_{O_{n-1}})}
\quad, \ \forall n\ge 1
\end{equation}
For $n\ge 2$, one can consider $P_{H_{n-1}\otimes O_{n-1};H_{n-1}}$ as an operator
\begin{equation}\label{act-H-tens-O-id}
P_{H_{n-1}\otimes O_{n-1};H_{n-1}}\colon
L^{\infty}_{\mathbb{C}}(S_{H_{n}})\otimes L^{\infty}_{\mathbb{C}}(S_{O_{n-1}})
\to L^{\infty}_{\mathbb{C}}(S_{H_{n-1}})\otimes L^{\infty}_{\mathbb{C}}(S_{O_{n-2}})
\end{equation}
whose range satisfies
\begin{equation}\label{range-P(Hn)-O(n-1);H(n-1)}
P_{H_{n-1}\otimes O_{n-1};H_{n-1}}\left(
L^{\infty}_{\mathbb{C}}(S_{H_{n}})\otimes L^{\infty}_{\mathbb{C}}(S_{O_{n-1}})\right)
\subseteq L^{\infty}_{\mathbb{C}}(S_{H_{n-1}})\otimes 1_{L^{\infty}_{\mathbb{C}}(S_{O_{n-2}})}
\end{equation}
Similarly
\begin{equation}\label{range-P(H1)-O(0);H(0)-O(0)}
P_{H_{0}\otimes O_{0};H_{0}}\left(
L^{\infty}_{\mathbb{C}}(S_{H;1})\otimes L^{\infty}_{\mathbb{C}}(S_{O_{0}})\right)
\subseteq L^{\infty}_{\mathbb{C}}(S_{H_{0}})\otimes 1_{L^{\infty}_{\mathbb{C}}(S_{O_{0}})}
\end{equation}
Define further the algebra
\begin{equation}\label{df-A(H-tens-O)}
\mathcal{A}_{H\otimes O} :=
\left(L^{\infty}_{\mathbb{C}}(S_{H_{0}})\otimes L^{\infty}_{\mathbb{C}}(S_{O_{0}})\right)\otimes
\bigotimes_{n\in\mathbb{N}\setminus \{0\}}
\left(L^{\infty}_{\mathbb{C}}(S_{H_{n}})\otimes L^{\infty}_{\mathbb{C}}(S_{O_{n-1}})\right)
\end{equation}
(notice that that $L^{\infty}_{\mathbb{C}}(S_{O_{0}})$ appears \textbf{twice} in the right hand side
of \eqref{df-A(H-tens-O)}) and the embeddings
\begin{equation}\label{df-O-H-emb0}
j_{H\otimes O;0} \colon F_{0}\in
L^{\infty}_{\mathbb{C}}(S_{H_{0}})\otimes L^{\infty}_{\mathbb{C}}(S_{O_{0}})
\to F_{0}\otimes 1_{\{0\}^c}\in\mathcal{A}_{H\otimes O}
\end{equation}
and, for $n\ge 1$,
\begin{equation}\label{df-O-H-emb1}
j_{H\otimes O;n} := j_{H_{n}}\otimes j_{O_{n-1}}
\colon L^{\infty}_{\mathbb{C}}(S_{H_{n}})\otimes L^{\infty}_{\mathbb{C}}(S_{O_{n-1}})
\to \mathcal{A}_{H\otimes O}
\end{equation}
where $j_{H_{n}}$ (resp. $j_{O_{n-1}}$) is the natural embedding of $L^{\infty}_{\mathbb{C}}(S_{H_{n}})$
(resp. $L^{\infty}_{\mathbb{C}}(S_{O_{n-1}})$) into $\mathcal{A}_{H\otimes O}$.
Now, for each $n\in\mathbb{N}$, denote
\begin{align}\label{df-hat-A(HO;[0,n])}
\hat{\mathcal{A}}_{H\otimes O;[0,n]} :=&\left(L^{\infty}_{\mathbb{C}}(S_{H_{0}})\otimes 1_{L^{\infty}_{\mathbb{C}}(S_{O_{0}})}\right)\otimes\\
&\bigotimes_{k=1}^{n}
\left(L^{\infty}_{\mathbb{C}}(S_{H_{k}})\otimes L^{\infty}_{\mathbb{C}}(S_{O_{k-1}})\right)
\otimes (1_{L^{\infty}_{\mathbb{C}}(S_{H_{n+1}})}\otimes L^{\infty}_{\mathbb{C}}(S_{O_{n}}))\notag
\end{align}
Clearly
\begin{align*}
\hat{\mathcal{A}}_{H\otimes O;[0,n]} \subset \mathcal{A}_{H\otimes O;[0,n]}:=& j_{H\otimes O;0}\left(L^{\infty}_{\mathbb{C}}(S_{H_{0}})\otimes L^{\infty}_{\mathbb{C}}(S_{O_{0}})\right) \vee\\
&\bigvee_{k\in [1,n]}j_{H\otimes O;n}\left(L^{\infty}_{\mathbb{C}}(S_{H_{k}})\otimes L^{\infty}_{\mathbb{C}}(S_{O_{k-1}})\right)
\end{align*}
and $(\hat{\mathcal{A}}_{H\otimes O;[0,n]})_{n\ge 0}$ is an increasing sequence of sub--algebras of
$\mathcal{A}_{H\otimes O;[0,n]}$.\\ Denote
\begin{equation}\label{df-hat-A)}
\hat{\mathcal{A}}_{H\otimes O} := \bigcup_{n\in\mathbb{N}}\hat{\mathcal{A}}_{H\otimes O;[0,n]}
\subset \mathcal{A}_{H\otimes O}
\end{equation}
and, for each $n\in\mathbb{N}$, denote
\begin{equation}\label{df-hat-P(H-otimes-O)}
\hat{P}_{H\otimes O;[0,n]}  := P_{H\otimes O}\Big|_{\hat{\mathcal{A}}_{H\otimes O;[0,n]}}
\quad;\quad  \hat{P}_{H\otimes O}  := P_{H\otimes O}\Big|_{\hat{\mathcal{A}}_{H\otimes O}}
\end{equation}
\noindent The following theorem shows that any hidden classical Markov process is the restriction,
to an appropriate sub--algebra, of a classical Markov process in the usual sense.
\begin{theorem}\label{th:(O,H)-is-a-QMC}{\rm
In the assumptions 1) and 2) at the beginning of this section
and for $p_{O_{0}}$ an \textbf{arbitrary probability measure on} $S_{O_{0}}$, let $P_{H\otimes O}$
be the state of the Markov process on $\mathcal{A}_{H\otimes O} $ characterized by the pair\\ $(p_{H_{0}}\otimes p_{O_{0}} \ , \ (P_{H_{n-1}\otimes O_{n-1};H_{n-1}})_{n\ge 1})$
(see Corollary \ref{thm:struct-Mark-st}).
Then
\begin{equation}\label{P(O)}
P_{H,O}
\equiv \hat{P}_{H\otimes O}
\end{equation}
where

$\bullet$ $P_{H,O}$ is defined in Corollary \ref{cor:joint-exps-H-Mark} and $\hat{P}_{H\otimes O}$ is defined by \eqref{df-hat-P(H-otimes-O)};

$\bullet$ $\equiv$ denotes stochastic equivalence (see Definition \ref{df:class-stoc-proc}).

\noindent In other words, in the notations \eqref{df-H-embeds-cl-alg} and \eqref{df-O-embeds-cl}, for
all $n\in\mathbb{N}$, $m\in \{0, \dots, n\}$, $f_{m}\in L^{\infty}_{\mathbb{C}}(S_{O_{m}})$ and
$g_{m}\in L^{\infty}_{\mathbb{C}}(S_{H_{m}})$:
\begin{align}\label{joint-probs-O-H-Mark-ker}
&P_{H,O}\left(\prod_{m=0}^{n}f_{m}(O_{m})g_{m}(H_{m})\right)\\
=&(p_{H_{0}}\otimes p_{O_{0}})\left((g_{0}\otimes 1_{O_{0}})P_{H_{0}\otimes O_{0};H_{0}}
\left((g_{1}\otimes f_{0})\cdots\right.\right.\notag\\
&\left.\left.\cdots\left( (g_{n-1}\otimes f_{n-2})
P_{H_{n-1}\otimes O_{n-1};H_{n-1}}\left( (g_{n}\otimes f_{n-1})
P_{H_{n}\otimes O_{n};H_{n}}(1_{H_{n+1}}\otimes f_{n})
\right)\right)\right)\right)\notag
\end{align}
Conversely, given two sequences of Markov operators $(P_{H_{n}})$ and $(B_{O_{n}})$ (not necessarily given by regular Markov kernels) and two arbitrary states $p_{H_{0}}$ on $L^{\infty}_{\mathbb{C}}(S_{H_{0}})$ and $p_{O_{0}}$ on $L^{\infty}_{\mathbb{C}}(S_{O;0})$, defining $P_{H_{n-1}\otimes O_{n-1};H_{n-1}}:= P_{H_{n-1}}\otimes B_{O_{n-1}}$ ($n\ge 1$),
the right hand side of \eqref{joint-probs-O-H-Mark-ker} \textbf{does not depend on} $p_{O_{0}}$ and defines the joint expectations of a unique state $P_{H\otimes O}$ on $\mathcal{A} _{H\otimes O} $ whose restriction to $\hat{\mathcal{A}} _{H\otimes O}$ coincides with $P_{H,O}$. \\
}\end{theorem}
\textbf{Proof}.
Because of \eqref{act-H-tens-O-id} and the Remark after Corollary \ref{thm:struct-Mark-st}
the pair\\ $(p_{H_{0}}\otimes p_{O_{0}} \ , \ (P_{H_{n-1}\otimes O_{n-1};H_{n-1}})_{n\ge 1}$
defines a unique Markov state on $\mathcal{A}_{H\otimes O} $.
Given the hidden Markov process, in terms of the regular Markov kernels for the $P_{H_{n}}$,
the joint expectations \eqref{joint-probs-O-H-Mark} become
\begin{equation}\label{joint-probs-O-H-Mark2}
P_{H,O}\left(\prod_{m=0}^{n}f_{m}(O_{m})g_{m}(H_{m})\right)
\end{equation}
$$
= p_{H_{0}}\left(E_{H_{0}}\left(B_{O,H,0}(f_{0})g_{0}\right)
P_{H_{0}}\left(B_{O,H,1}(f_{1})g_{1}\right)\left(\cdots
P_{H_{n-1}}\left(B_{O,H,n}(f_{n})g_{n}\right)\right)\right)
$$
$$
=\int_{\prod_{m=0}^{n}S_{H_{m}}}p_{H_{0}}(dh_{0})P_{H_{0}}(h_{0}, dh_{1})
\cdots P_{H_{n-2}}(h_{n-2}, dh_{n-1}) P_{H_{n-1}}(h_{n-1}, dh_{n})
$$
$$
\prod_{m=0}^{n}B_{O_{m},H_{m}}(f_{m})(h_{m})
\prod_{m=0}^{n}g_{m}(h_{m})
$$
$$
=\int_{\prod_{m=0}^{n}S_{H_{m}}}p_{H_{0}}(dh_{0})P_{H_{0}}(h_{0}, dh_{1})P_{H_{n-2}}(h_{n-2}, dh_{n-1})
\cdots P_{H_{n-1}}(h_{n-1}, dh_{n})
$$
$$
\prod_{m=0}^{n}\int_{S_{O;m}}f_{m}(o_{m}) B_{O_{m}}(h_{m}, do_{m})
\prod_{m=0}^{n}g_{m}(h_{m})
$$
$$
=\int_{\prod_{m=0}^{n}S_{H_{m}}}p_{H_{0}}(dh_{0})P_{H_{0}}(h_{0}, dh_{1})
\cdots P_{H_{n-2}}(h_{n-2}, dh_{n-1}) P_{H_{n-1}}(h_{n-1}, dh_{n})
$$
$$
\int_{\prod_{m=0}^{n}S_{O;m}}
B_{O_{0}}(h_{0}, do_{0})\cdots B_{O_{n}}(h_{n}, do_{n})
\prod_{m=0}^{n}f_{m}(o_{m})\prod_{m=0}^{n}g_{m}(h_{m})
$$
$$
=\int_{\prod_{m=0}^{n}S_{H_{m}}\times S_{O;m}}
g_{0}(h_{0})\prod_{m=1}^{n-1}g_{m}(h_{m})f_{m-1}(o_{m-1})f_{n}(o_{n})
$$
$$
p_{H_{0}}(dh_{0}) (P_{H_{0}}(h_{0}, dh_{1})B_{O_{0}}(h_{0}, do_{0}))
(P_{H_{0}}(h_{1}, dh_{2})B_{O_{1}}(h_{1}, do_{1}))
\cdots
$$
$$
\cdots(P_{H_{n-2}}(h_{n-2}, dh_{n-1})B_{O_{n-2}}(h_{n-2}, do_{n-2}))
(P_{H_{n-1}}(h_{n-1}, dh_{n})B_{O_{n-1}}(h_{n-1}, do_{n-1}))
$$
$$
B_{O_{n}}(h_{n}, do_{n})
$$
$$
=\int_{\prod_{m=0}^{n}S_{H_{m}}\times S_{O;m}}
p_{H_{0}}(dh_{0})g_{0}(h_{0})
(P_{H_{0}}(h_{0}, dh_{1})B_{O_{0}}(h_{0}, do_{0}))g_{1}(h_{1})f_{0}(o_{0})
$$
$$
(P_{H_{0}}(h_{1}, dh_{2})B_{O_{1}}(h_{1}, do_{1}))g_{2}(h_{2})f_{1}(o_{1})
\cdots
$$
$$
\cdots
(P_{H_{n-2}}(h_{n-2}, dh_{n-1})B_{O_{n-2}}(h_{n-2}, do_{n-2}))g_{n-1}(h_{n-1})f_{n-2}(o_{n-2})
$$
$$
(P_{H_{n-1}}(h_{n-1}, dh_{n})B_{O_{n-1}}(h_{n-1}, do_{n-1}))g_{n}(h_{n})f_{n-1}(o_{n-1})
$$
$$
B_{O_{n}}(h_{n}, do_{n})f_{n}(o_{n})
$$
Moreover, in agreement with \eqref{act-H-tens-O},
\begin{align*}
&\int_{S_{O;m-1}}B_{O_{n}}(h_{n}, do_{n})f_{n}(o_{n})
=(B_{O_{n}}f_{n})(o_{n})\\
=&\left(P_{H_{n}}\otimes B_{O_{n}}\right)(1_{H_{n+1}}\otimes f_{n})(h_{n}) =P_{H_{n}\otimes O_{n};H_{n}}(1_{H_{n+1}}\otimes f_{n})(h_{n})
\end{align*}
By assumption $p_{O_{0}}$ is a probability measure on $S_{O_{0}}$, hence
for any\\ $G_{0}\in L^{\infty}_{\mathbb{C}}(S_{H_{0}})$,
$$
\int_{S_{H_{0}}} p_{H_{0}}(dh_{0})G_{0}(h_{0})
= p_{H_{0}}(G_{0})
= (p_{H_{0}}\otimes p_{O_{0}})(G_{0}\otimes 1_{O_{0}})
$$
Therefore, in the notation \eqref {H-tens-O-ops},  \eqref{joint-probs-O-H-Mark2} becomes
\begin{equation}\label{joint-probs-O-H-Mark-ker3}
P_{H,O}\left(\prod_{m=0}^{n}f_{m}(O_{m})g_{m}(H_{m})\right)
\end{equation}
$$
=(p_{H_{0}}\otimes p_{O_{0}})\left((g_{0}\otimes 1_{O_{0}})P_{H_{0}\otimes O_{0};H_{0}}
\left((g_{1}\otimes f_{0})\cdots\right.\right.
$$
$$
\left.\left.
\cdots\left( (g_{n-1}\otimes f_{n-2})
P_{H_{n-1}\otimes O_{n-1};H_{n-1}}\left( (g_{n}\otimes f_{n-1})
P_{H_{n}\otimes O_{n};H_{n}}(1_{H_{n+1}}\otimes f_{n})
\right)\right)\right)\right)
$$
This proves \eqref{joint-probs-O-H-Mark-ker}. Conversely
given two sequences of Markov operators $(P_{H_{n}})$ and $(B_{O_{n}})$
(not necessarily given by regular Markov kernels) and two arbitrary states $p_{H_{0}}$ on $L^{\infty}_{\mathbb{C}}(S_{H_{0}})$ and $p_{O_{0}}$ on $L^{\infty}_{\mathbb{C}}(S_{O_{0}})$,
the joint expectations on the right hand side of \eqref{joint-probs-O-H-Mark-ker}
are projective because, putting $g_{n}=1_{H_{n}}$ and
$f_{n}=1_{O_{n}}$, \eqref{joint-probs-O-H-Mark-ker} becomes
\begin{align}\label{projectivity1}
&(p_{H_{0}}\otimes p_{O_{0}})\left((g_{0}\otimes 1_{O_{0}})P_{H_{0}\otimes O_{0};H_{0}}
\left((g_{1}\otimes f_{0})\cdots\right.\right.\\
&\left.\left.\cdots\left( (g_{n-1}\otimes f_{n-2})
P_{H_{n-1}\otimes O_{n-1};H_{n-1}}\left( 1_{O_n}\otimes f_{n-1})
P_{H_{n}\otimes O_{n};H_{n}}(1_{H_{n+1}}\otimes 1_{O_{n}})
\right)\right)\right)\right)\notag\\
=&(p_{H_{0}}\otimes p_{O_{0}})\left((g_{0}\otimes 1_{O_{0}})P_{H_{0}\otimes O_{0};H_{0}}
\left((g_{1}\otimes f_{0})\cdots\right.\right.\notag\\
&\left.\left. \cdots\left( (g_{n-1}\otimes f_{n-2})
P_{H_{n-1}\otimes O_{n-1};H_{n-1}}\left( (1_{O_{n}}\otimes f_{n-1})\right)\right)\right)\right)\notag
\end{align}
because $P_{H_{n+1}\otimes O_{n};H_{n}}$ is a Markov operator. Since \eqref{projectivity1}
is \eqref{joint-probs-O-H-Mark-ker} with $n$ replaced by $n-1$, projectivity holds.
Therefore \eqref{projectivity1} defines a unique Markov state $P$ on $\mathcal{A}_{H\otimes O} $.
We prove that these expectations do not depend on $p_{O_{0}}$ by induction.
The statement is true for $n=0$ because
$$
(p_{H_{0}}\otimes p_{O_{0}})(g_{0}\otimes 1_{O_{0}})
= p_{H_{0}}(g_{0})p_{O_{0}})(1_{O_{0}})
= p_{H_{0}}(g_{0})
$$
Supposing that the statement is true for $n-1$, one has
$$
(p_{H_{0}}\otimes p_{O_{0}})\left((g_{0}\otimes 1_{O_{0}})P_{H_{0}\otimes O_{0};H_{0}}
\left((g_{1}\otimes f_{0})\cdots\right.\right.
$$
$$
\left.\left.
\cdots\left( (g_{n-1}\otimes f_{n-2})
P_{H_{n-1}\otimes O_{n-1};H_{n-1}}\left( (g_{n}\otimes f_{n-1})
P_{H_{n}\otimes O_{n};H_{n}}(1_{H_{n+1}}\otimes f_{n})
\right)\right)\right)\right)
$$
$$
=(p_{H_{0}}\otimes p_{O_{0}})\left((g_{0}\otimes 1_{O_{0}})P_{H_{0}\otimes O_{0};H_{0}}
\left((g_{1}\otimes f_{0})\cdots\right.\right.
$$
$$
\left.\left.
\cdots\left( (g_{n-1}\otimes f_{n-2})
P_{H_{n-1}\otimes O_{n-1};H_{n-1}}\left( (g_{n}\otimes f_{n-1})
(P_{H_{n}}\otimes B_{O_{n}})(1_{H_{n+1}}\otimes f_{n})
\right)\right)\right)\right)
$$
$$
=(p_{H_{0}}\otimes p_{O_{0}})\left((g_{0}\otimes 1_{O_{0}})P_{H_{0}\otimes O_{0};H_{0}}
\left((g_{1}\otimes f_{0})\cdots\right.\right.
$$
\begin{equation}\label{indep-p(O;0)}
\left.\left.
\cdots\left( (g_{n-1}\otimes f_{n-2})
P_{H_{n-1}\otimes O_{n-1};H_{n-1}}\left((g_{n}\otimes B_{O_{n}}(f_{n})f_{n-1})\right)
\right)\right)\right)
\end{equation}
and, by the induction assumption, the right hand side of \eqref{indep-p(O;0)} does not
depend on $p_{O_{0}}$. Therefore, by induction, the right hand side of \eqref{joint-probs-O-H-Mark-ker}
does not depend on $p_{O_{0}}$.
Finally, if the operators $P_{H_{n}}$ and $B_{O_{n}}$ are given by regular Markov kernels,
as in \eqref{trans-ops-class-HMP} and \eqref{trans-ops-class-BOH}, the same arguments used
to prove the identity \eqref{joint-probs-O-H-Mark-ker} show that the joint expectations
\begin{equation}\label{joint-probs-O-H-Mark-ker4}
(p_{H_{0}}\otimes p_{O_{0}})\left((g_{0}\otimes 1_{O_{0}})P_{H_{0}\otimes O_{0};H_{0}}
\left((g_{1}\otimes f_{0})\cdots\right.\right.
\end{equation}
$$
\left.\left.
\cdots\left( (g_{n-1}\otimes f_{n-2})
P_{H_{n-1}\otimes O_{n-1};H_{n-1}}\left( (g_{n}\otimes f_{n-1})
P_{H_{n}\otimes O_{n};H_{n}}(1_{H_{n+1}}\otimes f_{n})
\right)\right)\right)\right)
$$
coincide with those of the unique hidden Markov process with Markov chain defined
by $(p_{H_{0}},(P_{H_{n}}))$ and emission operator sequence $(B_{O_{n}})$.
$\qquad\square$

\subsection{Classical backward hidden Markov process}

\begin{definition}\label{df:class-backw-HMP}{\rm
In the notations of Definition \ref{df:class-HP}, if the conditional independence
assumption \eqref{cond-exps-O|H} is replaced by the following:
$$
E_{H_{0},O_{0}, \dots,H_{n}, O_{n}}\left( j_{O_{n+1}}(f_{n+1})j_{H_{n+1}}(g_{n+1})\right)
$$
\begin{equation}\label{cond-exps-O(n+1)|Hn}
= E_{H_{n}}\left(j_{O_{n+1}}(f_{n+1})\right)E_{H_{n}}\left(j_{H_{n+1}}(g_{n+1})\right)
\end{equation}
the $(H,O)$--process is called a \textbf{backward hidden Markov process}
}\end{definition}
\textbf{Remark}.
In order to better understand the meaning of condition \eqref{cond-exps-O(n+1)|Hn} it is convenient to split it into two conditions: one is the usual Markov condition for the $HO$--process, namely,
$$
E_{H_{0},O_{0}, \dots,H_{n}, O_{n}}\left( j_{O_{n+1}}(f_{n+1})j_{H_{n+1}}(g_{n+1})\right)
$$
\begin{equation}\label{HO-proc-Markov}
= E_{H_{n}, O_{n}}\left(j_{O_{n+1}}(f_{n+1})j_{H_{n+1}}(g_{n+1})\right)
\end{equation}
The other one is an additional conditional independence assumption on the $HO$--process, namely,
$$
 E_{H_{n}, O_{n}}\left(j_{O_{n+1}}(f_{n+1})j_{H_{n+1}}(g_{n+1})\right)
$$
\begin{equation}\label{cond-indep-HO-proc}
= E_{H_{n}}\left(j_{O_{n+1}}(f_{n+1})\right)E_{H_{n}}\left(j_{H_{n+1}}(g_{n+1})\right)
\end{equation}
Notice that, with respect to \eqref{HO-proc-Markov}, condition \eqref{cond-indep-HO-proc}
\textbf{breaks the symmetry} between the $H$--process and the $O$--process. Moreover, in the
notations  \eqref{df-H-embeds-cl-alg}, \eqref{df-O-embeds-cl}, it suggests to introduce the Markov operators
\begin{equation}\label{df-O|H-Mark-op}
P_{O_{n+1},H_{n}} := j_{H_{n}}^{-1}\circ E_{H_{n}}\circ j_{O_{n+1}}
 \colon L^{\infty}_{\mathbb{C}}(S_{O_{n+1}}) \to L^{\infty}_{\mathbb{C}}(S_{H_{n}})
\end{equation}
\begin{equation}\label{df-H|H-Mark-op}
P_{H_{n+1},H_{n}} := j_{H_{n}}^{-1}\circ E_{H_{n}}\circ j_{H_{n+1}}
 \colon L^{\infty}_{\mathbb{C}}(S_{H_{n+1}}) \to L^{\infty}_{\mathbb{C}}(S_{H_{n}})
\end{equation}
Furthermore, while in condition \eqref{cond-exps-O|H} the observable process is conditioned on
the underlying one \textbf{at the same time}, in condition \eqref{cond-indep-HO-proc} the
observable process is conditioned on the underlying one \textbf{at the previous time}.\\
In the statement of the following theorem, to anticipate the link with the quantum extension,
we use the notations:
\begin{equation}\label{df-st-alg-n-cl2}
\mathcal{B}_{H_{n}} := L^{\infty}_{\mathbb{C}}(S_{H_{n}}) \quad;\quad
\mathcal{B}_{O_{n}} := L^{\infty}_{\mathbb{C}}(S_{O_{n}})
\end{equation}
\begin{theorem}\label{th:joint-exps-backw(O,H)}{\rm
Given a backward hidden Markov process in the sense of Definition \ref{df:class-backw-HMP},
the associated $H$--process is a Markov process.\\
Define the sequences of Markov operators $(P_{O_{n+1},H_{n}})_{n\ge 0})$ and
$(P_{H_{n+1},H_{n}})_{n\ge 0}$ respectively by \eqref{df-O|H-Mark-op} and \eqref{df-H|H-Mark-op}.
Denote, for $n\ge 1$,  $f_{n}\in \mathcal{B}(O_{n})$ and $g_{n}\in \mathcal{B}(H_{n})$
the transition expectation
$\mathcal{E}_{O_{n},H_{n};H_{n-1}}\colon \mathcal{B}(O_{n}) \otimes\mathcal{B}(H_{n})
\to \mathcal{B}(H_{n-1})$
(completely positive identity preserving map):
\begin{equation}\label{df-calE(On,Hn;Hn-1)}
\mathcal{E}_{O_{n},H_{n};H_{n-1}}(f_{n}\otimes g_{n})
:= P_{O_{n},H_{n-1}}(f_{n})P_{H_{n},H_{n-1}}(g_{n})
\end{equation}
and the initial distribution of the $H$--process
\begin{equation}\label{df-PH0}
P^{(0)}_{H_{0}} := P_{H,O}\big|_{\mathcal{A}_{H_{0}}}
\end{equation}
Then the probability distributions $P_{H,O}$ of the $(H,O)$--process are uniquely determined by the
triple $(P^{(0)}_{H_{0}},  P_{O_{0},H_{0}}, (\mathcal{E}_{O_{n},H_{n};H_{n-1}})_{n\ge 1})$ through
the identity
\begin{equation}\label{joint-exps-(O(n+1),Hn)-E3}
P_{H,O}\left(\prod_{m=0}^{n}j_{O_{m}}(f_{m})j_{H_{m}}(g_{m})\right)
\end{equation}
$$
=P_{H_{0}}\left(P_{O_{0},H_{0}}(f_{0})g_{0}\mathcal{E}_{O_{1},H_{1};H_{0}}(f_{1}\otimes \left(
g_{1}\mathcal{E}_{O_{2},H_{2};H_{1}}(f_{2}\otimes g_{2}) \cdots
\left( \mathcal{E}_{O_{n-2},H_{n-2};H_{n-3}}\left(
\right.\right.\right.\right.
$$
$$
\left.\left.\left. \left.\left.
(f_{n-2}\otimes \left(g_{n-2}\mathcal{E}_{O_{n-1},H_{n-1};H_{n-2}}\left(
f_{n-1}\otimes \left(g_{n-1}\mathcal{E}_{O_{n},H_{n};H_{n-1}}(f_{n}\otimes g_{n})\right)\right)
\right)\right)\right)\right)\right)\right)
$$
for $n\ge 1$, $m\in \{0, \dots, n\}$, $f_{m}\in \mathcal{B}(O_{m})$ and
$g_{m}\in \mathcal{B}(H_{m})$.\\
Conversely, let be given a triple
$(P_{H_{0}},  P_{O_{0},H_{0}}, (\mathcal{E}_{O_{n},H_{n};H_{n-1}})_{n\ge 1})$
such that $P_{H_{0}}$ is a (arbitrary) state on $\mathcal{B}_{H_{0}}$,
\begin{equation}\label{df-O0|H0-Mark-op-cl}
P_{O_{0},H_{0}} \colon
 \colon \mathcal{B}(O_{0}) \to \mathcal{B}(H_{0})
\end{equation}
is a (arbitrary) Markov operator and each
$$
\mathcal{E}_{O_{n},H_{n};H_{n-1}} \colon \mathcal{B}(O_{n})\otimes \mathcal{B}(H_{n})
\to \mathcal{B}(H_{n-1})
$$
is a transition expectation uniquely determined by the sequence of (arbitrary) Markov operators
\begin{equation}\label{df-P(O(n+1)H(n)}
P_{O_{n+1},H_{n}} \colon \mathcal{B}(O_{n+1}) \to \mathcal{B}(H_{n})
\end{equation}
\begin{equation}\label{df-P(H(n+1)H(n)}
P_{H_{n+1},H_{n}} \colon \mathcal{B}(H_{n+1}) \to \mathcal{B}(H_{n})
\end{equation}
through the identity \eqref{df-calE(On,Hn;Hn-1)}.
Then the right hand side of \eqref{joint-exps-(O(n+1),Hn)-E3} defines a unique state $P_{H,O}$
on $\mathcal{A}_{H,O}$ such that the quadruple \eqref{df-HO-proc-clas} is a backward hidden Markov
process.
}\end{theorem}
\textbf{Proof}.
\begin{equation}\label{joint-exps-(O(n+1),Hn)-E1}
P\left(\prod_{m=0}^{n}j_{O_{m}}(f_{m})j_{H_{m}}(g_{m})\right)
=P_{H_{0},O_{0}, \dots,H_{n},O_{n}}\left(\prod_{m=0}^{n}j_{O_{m}}(f_{m})j_{H_{m}}(g_{m})\right)
\end{equation}
$$
=P_{H_{0},O_{0}, \dots,H_{n-1},O_{n-1}}\left(\prod_{m=0}^{n-1}j_{O_{m}}(f_{m})j_{H_{m}}(g_{m})
E_{H_{0},O_{0},\dots,H_{n-1},O_{n-1}}\left(j_{O_{n}}(f_{n})j_{H_{n}}(g_{n})\right)\right)
$$
$$
\overset{\eqref{HO-proc-Markov}}{=}
P_{H_{0},O_{0}, \dots,H_{n-1},O_{n-1}}\left(\prod_{m=0}^{n-1}j_{O_{m}}(f_{m})j_{H_{m}}(g_{m})
E_{H_{n-1},O_{n-1}} \left(j_{O_{n}}(f_{n})j_{H_{n}}(g_{n})\right) \right)
$$
$$
\overset{\eqref{cond-indep-HO-proc}}{=}
P_{H_{0},O_{0}, \dots,H_{n-1},O_{n-1}}\left(\prod_{m=0}^{n-2}j_{O_{m}}(f_{m})j_{H_{m}}(g_{m})
\right.
$$
$$
\left. f_{m} (O_{n-1})
\left(g_{m}(H_{n-1})E_{H_{n-1}}\left(j_{O_{n}}(f_{n})\right)E_{H_{n-1}}\left(j_{H_{n}}(g_{n})\right)
\right)\right)
$$
$$
=P_{H_{0},O_{0}, \dots,H_{n-1},O_{n-1}}\left(\prod_{m=0}^{n-2}j_{O_{m}}(f_{m})j_{H_{m}}(g_{m})
\right.
$$
$$
\left. j_{O_{n-1}}(f_{n-1})
j_{H_{n-1}}\left(g_{n-1}P_{O_{n},H_{n-1}}(f_{n})P_{H_{n},H_{n-1}}(g_{n})\right) \right)
$$
Thus, with the notation \eqref{df-calE(On,Hn;Hn-1)}, for $n=1$ one obtains
$$
P\left(j_{O_{0}}(f_{0})j_{H_{0}}(g_{0})j_{O_{1}}(f_{1})j_{H_{1}}(g_{1})\right)
=P_{H_{0},O_{0}}\left(j_{O_{0}}(f_{0})j_{H_{0}}(g_{0}\mathcal{E}_{O_{1},H_{1};H_{0}}
(f_{1}\otimes g_{1}))\right)
$$
$$
=P\left(j_{H_{0}}(P_{O_{0},H_{0}}(f_{0})g_{0}\mathcal{E}_{O_{1},H_{1};H_{0}}(f_{1}\otimes g_{1}))\right)
$$
\begin{equation}\label{joint-exps-(O1,H1;H0)}
=P_{H_{0}}\left(P_{O_{0},H_{0}}(f_{0})g_{0}\mathcal{E}_{O_{1},H_{1};H_{0}}(f_{1}\otimes g_{1})\right)
\end{equation}
and, for arbitrary $n$,
\begin{equation}\label{joint-exps-(O(n+1),Hn)-E2}
P\left(\prod_{m=0}^{n}j_{O_{m}}(f_{m})j_{H_{m}}(g_{m})\right)
=P_{H_{0},O_{0}, \dots,H_{n-1},O_{n-1}}\left(   \right.
\end{equation}
$$
\left.\prod_{m=0}^{n-2}j_{O_{m}}(f_{m})j_{H_{m}}(g_{m})
 j_{O_{n-1}}(f_{n-1})
j_{H_{n-1}}\left(g_{n-1}\mathcal{E}_{O_{n},H_{n};H_{n-1}}(f_{n}\otimes g_{n})   \right) \right)
$$
$$
=P_{H_{0},O_{0}, \dots,H_{n-2},O_{n-2}}\left(\prod_{m=0}^{n-3}j_{O_{m}}(f_{m})j_{H_{m}}(g_{m})
\right.
$$
$$
\left. j_{O_{n-2}}(f_{n-2})j_{H_{n-2}}\left(g_{n-2}\mathcal{E}_{O_{n-1},H_{n-1};H_{n-2}}\left(
f_{n-1}\otimes \left(g_{n-1}\mathcal{E}_{O_{n},H_{n};H_{n-1}}(f_{n}\otimes g_{n})\right)\right)
\right)\right)
$$
$$
=P_{H_{0},O_{0}, \dots,H_{n-3},O_{n-3}}\left(\prod_{m=0}^{n-3}j_{O_{m}}(f_{m})j_{H_{m}}(g_{m})
j_{H_{n-3}}\left( \mathcal{E}_{O_{n-2},H_{n-2};H_{n-3}}\left(
\right.\right.\right.
$$
$$
\left.\left.\left. j_{H_{n-3}}\left( \mathcal{E}_{O_{n-2},H_{n-2};H_{n-3}}\left(
(f_{n-2}\otimes \left(g_{n-2}\mathcal{E}_{O_{n-1},H_{n-1};H_{n-2}}\left(
\right.\right.\right.\right.\right.\right.\right.
$$
$$
\left.\left.\left.
f_{n-1}\otimes \left(g_{n-1}\mathcal{E}_{O_{n},H_{n};H_{n-1}}(f_{n}\otimes g_{n})\right)\right)
\right)\right)
$$
Suppose by  induction that \eqref{joint-exps-(O(n+1),Hn)-E3} holds for $n\in\mathbb{N}$.
Then
\begin{equation}\label{joint-exps-(O(n+1),Hn)-E4}
P\left(\prod_{m=0}^{n+1}j_{O_{m}}(f_{m})j_{H_{m}}(g_{m})\right)
\end{equation}
$$
=P\left(\prod_{m=0}^{n}j_{O_{m}}(f_{m})j_{H_{m}}(g_{m})
j_{O_{n+1}}(f_{n+1})j_{H_{n+1}}(g_{n+1})\right)
$$
$$
=P\left(\prod_{m=0}^{n}j_{O_{m}}(f_{m})j_{H_{m}}(g_{m})
E_{H_{0},O_{0},\dots,H_{n},O_{n}}\left(
j_{O_{n+1}}(f_{n+1})j_{H_{n+1}}(g_{n+1})\right)\right)
$$
$$
=P\left(\prod_{m=0}^{n}j_{O_{m}}(f_{m})j_{H_{m}}(g_{m})
E_{O_{n+1},H_{n}}\left(j_{O_{n+1}}(f_{n+1})\right)
E_{H_{n+1},H_{n}}\left(j_{H_{n+1}}(g_{n+1})\right)\right)
$$
$$
=P\left(\prod_{m=0}^{n-1}j_{O_{m}}(f_{m})j_{H_{m}}(g_{m})\right.
$$
$$
\left. j_{O_{n}}(f_{n})j_{H_{n}}\left(
g_{n}P_{H_{n+1},O_{n+1};H_{n}}(f_{n+1})
P_{H_{n+1},H_{n}}\left(j_{H_{n+1}}(g_{n+1})\right)\right)\right)
$$
$$
=P\left(\prod_{m=0}^{n}j_{O_{m}}(\widehat{f}_{m})j_{H_{m}}(\widehat{g}_{m})\right)
$$
where
$$
\widehat{f}_{m} := f_{m}
\quad;\quad \widehat{g}_{m} :=\begin{cases}
                                g_{m}, & \mbox{if } 0\le m\le n-1 \\
g_{n}P_{H_{n+1},O_{n+1};H_{n}}(f_{n+1})
P_{H_{n+1},H_{n}}\left(j_{H_{n+1}}(g_{n+1})\right), & \mbox{if }  m=n
                              \end{cases}
$$
Therefore by induction \eqref{joint-exps-(O(n+1),Hn)-E3} holds.
$\qquad\square$

\subsubsection{Time--consecutive hidden process}\label{sec:Time-cons-HP}

\begin{theorem}\label{Comparison-2-dfs-HMP}{\rm
Condition \eqref{cond-exps-O(n+1)|Hn} is equivalent to the fact that:\\
(i) the process $(H_{n})$ is a Markov process;\\
(ii) the following conditional independence condition holds:
\begin{equation}\label{cond-exp-HMP}
E_{H_{0}, \dots,H_{n}}\left(\prod_{m=0}^{n}j_{O_{m}}(f_{m})\right)
=  \prod_{m=0}^{n}E_{H_{m-1}}\left(j_{O_{m}}(f_{m})\right)
\end{equation}
for all $n\in\mathbb{N}$ and $f_{m}\in\mathcal{B}(O_{m})$ ($m\in\{0,\dots,n\}$)
and with the convention that
\begin{equation}\label{conv-E(H-1)}
E_{H_{-1}}(j_{O_{0}}(f_{0})):=E_{H_{0}}(j_{O_{0}}(f_{0}))
\end{equation}
}\end{theorem}
\textbf{Proof}.
From Theorem \ref{th:joint-exps-backw(O,H)}, we know that condition \eqref{cond-exps-O(n+1)|Hn} is
equivalent to the fact that the joint distributions of the $(H,O)$--process are given by
\eqref{joint-exps-(O(n+1),Hn)-E3}, i.e. to the fact that, for any $n\in\mathbb{N}$,
$g_{m}\in\mathcal{B}(H_{m})$ and $f_{m}\in\mathcal{B}(O_{m})$ ($m\in\{0,\dots,n\}$), one has
\begin{equation}\label{joint-exps-(O(n+1),Hn)-E3a}
P_{H,O}\left(\prod_{m=0}^{n}j_{O_{m}}(f_{m})j_{H_{m}}(g_{m})\right)
\end{equation}
$$
\overset{\eqref{joint-exps-(O(n+1),Hn)-E3}}{=}
P_{H_{0}}\left(P_{O_{0},H_{0}}(f_{0})g_{0}\mathcal{E}_{O_{1},H_{1};H_{0}}(f_{1}\otimes \left(
g_{1}\mathcal{E}_{O_{2},H_{2};H_{1}}(f_{2}\otimes g_{2}) \cdots
\left( \mathcal{E}_{O_{n-2},H_{n-2};H_{n-3}}\left(
\right.\right.\right.\right.
$$
$$
\left.\left.\left. \left.\left.
(f_{n-2}\otimes \left(g_{n-2}\mathcal{E}_{O_{n-1},H_{n-1};H_{n-2}}\left(
f_{n-1}\otimes \left(g_{n-1}\mathcal{E}_{O_{n},H_{n};H_{n-1}}(f_{n}\otimes g_{n})\right)\right)
\right)\right)\right)\right)\right)\right)
$$
$$
\overset{\eqref{df-calE(On,Hn;Hn-1)}}{=}
P_{H_{0}}\left(P_{O_{0},H_{0}}(f_{0})g_{0}(P_{O_{1},H_{0}}(f_{1}) \left(
g_{1}P_{O_{2},H_{1}}(f_{2})P_{H_{2},H_{1}}(g_{2})  \cdots
\left(P_{H_{n-2};H_{n-3}}\left(
\right.\right.\right.\right.
$$
$$
\left.\left.\left. \left.\left.
(P_{O_{n-2},H_{n-3}}(f_{n-2}) \left(g_{n-2}P_{H_{n-1};H_{n-2}}\left(
P_{O_{n-1},H_{n-2}}(f_{n-1})\left(g_{n-1}P_{O_{n},H_{n-1}}(f_{n})P_{H_{n},H_{n-1}}(g_{n})
\right)\right)\right)\right)\right)\right)\right)\right)
$$
Since we know that $(H_{n})$ is a Markov process, this is equal to
$$
P_{H_{0}, \dots,H_{n}}\left(j_{H_{0}}(g_{0})j_{H_{1}}(g_{1})\cdots j_{H_{n}}(g_{n})
j_{H_{0}}(P_{O_{0},H_{0}}(f_{0}))j_{H_{0}}(P_{O_{1},H_{0}}(f_{1})) \cdots  j_{H_{n-1}}(P_{O_{n},H_{n-1}}(f_{n}))\right)
$$
$$
= P_{H_{0}, \dots,H_{n}}\left(j_{H_{0}}(g_{0})j_{H_{1}}(g_{1})\cdots j_{H_{n}}(g_{n})
j_{H_{0}}(P_{O_{0},H_{0}}(f_{0}))
\prod_{m=1}^{n}j_{H_{m-1}}(P_{O_{m},H_{m-1}}(f_{m}))\right)
$$
\begin{equation}\label{cond-exp-HMP1}
= P_{H_{0}, \dots,H_{n}}\left(\prod_{m=0}^{n}j_{H_{m}}(g_{m})
\prod_{m=0}^{n}E_{H_{m-1}}(j_{O_{m}}(f_{m}))\right)
\end{equation}
with the convention \eqref{conv-E(H-1)}.
On the other hand
$$
P_{H,O}\left(\prod_{m=0}^{n}j_{O_{m}}(f_{m})j_{H_{m}}(g_{m})\right)
= P_{H,O}\left(\prod_{m=0}^{n}j_{H_{m}}(g_{m})
\prod_{m=0}^{n}j_{O_{m}}(f_{m})\right)
$$
\begin{equation}\label{cond-exp-HMP2}
= P_{H_{0}, \dots,H_{n}}\left(\prod_{m=0}^{n}j_{H_{m}}(g_{m})
E_{H_{0}, \dots,H_{n}}(\prod_{m=0}^{n}j_{O_{m}}(f_{m}))\right)
\end{equation}
and, since the $g_{m}$ ($m\in \{0,1,\dots,n\}$) are arbitrary, the identities \eqref{cond-exp-HMP1}
\eqref{cond-exp-HMP2} are equivalent to \eqref{cond-exp-HMP}.
This shows that, under the assumption that the process $(H_{n})$ is a Markov process,
condition \eqref{cond-exps-O(n+1)|Hn} is equivalent to \eqref{cond-exp-HMP}.
$\qquad\square$\\

\noindent\textbf{Remark}.
The fact that, contrarily to what happens for condition \eqref{cond-exps-O(n+1)|Hn},
condition \eqref{cond-exp-HMP} does not imply that the process $(H_{n})$ is Markov,
suggests the following definition.
\begin{definition}\label{df:time-cons-HMP}{\rm
In the notations of Definition \ref{df:class-HP}, the classical stochastic
$(H,O)$--process
\begin{equation}\label{df-HO-proc-TC}
\left( \mathcal{A}_{H,O}, P_{H,O}, (j_{H_{n}}\otimes j_{O_{n}}),
(\mathcal{B}_{H_{n}}\otimes \mathcal{B}_{O_{n}})  \right)
\end{equation}
is called a \textbf{time--consecutive hidden process} if condition \eqref{cond-exp-HMP} is
satisfied. If the hidden process $(H_{n})$ is Markov, one speaks of \textbf{time--consecutive
hidden Markov process} (as opposed to the same time HMP defined by condition \eqref{cond-exps-O|H}).
}\end{definition}

\subsection{Classical generalized hidden processes}

In this section we introduce a natural extension of the classical hidden processes, discussed in the previous section and in turn generalizing hidden Markov processes. This extension will help
to better understand the quantum extension of HMP introduced in section \ref{sec:QHMP}.\\
Theorem \ref{th:joint-exps-same-time(O,H)} shows that classical hidden processes are
characterized by the property that their joint probabilities have the form given by the right
hand side of \eqref{joint-exps-(O,H)-B-expl}. Now notice that, using the homomorphism property
of $j_{H_m}$, \eqref{joint-exps-(O,H)-B-expl} can be written in the form
\begin{equation}\label{joint-exps-(O,H)-B-expl2}
P_{H,O}\left(\prod_{m=0}^{n}f_{m}(O_{m})g_{m}(H_{m})\right)
=P_{H}\left(\prod_{m=0}^{n}j_{H_{m}}(B_{O_{m},H_{m}}(f_{m})g_{m})\right)
\end{equation}
where $n\in\mathbb{N}$, $m\in \{0, \dots, n\}$, $f_{m}\in L^{\infty}_{\mathbb{C}}(S_{O_{m}})$ and
$g_{m}\in L^{\infty}_{\mathbb{C}}(S_{H_{m}})$.\\
Therefore, introducing the operator
\begin{equation}\label{non-pos-map}
\widehat{B}_{{O,H,n}}(f_{n}\otimes g_{n})
:= B_{O_{O,H,n}}(f_{n})g_{n}
\end{equation}
which is a Markov operator (transition expectation, see Definition \ref{df-trans-exp})
$$
\widehat{B}_{O,H,n} \colon
L^{\infty}_{\mathbb{C}}(S_{O_{n}})\otimes L^{\infty}_{\mathbb{C}}(S_{H_{n}}) \to
L^{\infty}_{\mathbb{C}}(S_{H_{n}})
\equiv 1_{L^{\infty}_{\mathbb{C}}(S_{O_{n}})}\otimes  L^{\infty}_{\mathbb{C}}(S_{H_{n}})
$$
\eqref{joint-exps-(O,H)-B-expl2} becomes
\begin{equation}\label{joint-exps-(O,H)-B-expl3}
P_{H,O}\left(\prod_{m=0}^{n}f_{m}(O_{m})g_{m}(H_{m})\right)
=P_{H}\left(\prod_{m=0}^{n}j_{H_{m}}(\widehat{B}_{{O,H,m}}(f_{m}\otimes g_{m}))\right)
\end{equation}
We will see in section \ref{sec:QHMP} that the quantum extension of HMP simply replaces the
algebras $L^{\infty}_{\mathbb{C}}(S_{O_{n}})$ and $L^{\infty}_{\mathbb{C}}(S_{H_{n}})$ by
arbitrary non--commutative $*$--algebras and allows the operator $\widehat{B}_{{O,H,n}}$ to be an
arbitrary Markov operator. In the definition below, we introduce a further generalization of HMP
which might, at first sight, seem to be artificial, but it is necessary to include the new class of hidden (Markov) processes that arise as diagonal restrictions of quantum Markov chains (see Section
\ref{sec:HP-Diag-restr-De-pres-QMC} below).
\eject
\begin{definition}\label{df:class-gen-hidd-proc}{\rm
A \textbf{classical generalized hidden process} is defined by a pair of stochastic processes
$(O,H)\equiv (O_{n},H_{n})_{n\in \mathbb{N}}$ on the same probability space characterized by
the joint expectations:
\begin{equation}\label{joint-exps-(O,H)-B-expl-gen}
P_{H,O}\left(\prod_{m=0}^{n}f_{m}(O_{m})g_{m}(H_{m})\right)
=P_{H}\left(\prod_{m=0}^{n}j_{H_{h_{n;m}}}(\widehat{B}_{{O,H,m}}(f_{m}\otimes g_{m}))\right)
\end{equation}
where the operators $\widehat{B}_{{O,H,m}}$ are general Markov operators, i.e. not necessarily
of the form \eqref{non-pos-map}, and for each $n\in\mathbb{N}$,
$$
h_{n; \, \cdot \, } \colon m\in\{0,1,\dots, n\} \to h_{n;m}\in\{0,1,\dots, n\}
$$
is a map (not necessarily $1$--to--$1$).
}\end{definition}
\noindent\textbf{Remark}.
When each $h_{n; \, \cdot \, }$ is the \textbf{identity function} ($h_{n;m}=m \, , \, \forall m$),
and $\widehat{B}_{{O,H,m}}$ has the special form \eqref{non-pos-map}, the joint expectations
\eqref{joint-exps-(O,H)-B-expl-gen} take the form \eqref{joint-exps-(O,H)-B-expl2}
which characterize the hidden processes (see Definition \ref{df:class-HP}).\\

\noindent If $k_{n;m}:=m-1$ for $m\ge 1$ and $\widehat{B}_{{O,H,m}}$ has the special
form \eqref{non-pos-map}, the joint expectations \eqref{joint-exps-(O,H)-B-expl-gen} take the
form \eqref{joint-exps-(O(n+1),Hn)-E3} which characterize the backward hidden processes
(see Definition \ref{df:class-backw-HMP}).\\

\noindent When $h_{n; \, \cdot \, }$ is the \textbf{constant function}
$h_{n;m}=n-1 \, , \, \forall m$, and $\widehat{B}_{{O,H,m}}$ is given by \eqref{non-pos-map},
the joint expectations \eqref{joint-exps-(O,H)-B-expl-gen} extend those that one finds by restricting an $e$--diagonal quantum Markov chain to the $e$--diagonal sub--algebra (see sections
\ref{sec:QMC-HMP-max-obs}, \ref{sec:diag-restr-QHMP} and \ref{sec:HP-Diag-restr-De-pres-QMC}
for more details).

\section{Algebraic (classical or quantum) stochastic process}\label{sec:Alg-clas-q-stoc-proc}

\noindent The \textbf{transition from classical to algebraic} stochastic process(which include
both classical and quantum) is achieved replacing, in \eqref{alg-class-stoch-proc},
$L^{\infty}_{\mathbb{C}}(\Omega,\mathcal{F}, P)$ and the $L^{\infty}_{\mathbb{C}}(S_{n})$ by
arbitrary, not necessarily commutative, $*$--algebras.
\begin{definition}\label{df:alg-stoch-proc}{\rm
An \textbf{algebraic stochastic process} indexed by a set $T$ is a quadruple
$$
(\mathcal{A},\varphi, (\mathcal{B}_{n})_{n\in T}, (j_{n})_{n\in T})
$$
where:\\
-- $\mathcal{A}$, called the \textbf{sample algebra} of the process, is a $*$--algebra.\\
-- Each $\mathcal{B}_{n}$, called the \textbf{state algebra} at $n\in T$, is a $*$--algebra.\\
-- For each $n\in T$,
\begin{equation}\label{df-jn-tens2}
j_{n} : b_{n}\in \mathcal{B}_{n} \to j_{n}(b_{n})\in \mathcal{A}
\end{equation}
is a $*$--homomorphism.
The pair $(\mathcal{A},\varphi)$ is called an \textbf{algebraic probability space} (classical if
$\mathcal{A}$ is commutative, quantum if it is non--commutative).
}\end{definition}
\textbf{Remark}.
Each $*$--algebra $\mathcal{B}_{n}$ has a set (possibly infinite) of algebraically independent
hermitean generators
$$
B_{n} := \{b_{n;1}, b_{n;2}, \dots, b_{n;d_{n}} \} \quad,\quad n\in T \ , \
d_{n}\in \mathbb{N}\cup\{+\infty\}
$$
Therefore, the assignment of the family of $*$--homomorphism $(j_{n})_{n\in T}$ is equivalent
to give the set
\begin{equation}\label{df:op-proc}
 \{X_{n;h} := j_{n}(b_{n;h}) \colon  h\in \{1,\dots, d_{n}\} \ , \ n\in T\} \subseteq \mathcal{A}
\end{equation}
of hermitean operators (which algebraically generates $j_{n}(B_{n})$). It is known that, if $(\mathcal{A} , \varphi)$ is an algebraic probability space, then any hermitean element of
$\mathcal{A}$ can be identified, up to moment equivalence, to a classical real valued random variable. The family \eqref{df:op-proc} is called an \textbf{operator stochastic process}.
In many cases in quantum probability, one deals with operator stochastic process.
However the operator stochastic processes \textbf{only include processes with finite moments
of all orders}, while Definition \ref{df:alg-stoch-proc} includes all stochastic processes.
Finally it is important to notice that, in quantum probability, one often uses sets of generators
that contain not only hermitean elements (but it is always possible to go back to the hermitean case).
\begin{definition}\label{df:tens-KR-alg-stoch-proc}{\rm
An algebraic stochastic process in \textbf{tensor Kolmogorov representation} is given by a quadruple
of the form
\begin{equation}\label{tens-KR-alg-stoch-proc}
\left(\mathcal{A}:=\bigotimes_{T}\mathcal{B}_{n},\varphi , (\mathcal{B}_{n})_{n\in T},
(j_{n})_{n\in T}\right)
\end{equation}
with embeddings given by
\begin{equation}\label{df-tens-emb-An}
j_{n} \colon b_{n} \in \mathcal{B}_{n} \to b_{n}\otimes 1_{T\setminus \{n\}}
\in j_{n}(\mathcal{B}_{n}) := \mathcal{A}_{n}
\end{equation}
where $1_{T\setminus \{n\}}$ is the identity of $\bigotimes_{m\in T\setminus \{n\}}\mathcal{B}_{m}$.
}\end{definition}
\noindent  In the following we assume that $T$ is an \textbf{at most countable set}
and we only consider algebraic stochastic processes in \textbf{tensor Kolmogorov representation}
(from now on simply called \textbf{algebraic stochastic processes} unless a specification
is appropriate).
For these processes we define the associated local algebras.
\begin{definition}\label{df:loc-alg}{\rm
For any algebraic stochastic process of the form \eqref{tens-KR-alg-stoch-proc}, and for any
$F\subset_{fin} T$, define the \textbf{local algebra on} $F$ by
\begin{equation}\label{df-loc-alg-gen}
\mathcal{A}_{F} := \bigvee_{n\in F}j_{n}(\mathcal{B}_{n})
\left(\hbox{ $*$--algebra generated by } \{j_{n}(\mathcal{B}_{n}) \colon n\in F\}\right)
\subseteq \mathcal{A}
\end{equation}
\textbf{Remark}.
Contrarily to what happens in the classical case, not every quantum stochastic process
admits a tensor Kolmogorov representation (in quantum physics there are many examples of
such processes), but many of the following notions can be extended to the whole class of
quantum stochastic process.\\
The family $\{\mathcal{A}_{F} \colon F\in\mathcal{F}_{fin}(T)\}$ will be called the family of
\textbf{local algebras}. The restriction of $\varphi$ on $\mathcal{A}_{F}$ will be denoted
$\varphi_{F}$.
\begin{equation}\label{identificat-loc-alg}
\mathcal{A}_{F} := \bigvee_{n\in F}j_{n}(\mathcal{B}_{n}) \equiv \bigotimes_{n\in F}\mathcal{B}_{n}
\end{equation}
where $\equiv $ means that the two notations are used indifferently.
}\end{definition}
\noindent The family of states $\{\varphi_{F} \colon F\in\mathcal{F}_{fin}(T)\}$ satisfies
the analogue of the Kolmogorov compatibility conditions, i.e.:\\
1) Each $\varphi_{F}$ is a state on $\mathcal{A}_{F}$.\\
2) For any $F\subset G\in\mathcal{F}_{fin}(T)$,
\begin{equation}\label{alg-proj-cond-st}
\varphi_{G}\big|_{\mathcal{A}_{F}} = \varphi_{F}
\end{equation}
\begin{definition}{\rm
Let be given, for each $F\in\mathcal{F}_{fin}(T)$, a state $\varphi_{F}$ on $\mathcal{A}_{F}$.
If the family $(\varphi_{F})_{F\in\mathcal{F}_{fin}(T)}$ satisfies conditions 1) and 2) above, it is
called \textbf{projective}.
If, for each $a\in\mathcal{A}$, there exists $F_0(a)$ such that, for all $F\supset F_0(a)$,
$\varphi_{F}(a)=\varphi_{F_0}(a)$, one says that the limit $\lim_{F\uparrow T}\varphi_{F}$
exists \textbf{in the strongly finite sense on} $\mathcal{A}$. $\mathcal{F}_{fin}(T)$ is an
increasing net for the partial order induced by inclusion and $F\uparrow T$ is understood in
the sense of this partial order.
}\end{definition}
The following is a weak algebraic formulation of Kolmogorov compatibility theorem.
\begin{theorem}\label{proj-fam-equiv-st}{\rm
Let the local algebras $\mathcal{A}_{F}$ be given by \eqref{df-loc-alg-gen} and let, for each
$F\subset_{fin} T$ be given a state $\varphi_{F}$ on $\mathcal{A}_{F}$ so that
 $(\varphi_{F})_{F\in\mathcal{F}_{fin}(T)}$ is a projective family of states. Then the limit
\begin{equation}\label{proj-impl-st}
\lim_{F\uparrow T}\varphi_{F} =: \varphi
\end{equation}
exists in the strongly finite sense on on $\mathcal{A}$.
}\end{theorem}
\textbf{Proof}.
Since the tensor product in \eqref{tens-KR-alg-stoch-proc} is algebraic, every element of
$\mathcal{A}$ belongs to some $\mathcal{A}_{F}$. Therefore, by projectivity, for each
$a\in\mathcal{A}$, there exists $F_0(a)$ such that, for all $F\supset F_0(a)$,
$\varphi_{F}(a)=\varphi_{F_0}(a)$. This means that the limit \eqref{proj-impl-st}
exists in the strongly finite sense on $\mathcal{A}$. $\qquad\square$\\


\subsection{Quantum Markov chains (homogeneous, backward)}\label{sec:QMC}

In this section we briefly review, mainly to fix the notations, the construction of quantum
Markov chains \cite{Ac74-Camerino}, \cite{Ac74-FAA}.
For simplicity we limit our considerations to the backward homogeneous case: this contains all the main ideas and once this is understood, the many possible variants will present no difficulties.\\
Let be $H$ a separable Hilbert space and $\mathcal{B}=\mathcal{B}(H)$ the algebra of bounded
linear operators on $H$. Denote
$$
\mathcal{A}:=\bigotimes_{\mathbb{N}}\mathcal{B}
$$
where, $\bigotimes$ denotes the \textbf{algebraic tensor product} (in this paper we do not discuss
topological aspects).
For each $n\in \mathbb{N}$, define the natural embedding ($*$--homomorphism) called the
$n$--th \textbf{tensor embedding} in the following way:
\begin{equation}\label{df-jn-tens}
j_n\colon b \in \mathcal{B} \ \hookrightarrow \
j_n(b) := \left(\bigotimes_{\{n\}^c}1_{\mathcal{B}}\right) \otimes b
\qquad,\qquad \forall b\in \mathcal{B}
\end{equation}
\begin{definition}\label{df-trans-exp}{\rm
A linear map $\mathcal{E}$ from $\mathcal{B}\otimes \mathcal{B}$ to $\mathcal{B}$ is called
a \textbf{transition expectation} if it is completely positive and identity preserving (i.e. a 
Markov operator).
}\end{definition}
We will study transition expectations $\mathcal{E}$ of the form
\begin{equation}\label{df-cal-E-Tr-CDA}
\mathcal{E}(x)=\overline{\hbox{Tr}}_2\left(\sum_{r\in D_{\mathcal{E}}}K^*_r\ x\ K_r\right)
\quad,\quad x\in \mathcal{B}\otimes \mathcal{B}
\end{equation}
where

$\bullet$ the partial trace with respect to the second factor $\overline{\hbox{Tr}}_{2}$ is the operator valued weight, in the sense of Haagerup,
defined by
\begin{equation}\label{df-partial-Tr}
\overline{\hbox{Tr}}_{2}(a\otimes b) := a\hbox{Tr}(b)   \quad,\quad a,b\in \mathcal{B}
\end{equation}

$\bullet$ $D_{\mathcal{E}}$ is a sub--set of $\mathbb{N}$ (finite if $H$ is finite dimensional) and, for $j\in D_{\mathcal{E}}$, $K_j\in \mathcal{B}\otimes \mathcal{B}$.\\

\noindent\textbf{Remark}.
Every $K_{r}$ in \eqref{df-cal-E-Tr-CDA} can be written in two forms
\begin{equation}\label{Kr-in-matr-uns}
K_{r}=\ \sum_{i,j\in D}e_{ij}\otimes\ K_{r;ij}
=\ \sum_{i,j\in D}K'_{r;ij}\otimes\  e_{ij}
\quad,\quad r\in D_{\mathcal{E}}
\end{equation}
where the $K_{r;ij}, K'_{r;ij} \in \mathcal{B}$. Therefore, for $a, b\in \mathcal{B}$,
$$
\mathcal{E}(a\otimes b)
=\overline{\hbox{Tr}}_2\left(\sum_{r\in D_{\mathcal{E}}}\left(
\sum_{i,j\in D}(e_{ij}\otimes\ K_{r;ij})^*\ (a\otimes b) \
\sum_{i',j'\in D}e_{i'j'}\otimes\ K_{r;i'j'}\right)\right)
$$
$$
=\sum_{r\in D_{\mathcal{E}}}\sum_{i,j,i',j'\in D}\overline{\hbox{Tr}}_2\left(
e_{ji}ae_{i'j'}\otimes\ (K_{r;ij}^*\ b \ K_{r;i'j'}))\right)
$$
\begin{equation}\label{df-cal-E-spec1}
=\sum_{r\in D_{\mathcal{E}}}\sum_{i,j,i',j'\in D}
e_{ji}ae_{i'j'} \ \hbox{Tr}\left(K_{r;ij}^*\ b \ K_{r;i'j'}\right)
\end{equation}
\noindent The structure of the operators $K_{r;ij}$ in \eqref{Kr-in-matr-uns} is given by the
following lemma which also provides a simple rule to construct conditional density amplitudes.
\begin{lemma}\label{|K(jk)|2-stoch-matr}{\rm
A family of operators of the form \eqref{Kr-in-matr-uns}, for some ONB $e\equiv (e_{j})_{j\in D}$,
is such that the operator $\mathcal{E}$ defined by \eqref{df-cal-E-Tr-CDA} is a transition
expectation if and only if
\begin{equation}\label{conds-on-K(r;ij')}
\sum_{r\in D_{\mathcal{E}}}\sum_{i\in D} \overline{\hbox{Tr}}(K_{r;ij}^*K_{r;ij'})
=\delta_{j,j'}
\quad,\quad \forall j,j'\in D
\end{equation}
}\end{lemma}
\textbf{Proof}.
The complete positivity of $\mathcal{E}$ follows by inspection of the right hand side
of \eqref{df-cal-E-Tr-CDA}.
The conditon $\mathcal{E}(1\otimes 1)=1$ is equivalent to
\begin{align*}
1_{\mathcal{B}}=&\mathcal{E}(1_{\mathcal{B}}\otimes 1_{\mathcal{B}})=\overline{\hbox{Tr}}_2(\sum_{r\in D_{\mathcal{E}}}K^*_rK_r)\\
=&\sum_{r\in D_{\mathcal{E}}}\overline{\hbox{Tr}}_2((
\sum_{i,j\in D}e_{ji}\otimes\ K_{r;ij}^*)
(\sum_{i',j'\in D}e_{i'j'}\otimes\ K_{r;i'j'}))\\
=&\sum_{r\in D_{\mathcal{E}}}\sum_{i,j,j'\in D}\overline{\hbox{Tr}}_2(e_{jj'}\otimes\
K_{r;ij}^*K_{r;ij'})
=\sum_{j,j'\in D}e_{jj'}\sum_{r\in D_{\mathcal{E}}}\sum_{i\in D} \overline{\hbox{Tr}}(K_{r;ij}^*K_{r;ij'})
\end{align*}
Clearly, this is equivalent to
$$\sum_{r\in D_{\mathcal{E}}}\sum_{i\in D} \overline{\hbox{Tr}}(K_{r;ij}^*K_{r;ij'})
=\delta_{j,j'}\,,\quad \forall j,j'\in D
$$
which is \eqref{conds-on-K(r;ij')}.
$\qquad\square$

\subsubsection{Conditional density amplitudes and operator valued isometries}

A special class of transition expectations of the form \eqref{df-cal-E-Tr-CDA}, important in
many applications, is obtained when
$$
|D_{\mathcal{E}}| = 1
$$
In this case, letting
\begin{equation}\label{df-cond-dens-ampl}
K=\ \sum_{i,j\in D}e_{ij}\otimes\ K_{ij} \ ;\quad
\overline{\hbox{Tr}}_{2}(K^*K) := 1_{\mathcal{B}}
\end{equation}
the index $r$ is absent in \eqref{df-cal-E-Tr-CDA} which becomes
\begin{equation}\label{df-cal-E-Tr-CDA-1K}
\mathcal{E}(x)=\overline{\hbox{Tr}}_2\left(\sum_{r\in D_{\mathcal{E}}}K^*_r\ x\ K_r\right)
\quad,\quad x\in \mathcal{B}\otimes \mathcal{B}
\end{equation}
and \eqref{conds-on-K(r;ij')} becomes
\begin{equation}\label{orth-cond-on-K(ij)}
\sum_{i\in D} \overline{\hbox{Tr}}(K_{ij}^*K_{ij'})
=\delta_{j,j'}\ ,\quad \forall j,j'\in D
\end{equation}
In this case \eqref{df-cal-E-spec1} becomes
\begin{equation}\label{df-cal-E-spec2}
\mathcal{E}(a\otimes b)
=\sum_{i,j,i',j'\in D}
e_{ij}ae_{i'j'} \ \hbox{Tr}\left(K_{ij}^*\ b \ K_{i'j'}\right)
\end{equation}
\begin{definition}\label{df:trans-ampl-matr}{\rm
An operator $K\in \mathcal{B}\otimes \mathcal{B}$ satisfying \eqref{df-cond-dens-ampl} is called a \textbf{conditional density amplitude} (CDA) with respect to the trace. (Since in this paper we
will only consider CDA of this type, we simply call them CDA.)\\
A matrix $K\equiv(K_{jk})$ with entries in $\mathcal{B}$ satisfying \eqref{orth-cond-on-K(ij)}
is called a \textbf{transition amplitude matrix}
}\end{definition}
\textbf{Remark}.
Definition \ref{df:trans-ampl-matr} reflects the standard use in quantum mechanics of calling
\textit{amplitudes} the entries of a matrix $K=(K_{ij})$, such that the numbers
\begin{equation}\label{CCDA-26.7}
p_{ij}=| K_{ij} |^2
\end{equation}
are entries of a stochastic matrix
($p_{i,j}\ge 0 , \, \forall i,j\in D \, , \,\sum_{j\in D} p_{i,j}= 1$).\\
Recalling that a matrix $V\equiv( V_{ij})_{i,j\in D}$ is an isometry if and only if
\begin{equation}\label{iso-cond-V}
(V^*V)_{j'j}=\sum_{i\in D}\overline{V}_{ij}V_{ij'} = \delta_{j'j}
\end{equation}
and comparing \eqref{iso-cond-V} with \eqref{orth-cond-on-K(ij)}, one sees that the 
$\mathcal{B}$--valued matrix $(K_{ij})$ can be considered as a kind of 
\textit{operator valued isometry}.\\ Denoting
\begin{equation}\label{phk=|Vkh|2}
p_{hk} := | V_{kh} |^2 \in \mathbb{R}_{+}
\end{equation}
one can identify $V_{kh}$ with a square root $\sqrt{p_{hk}}$ of $p_{hk}$ and \eqref{iso-cond-V}
implies that
$$
\sum_{k\in D}p_{hk}
=\sum_{k\in D}\overline{V}_{kh}V_{kh} = 1
$$
i.e. $P\equiv (p_{hk})_{h,k\in D}$ is a stochastic matrix. Similarly, defining
\begin{equation}\label{Pij=|Kij|2}
p_{hk} =  \hbox{Tr}(|K_{kh}|^2) \in \mathbb{R}_{+}
\end{equation}
and putting $j=j'$ in \eqref{orth-cond-on-K(ij)}, one obtains
\begin{equation}\label{orth-cond-on-K(ij)2}
\sum_{i\in D} \overline{\hbox{Tr}}(K_{ij}^*K_{ij})
=\sum_{i\in D} \overline{\hbox{Tr}}(|K_{ij}|^2)
= \sum_{i\in D}p_{ji}
= 1 \ ,\quad \forall j\in D
\end{equation}
again $P\equiv (p_{hk})_{h,k\in D}$ is a stochastic matrix.\\

\noindent A (homogeneous, backward) quantum Markov chain (or quantum Markov state) on
$\mathcal{A}:=\bigotimes _\mathbb{N}\mathcal{B}$ is uniquely determined by a pair
$(\varphi_0,\mathcal{E})$, where $\varphi_0$ is a state on $\mathcal{B}$ (initial state) and
$\mathcal{E}:\mathcal{B}\otimes \mathcal{B}\to \mathcal{B}$ a transition expectation,
in the sense that its joint expectations are given by
$$
\varphi\bigl(j_0(a_0)j_1(a_1)\cdots j_n(a_n)\bigr)
$$
\begin{equation}\label{joint-exp-QMC}
=\varphi_0(\mathcal{E}(a_0\otimes \mathcal{E}(a_1\otimes \cdots \otimes \mathcal{E}(a_{n-1}
\otimes \mathcal{E}(a_n\otimes 1))\cdots ))
\end{equation}
When $\hbox{dim}(H)=d<\infty$ we will use the notation
$$
D := \{1,\dots, d\}
$$
To each ortho--normal basis (o.n.b.) $e\equiv(e_{h})_{h\in D}$ of $H$, one can associate a system of
matrix units $(e_{h,k})$ where, for all $h, k\in D$,
\begin{equation}\label{notat-ee^*}
e_{h,k}(\xi) :=
e_{k}e_{h}^*(\xi) := \langle e_{h}, \xi\rangle e_{k}  \quad,\quad\forall \xi\in H
\end{equation}
and its matrix, in the $e$--basis, has all entries equal to zero with the exception of the $(h,k)$--th element which is $1$.
In this case, the quadruple \eqref{tens-KR-alg-stoch-proc} in Definition \ref{df:tens-KR-alg-stoch-proc} is reduced to
\begin{equation}\label{tens-KR-alg-stoch-proc-homog}
\left(\mathcal{A}:=\bigotimes_{\mathbb{N}}\mathcal{B},\varphi , \mathcal{B}_{n},
(j_{n})_{n\in T}\right)
\end{equation}
For each $j\in D_{\mathcal{E}}$, $K_{j}\in M_d(\mathbb{C})\otimes M_d(\mathbb{C})$ can be
written in the form
\begin{equation}\label{(1.4)}
K=\ \sum_{h,h'\in D}e_{h,h'}\otimes\ K_{h,h'}
=\ \sum_{h,h'\in D}e_{h}e_{h'}^*\otimes\ K_{h,h'}
\end{equation}
In the following we denote $\mathcal{D}_{e}$ the $e$--diagonal sub--algebra of $M_d(\mathbb{C})$ defined by
\begin{equation}\label{df-De}
\mathcal{D}_{e} := \{\sum_{h=1}^d x_{h}e_{hh} \colon x_{h}\in\mathbb{C} , h\in D\}
\equiv \ell_{\mathbb{C}}^{\infty}(D)
\end{equation}
where $L_{\mathbb{C}}^{\infty}(D)$ denotes the space of all functions $f:D\to \mathbb{C}$.\\
The $e$--diagonal sub--algebra of $\mathcal{A}$ is defined by
\begin{equation}\label{df-calDe}
\mathcal{D}_{e} := \bigotimes_{\mathbb{N}} D_{e}
\equiv \bigotimes_{\mathbb{N}} L_{\mathbb{C}}^{\infty}(D)
\end{equation}
Therefore, if $\varphi\equiv(\varphi_0,\mathcal{E})$ is any quantum Markov state on $\mathcal{A}$,
for any diagonal algebra $\mathcal{D}_{e}$, the restriction of $\varphi$ on $\mathcal{D}_{e}$
defines a unique classical process $X\equiv (X_n)$, with state space $D$, characterized by the joint probabilities
$$
\hbox{Prob}\Bigl(X_0=i_0, X_1=i_1,\cdots ,X_n=i_n\Bigr)=
\varphi\bigl(e_{i_0,i_0}\otimes e_{i_1,i_1}\otimes \cdots \otimes
e_{i_{n},i_{n}}\bigr)=
$$
\begin{equation}\label{(2.1)}
=\varphi_0\left(\mathcal{E}(e_{i_0,i_0}\otimes \mathcal{E}(e_{i_1,i_1}\otimes
\cdots \otimes \mathcal{E}(e_{i_{n-1}}e_{i_{n-1}}^*\otimes
\mathcal{E}(e_{i_n,i_n}\otimes 1))\cdots ))\right)
\end{equation}
for any $n\in\mathbb{N}$, $\{i_h\}_{h=0}^n\subset \{1,\cdots, d\}$.
With the identifications \eqref{df-De}, \eqref{df-calDe} the restriction of the  embedding $j_n$,
defined by \eqref{df-jn-tens}, to $D_{e}$ can be identified to
\begin{equation}\label{(1.6)}
j_n(f):=f(X_n),\ \ \ \forall\ f\in \mathcal{D}
\end{equation}

\subsection{$e$--\textbf{diagonalizable} quantum Markov chains}\label{sec:e-diag-QMC}

Recall that any classical Markov chain
$(X_n)_{n\in\mathbb{N}}$ with state space $D=\{1,\cdots, d\}$ is determined by a pair
$(p^{(0)}, P\equiv(p_{i,j}))$, where $p^{(0)}$ is a probability measure on $D$ (initial distribution) and $P\equiv(p_{i,j})$ is a stochastic matrix and can be embedded
in a quantum Markov chain \cite{Ac74-Camerino}, \cite{Ac74-FAA}.
In fact, for any choice of the entry--wise square root $\sqrt{P}:=(\sqrt{p_{k,h}})$ of $P$ defining
\begin{equation}\label{K-e-diag}
K:= \sum_{h\in D}\sqrt{p_{k,h}} e_{k,k}\otimes e_{h,h}
\end{equation}
\begin{equation}\label{calE-e-diag}
\mathcal{E}(x):= \overline{\hbox{Tr}}_2(K^*\ x\ K) \, , \quad
x\in M_d(\mathbb{C})\otimes M_d(\mathbb{C})
\end{equation}
for any
$x\otimes y := \sum_{m,n\in D}x_{m}y_{n} e_{m,m}\otimes e_{n,n}\in D_{e}\otimes D_{e}$, one has,
\begin{align*}
\mathcal{E}(x)=&\overline{\hbox{Tr}}_2(K^*\ x\ K)\\
=&\sum_{m,n,h,h',k,k'\in D}x_{m}y_{n}\overline{\sqrt{p_{k,h}}} \sqrt{p_{k',h'}}\overline{\hbox{Tr}}_2( (e_{k,k}\otimes e_{h,h)} (e_{m,m}\otimes e_{n,n}) (e_{k',k'}\otimes e_{h',h'}))\\
=&\sum_{h,k\in D}x_{k}y_{h}\overline{\sqrt{p_{k,h}}} \sqrt{p_{k,h}}\overline{\hbox{Tr}}_2(
(e_{k,k}\otimes e_{h,h}) )=\sum_{k\in D}x_{k}\left(\sum_{h\in D} p_{k,h}y_{h}\right)  e_{k,k}\\
=& xP(y)
\end{align*}
\eqref{calE-e-diag} shows that $\mathcal{E}$ is a transition expectation on $M_d(\mathbb{C})$
that can be written in the form \eqref{df-cal-E-Tr-CDA} with $|D_{\mathcal{E}}|=1$ and satisfying:
\begin{equation}\label{calE(De-otimes-De)-in-De}
\mathcal{E}(D_{e}\otimes D_{e}) \subseteq D_{e}
\end{equation}
Therefore, for any initial state on $M_d(\mathbb{C})$
\begin{equation}\label{(1.9)}
\varphi_0(\cdot ):=\hbox{Tr}\bigl(W_{0} \ \cdot \ ) \ ;\quad W_{0}:= \sum_{j\in D}p_{j}e_{j,j}
\end{equation}
the pair $(\varphi_0,\mathcal{E})$ defines a unique quantum Markov chain whose restriction
on the $e$--diagonal sub--algebra $\mathcal{D}_{e}$ is the classical Markov chain
$(X_n)_{n\in\mathbb{N}}$ described above.\\
In conclusion: any classical Markov chain can be obtained (usually in many ways) as the restriction
of a quantum Markov chain.
\begin{definition}\label{df:diagonalizable-QMC}{\rm
The quantum Markov chains with transition expectations characterized by
\eqref{K-e-diag}, \eqref{calE-e-diag} are called $e$--\textbf{diagonalizable}.
}\end{definition}

\subsection{Algebraic (classical or quantum) hidden processes}\label{sec:QHMP}

The structure \label{joint-exp-O-H-Mark} of the joint expectations of a classical
hidden Markov chain in the form \eqref{joint-probs-O-H-Mark-ker}, i.e.
\begin{equation}\label{joint-exp-O-H-Mark2}
P_{H,O}\left(\prod_{m=0}^{n}f_{m}(O_{m})g_{m}(H_{m})\right)
\end{equation}
$$
= p_{H_{0}}\left(E_{H_{0}}\left(B_{O_{0}}(f_{0})g_{0}\right)
P_{H_{0}}\left(B_{O_{1}}(f_{1})g_{1}\right)\left(\cdots \right.\right.
$$
$$
\left.\left.
P_{H_{n-1}}(B_{O_{n-1}}(f_{n})g_{n-1})
P_{H_{n-1}}(B_{O_{n}}(f_{n})g_{n})\right)\right)
$$
naturally suggests a way to extend this notion to the quantum case.\\ In fact, the difficulty
to interpret \eqref{joint-exp-O-H-Mark2} in the case where the algebras
$L^{\infty}_{\mathbb{C}}(S_{H_{n}})$ and $L^{\infty}_{\mathbb{C}}(S_{O_{n}})$ are replaced by
non--commutative algebras $\mathcal{B}_{H_{n}}$ and $\mathcal{B}_{O_{n}}$ respectively, is that
maps of the form  \eqref{non-pos-map}, i.e.
$$
f_{n}\otimes g_{n} \mapsto B_{O_{n}}(f_{n})g_{n}
$$
are not positive in general.
This problem is a special case of the problem met in the quantum extension of the classical Markov chains and therefore it can be solved using the same idea introduced in \cite{Ac74-Camerino}, \cite{Ac74-FAA}, namely: replacing the map \eqref{non-pos-map} by
a transition expectation,
$$
\mathcal{E}_{O,H;n}
\colon \mathcal{B}_{O_{n}}\otimes \mathcal{B}_{H_{n}}\to \mathcal{B}_{H_{n}}\ ,\quad n\in\mathbb{N}
$$
In usual markovianity, all the algebras are commutative and one only considers transition expectations uniquely determined by a sequence of Markov operators
$$
B_{O_{n}}\colon \mathcal{B}_{O_{n}}\to  \mathcal{B}_{H_{n}}
$$
through the identity
\begin{equation}\label{E(O,H;n)-comm-case}
\mathcal{E}_{O,H;n}(f_{n}\otimes g_{n}):=B_{O_{n}}(f_{n})g_{n}  \ ;\quad f_{n}\in\mathcal{B}_{O_{n}} \ , \ g_{n}\in\mathcal{B}_{H_{n}}
\end{equation}
Thus the empirical rules for the transition from classical hidden processes to quantum hidden processes are the following:\\

\noindent
(i) Replace $L^{\infty}_{\mathbb{C}}(S_{H_{n}})$ and $L^{\infty}_{\mathbb{C}}(S_{O_{n}})$
respectively by arbitrary $*$--algebras $\mathcal{B}_{H_{n}}$ and $\mathcal{B}_{O_{n}}$.\\
(ii) Replace the operator $g_{n}\otimes f_{n}\in L^{\infty}_{\mathbb{C}}(S_{H_{n}})\otimes L^{\infty}_{\mathbb{C}}(S_{O_{n}})\mapsto B_{O_{n}}(f_{n})g_{n}$ by a transition expectation
\begin{equation}\label{df-B(O,H,m)2-q}
\mathcal{E}_{O,H,n} \colon\mathcal{A}_{H_{n}}\otimes \mathcal{A}_{O_{n}}\to
\mathcal{A}_{H_{n}}
\end{equation}
With these replacements, the sample algebra of the classical process\\
$\mathcal{A}_{H,O} := \bigotimes_{n\in \mathbb{N}}
\left( L^{\infty}_{\mathbb{C}}(S_{H_{n}})\otimes L^{\infty}_{\mathbb{C}}(S_{O_{n}})\right)$
(see \eqref{df-alg-HO}) is replaced by
\begin{equation}\label{df-HO-alg-q}
\mathcal{A}_{H,O} :=
\bigotimes_{n\in\mathbb{N}}(\mathcal{B}_{O_{n}}\otimes \mathcal{B}_{H_{n}})
\end{equation}
and the global emission operator from $\mathcal{A}_{H,O}$ to the hidden algebra
\begin{equation}\label{df-hidd-alg-q}
\mathcal{A}_{H} := \bigotimes_{n\in\mathbb{N}}\mathcal{B}_{H_{n}}
\end{equation}
is replaced by
\begin{equation}\label{df-cal-E(O,H)clas}
\mathcal{E}_{O,H}:=\bigotimes_{n\in\mathbb{N}}\mathcal{E}_{O,H;n} \colon
\bigotimes_{n\in\mathbb{N}}(\mathcal{B}_{O_{n}}\otimes \mathcal{B}_{H_{n}})
=\mathcal{A}_{H,O}\to \mathcal{A}_{H} = \bigotimes_{n\in\mathbb{N}}\mathcal{B}_{H_{n}}
\end{equation}
In analogy with the classical case, the transition expectation $\mathcal{E}_{O,H,n}$
(which is a particular Markov operator) is called the \textbf{$n$--th emission operator}.
Also in the quantum case the hidden process $H\equiv\{H_n\}_{n\in\mathbb{N}} $,
\textbf{is arbitrary}.
\begin{theorem}\label{th:joint-exps-(O,H)-q}{\rm
Let be given:\\
-- a quantum stochastic process called the hidden process,
\begin{equation}\label{df-QHP-calA(H)}
(\mathcal{A}_{H},P_{H}, (\mathcal{B}_{H_{n}})_{n\in \mathbb{N}}, (j_{H_{n}})_{n\in \mathbb{N}})
\end{equation}
where $\mathcal{A}_{H}$ is given by \eqref{df-hidd-alg-q} and the $j_{H_{n}}$ are the usual tensor embeddings;\\
-- a family of $*$--algebras $(\mathcal{B}_{O_{n}})_{n\in \mathbb{N}}$;\\
-- a family of emission operators
$\mathcal{E}_{O,H,n} \colon\mathcal{B}_{H_{n}}\otimes \mathcal{B}_{O_{n}}\to \mathcal{B}_{H_{n}}$
$ (n\in \mathbb{N})$ and the associated global emission operator given by
\eqref{df-cal-E(O,H)clas}.\\
There exists a, unique up to stochastic equivalence, quantum stochastic process
\begin{equation}\label{df-QSP-calA(H,O)}
(\mathcal{A}_{H,O},P_{H,O}, (\mathcal{B}_{H_{n}}\otimes \mathcal{B}_{O_{n}})_{n\in \mathbb{N}}, (j_{O_{n}}\otimes j_{H_{n}})_{n\in \mathbb{N}})
\end{equation}
such that $\mathcal{A}_{H,O}$ is given by \eqref{df-HO-alg-q}, the $j_{O_{n}}\otimes j_{H_{n}}$
are the usual tensor embeddings and
\begin{equation}\label{P(H,O)=PH-circ-E(O,H)-q}
P_{H,O} = P_{H} \circ \mathcal{E}_{O,H}
\end{equation}
or, more explicitly,
\begin{equation}\label{joint-exps-(O,H)-B-expl-q}
P_{H,O}\Big(\bigotimes_{m=0}^{n}j_{H_{m}}(g_{m})\otimes j_{O_{m}}(f_{m})\Big) =P_{H}\Big(\bigotimes_{m=0}^{n} j_{H_{m}}(\mathcal{E}_{O,H,m}(g_{m}\otimes f_{m}))\Big)
\end{equation}
for all $n\in\mathbb{N}$, $m\in \{0, \dots, n\}$, $f_{m}\in \mathcal{B}_{O_{n}}$ and $g_{m}\in \mathcal{B}_{H_{n}}$.}
\end{theorem}
\textbf{Proof}.
Since $\mathcal{E}_{O,H} \colon \mathcal{A}_{H,O}\to \mathcal{A}_{H}$ is a Markov operator
and $P_{H}$ is a state on $\mathcal{A}_{H}$,
$P_{H,O} = P_{H} \circ \mathcal{E}_{O,H}$ is a state on $\mathcal{A}_{H,O}$.
Therefore the quantum stochastic process \eqref{df-QSP-calA(H,O)} is well defined.
$\qquad\square$
\begin{definition}{\rm
The quantum stochastic process \eqref{df-QSP-calA(H,O)} constructed in Theorem \ref{th:joint-exps-(O,H)-q} is called a quantum hidden process with hidden process
given by \eqref{df-QHP-calA(H)} and observable process given by
\begin{equation}\label{df-QOP-calA(O)}
(\mathcal{A}_{O},P_{O}, (\mathcal{B}_{O_{n}})_{n\in \mathbb{N}}, (j_{O_{n}})_{n\in \mathbb{N}})
\end{equation}
where
\begin{equation}\label{df-obs-alg-q}
\mathcal{A}_{O} := \bigotimes_{n\in\mathbb{N}}\mathcal{B}_{O_{n}}
\equiv \bigotimes_{n\in\mathbb{N}}\mathcal{B}_{O_{n}}\otimes 1_{\mathcal{A}_{H}}
\end{equation}
$P_{O}$ is the restriction of $P_{H,O}$ to $\mathcal{A}_{O}$ with the identification
\eqref{df-obs-alg-q} (obtained by putting all the $g_{m}=1_{\mathcal{A}_{H}}$ in
\eqref{joint-exps-(O,H)-B-expl-q}) and the $j_{O_{n}}$ are the usual tensor embeddings.
}\end{definition}

\section{Quantum hidden Markov processes}\label{sec:Q-hi-Mark-proc}

As in the classical case, the quantum hidden Markov processes are those quantum hidden processes whose hidden process is a quantum Markov process (in the present paper Markov chain).\\
In this Section, we consider non necessarily homogeneous (backward) Markov chains
(see \cite{Ac74-Camerino}, \cite{Ac74-FAA}).   
They are states on the algebra $\mathcal{A}_{H}$ given by \eqref{df-hidd-alg-q},
uniquely determined by a pair
\begin{equation}\label{pair-def-H-QMC}
(P_{H_{0}} \ , \ (\mathcal{E}_{H_{n}})_{n\in\mathbb{N}})
\end{equation}
where $P_{H_{0}}$ is a state on $\mathcal{B}_{H_{0}}$ and
\begin{equation}\label{backw-trans-exp}
\mathcal{E}_{H_{n}}\colon \mathcal{B}_{H_{n}}\otimes \mathcal{B}_{H_{n+1}} \to \mathcal{B}_{H_{n}}
\quad,\quad n\in\mathbb{N}
\end{equation}
a family of backward transition expectations.
Such a pair determines the joint expectations through the identities
\begin{align}\label{joint-exp-QMC-PH}
&P_{H}\Big(\prod_{m=0}^{n}j_{H_{m}}(g_{m})\Big)\\
=&P_{H_{0}}\left(\mathcal{E}_{H_{0}}\left(g_{0}\otimes \mathcal{E}_{H_{1}}\left(\cdots\otimes \mathcal{E}_{H_{n-1}}\left(g_{n-1}\otimes
\mathcal{E}_{H_{n}}\left(g_{n}\otimes 1_{H_{n+1}}\right)\right)\right) \right)\right)\notag
\end{align}
for all $g_{m}\in\mathcal{B}_{H_{m}}$, $m\in\{1,\dots,n\}$.
\begin{theorem}\label{th:struct-QHMP}{\rm
In the notations of Theorem \ref{th:joint-exps-(O,H)-q}, suppose that the quantum hidden process
\eqref{df-QHP-calA(H)} is the backward Markov process characterized by the pair \eqref{pair-def-H-QMC}.
Then the process \eqref{df-QSP-calA(H,O)} is the hidden quantum Markov process whose joint
expectations are given by
\begin{align}\label{joint-exp-qHMP}
&P_{H,O}\Big(\bigotimes_{m=0}^{n}j_{H_{m}}(g_{m})\otimes j_{O_{m}}(f_{m})\Big)\\
=& P_{H_{0}}\left(\mathcal{E}_{H_{0}}\left(\mathcal{E}_{O,H;0} (f_{0}\otimes g_{0})\otimes\mathcal{E}_{H_1} (\mathcal{E}_{O,H;1} \left((f_{1}\otimes g_{1})\otimes \cdots \right.\right.\right. \notag\\
&\qquad \left.\left.\left.\otimes \mathcal{E}_{H_{n-1}}\left( \mathcal{E}_{O,H;n}(f_{n-1}\otimes g_{n-1})\otimes
\mathcal{E}_{H_{n}}\left(\mathcal{E}_{O,H;n}\left(f_{n}\otimes g_{n}\right)\otimes 1_{H_{n+1}}
\right)\right)\right)\right)\right)\notag
\end{align}
for all $n\in\mathbb{N}$, $m\in \{0, \dots, n\}$, $f_{m}\in \mathcal{B}_{O_{n}}$ and
$g_{m}\in \mathcal{B}_{H_{n}}$. Moreover, introducing the notation
\begin{equation}\label{notat-calE(m;fm)}
\mathcal{E}_{m;f_{m}} := \mathcal{E}_{O,H,m}(1_{H_{m}}\otimes f_{m})
\end{equation}
the joint expectations of the observable process, obtained from \eqref{joint-exp-qHMP} by
putting $g_{m}=1$ for all $m\in \{1,\dots,n\}$, are given by
\begin{equation}\label{joint-exps-O-proc-q}
P_{O}\left(\bigotimes_{m=0}^{n}j_{O_{m}}(f_{m})\right)
\end{equation}
$$
= P_{H_{0}}\left(\mathcal{E}_{H_{0}}\left(\mathcal{E}_{O,H;0} (f_{0}\otimes \mathcal{E}_{H_1}(\mathcal{E}_{1;f_{1}}\otimes  \cdots\right.\right.\left.\left.\left.\otimes \mathcal{E}_{H_{n-1}}\left(\mathcal{E}_{n-1;f_{n-1}}\otimes
\mathcal{E}_{H_{n}}\left(\mathcal{E}_{n;f_{n}}\otimes 1_{H_{n+1}}\right)\right)\right)\right)\right)
$$
for all $n\in\mathbb{N}$, $m\in \{0, \dots, n\}$, $f_{m}\in \mathcal{B}_{O_{n}}$.
}\end{theorem}
\textbf{Remark}. Formula \eqref{joint-exps-O-proc-q} highlights the double tier structure of quantum hidden Markov processes: the operators $\mathcal{E}_{j; \, \cdot \, }$ transform elements in $\mathcal{B}_{O_{n}}$ (in the commutative case functions $f_{j}$ of the $j$--th observable random variable $O_{j}$)
into functions of the $j$--th hidden random variable $H_{j}$) and on these transformed functions the transition expectations $\mathcal{E}_{H_{j}}$ act as for usual Markov chains (see
\eqref{joint-exp-qHMP}).\\

\noindent\textbf{Proof}. The identity \eqref{joint-exps-O-proc-q} is obtained from \eqref{joint-exp-qHMP} replacing in the
right hand side the $g_{m}$ by the $\mathcal{E}_{O,H;m} (f_{m}\otimes g_{m})$. In the notation \eqref{notat-calE(m;fm)}, the identity \eqref{joint-exps-O-proc-q} is
obtained replacing in the right hand side of \eqref{joint-exp-qHMP} the $g_{m}$ by
$1_{\mathcal{B}_{H_{m}}}$.
$\square$\\

\noindent\textbf{Remark}.
Recall that, if all the algebras $\mathcal{B}_{O_{n}}$ and $\mathcal{B}_{H_{n}}$ are commutative,
the transition expectations are given by
\begin{equation}\label{E(HO;n)-comm-case}
\mathcal{E}_{H_{n}}\left(\mathcal{E}_{O,H;n}\left(f_{n}\otimes g_{n}\right)\otimes 1_{H_{n+1}}
\right)
:= P_{H_{n}}\left(B_{O,H,n}(f_{n})g_{n}\right)
\end{equation}
In this case, the joint expectations \eqref{joint-exp-qHMP} coincide with \eqref{joint-probs-O-H-Mark},
i.e. with the classical hidden Markov process with Markov operators $(P_{H_{n}})$ and
emission operators $(B_{O,H,n})$.\\

\subsection{Diagonalizable Markov chains and associated classical hidden processes}
\label{Diag-QMC-Cl-HMP}

In this section we prove that a special class of quantum Markov chains is strictly related to
classical hidden Markov processes in the sense that: the restriction of any element in this class
to any diagonal (in particular commutative) sub--algebra produces a HMP. This means that the
assignment of a single diagonalizable quantum Markov chain automatically gives uncountably 
many classical hidden Markov processes. This fact was known since the very early times of the
theoru of QMC (see \cite{Ac91-Q-Kalman-filters}). \\

\noindent The above mentioned special class of Markov chains is obtained particularizing
the construction, described at the beginning of Section \ref{sec:Q-hi-Mark-proc}, by introducing
$3$ types of additional conditions.\\

\noindent \textbf{(i) Assumptions on the algebras}\\

\noindent All the algebras $\mathcal{A}_{H_{n}}$ and $\mathcal{A}_{O_{n}}$ are taken to be
isomorphic to a single algebra $\mathcal{B}$, \textbf{independent of} $n$
and isomorphic to the algebra of all bounded operators on a Hilbert space $\mathcal{H}$:
\begin{equation}\label{hid-obs-q-algs-equal}
\mathcal{A}_{H_{n}} \equiv \mathcal{A}_{O_{n}} \equiv \mathcal{B}:=\mathcal{B}(\mathcal{H})
\end{equation}
where $\equiv$ denotes $*$--isomorphism.
Thus, for the sample algebra of the underlying Markov process, one has the identification
$$
\mathcal{A}_{H} := \bigotimes_{\mathbb{N}}\mathcal{B}
$$
with the tensor embeddings
\begin{equation}\label{df-jHn-calA-Hn}
j_{H_{n}} \colon b\in \mathcal{B} \to j_{H_{n}}(b) \equiv b\otimes 1_{\{n\}^c}\in\mathcal{A}_{H}
\quad;\quad \mathcal{A}_{H_{n}} := j_{H_{n}}(\mathcal{B})
\end{equation}
where $1_{\{n\}^c}$ denotes the identity in $\bigotimes_{\mathbb{N}\setminus\{n\}}\mathcal{B}$.\\
Similarly we define the observable algebra
$$
\mathcal{A}_{O} \equiv \bigotimes_{\mathbb{N}}\mathcal{B}
$$
and the tensor embeddings
\begin{equation}\label{df-jOn-calA-On}
j_{O_{n}} \colon b\in \mathcal{B} \to j_{O_{n}}(b) \equiv b\otimes 1_{\{n\}^c}\in\mathcal{A}_{O}
\quad;\quad \mathcal{A}_{O_{n}} := j_{O_{n}}(\mathcal{B})
\end{equation}
where again $1_{\{n\}^c}$ denotes the identity in
$\bigotimes_{\mathbb{N}\setminus\{n\}}\mathcal{B}$.\\
The algebra of the $(H,O)$--process is then
\begin{equation}\label{identif-calA(H,O)}
\mathcal{A}_{H,O} := \mathcal{A}_{H_{n}} \otimes \mathcal{A}_{O_{n}}
\equiv\bigotimes_{\mathbb{N}} \ (\mathcal{A}_{H} \otimes \mathcal{A}_{O})
\equiv\bigotimes_{\mathbb{N}}\mathcal{B}\otimes\bigotimes_{\mathbb{N}}\mathcal{B}
\end{equation}
\textbf{(ii) Assumptions on the transition expectations}\\

\noindent We suppose that the backward transition expectations of the $H$--process do not depend
on $n\in\mathbb{N}$:
\begin{equation}\label{backw-trans-exp-comm}
\mathcal{E} \colon \mathcal{B}\otimes \mathcal{B} \to \mathcal{B}
\end{equation}
and that $\mathcal{E}$ has the form
\begin{equation}\label{purely-gen-calE}
\mathcal{E}(x) := \overline{\hbox{Tr}}_{2}(K^*x K) \quad,\quad x\in \mathcal{B}\otimes \mathcal{B}
\end{equation}
where $K\in \mathcal{B}\otimes \mathcal{B}$ is a conditional density amplitude commuting with
its right shift, namely:
\begin{equation}\label{[K-otimes1,1-otimes-K]=0}
[K\otimes 1 \ , \ 1\otimes K]=0
\end{equation}

\textbf{(iii) Assumptions on the initial state}\\

\noindent $\varphi_{0}=\hbox{Tr}(w_{0} \, \cdot \, )$ is a state on $\mathcal{B}$ satisfying
\begin{equation}\label{[w0,K0]=0}
[w_{0}\otimes 1,\ K]=0
\end{equation}
Moreover we assume that $w_{0}$ has non--degenerate spectrum, i.e. that
$$
w_{0} := \sum_{j\in D}w_{0;j}e_{jj} = \sum_{j\in D}w_{0;j}e_{j}e_{j}^*
$$
where $e\equiv (e_{j})_{j\in D}$ is an ortho--normal basis of $\mathcal{H}$.
We denote $(e_{hk}=e_{h}e_{k}^*)_{h,k\in S}$ the system of matrix units in $\mathcal{B}$ associated
to the basis $e$. In these notations, we also suppose the validity of the following
non-degeneracy condition:
\begin{equation}\label{Kkk-ne-0-a}
(e_{kk}\otimes 1)K\ne 0  \quad,\quad\forall k\in D
\end{equation}
\begin{theorem}\label{thm:char-diagonaliz-CDA-NDS}{\rm
Under the conditions (i), (ii), (iii) listed above, there exists a stochastic matrix
$P\equiv (p_{jk})$ such that
\begin{equation}\label{struct-K-NDS}
K = \sum_{j,k\in D}\sqrt{p_{jk}}  e_{jj}\otimes e_{kk}
\end{equation}
where, for each $j,k\in D$, $\sqrt{p_{jk}}$ is an arbitrary square root of $p_{jk}$.
Conversely, given a stochastic matrix $P\equiv (p_{jk})$, for any square root $\sqrt{p_{jk}}$
of $p_{jk}$ ($j,k\in D$), the matrix $K$, defined by the right hand side of \eqref{struct-K-NDS} is
a conditional density amplitude.
}\end{theorem}
\textbf{Proof}.
If $w_{0}$ has non--degenerate spectrum, then the commutant of $w_{0}\otimes 1$ in
$\mathcal{B}\otimes \mathcal{B}$ is
\begin{equation}\label{comm-w0-otimes-1-NDS}
\{\mathcal{B}\otimes \mathcal{B}\}' =
\{x\in \mathcal{B}\otimes \mathcal{B} \colon
x = \sum_{j\in D}(e_{jj}\otimes 1_{M})x(e_{jj}\otimes 1_{M}) \}
\end{equation}
\noindent Since $(e_{hk})_{h,k\in S}$ is a system of matrix units in $\mathcal{B}$, $K$
can be written in the form
\begin{equation}\label{CCDA-struct-K-a}
K = \sum_{h,k\in S} e_{hk}\otimes K_{hk}
\qquad,\quad P_{0;hk}\in\mathcal{B}
\end{equation}
with $K_{hk}\in \mathcal{B}$.
Given \eqref{CCDA-struct-K-a}, condition \eqref{Kkk-ne-0-a} can be equivalently written as
\begin{equation}\label{Kkk-ne-0}
K_{kk}\ne 0  \quad,\quad\forall k\in D
\end{equation}
Therefore \eqref{[w0,K0]=0} implies
$$
K = \sum_{j\in D} (e_{jj}\otimes 1_{M})K(e_{jj}\otimes 1_{M})
\overset{\eqref{CCDA-struct-K-a}}{=}
\sum_{j\in D}\sum_{h,k\in D} (e_{jj}\otimes 1_{M}) (e_{hk}\otimes K_{hk})(e_{jj}\otimes 1_{M})
$$
$$
= \sum_{j\in D}\sum_{h,k\in D} (e_{jj}e_{hk}e_{jj}\otimes K_{hk})
= \sum_{j,h,k\in D}\delta_{jh}\delta_{kj}(e_{jj}\otimes K_{hk})
$$
\begin{equation}\label{K-NDS}
= \sum_{j\in D}(e_{jj}\otimes K_{jj})
\end{equation}
\eqref{[K-otimes1,1-otimes-K]=0} implies that
$$
0= \sum_{j,k\in D}[(e_{jj}\otimes K_{jj})\otimes 1 \ , \ 1\otimes (e_{kk}\otimes K_{kk})]
$$
$$
= \sum_{j,k\in D} (e_{jj}\otimes K_{jj}e_{kk}\otimes K_{kk}
- e_{jj}\otimes e_{kk}K_{jj}\otimes K_{kk})
$$
$$
= \sum_{j\in D}e_{jj}\otimes\sum_{k\in D} (K_{jj}e_{kk}\otimes K_{kk}
- e_{kk}K_{jj}\otimes K_{kk})
$$
$$
\iff 0= \sum_{k\in D} (K_{jj}e_{kk}\otimes K_{kk} - e_{kk}K_{jj}\otimes K_{kk})
\quad,\quad\forall j\in D
$$
Multiplying on
the right by $e_{hh}\otimes 1$ with $h\ne k$, one obtains
$$
0= \sum_{k\in D} e_{kk}K_{jj}e_{hh}\otimes K_{kk}
\quad,\quad\forall j\in D
$$
$$
\iff 0=  e_{kk}K_{jj}e_{hh}\otimes K_{kk}
\quad,\quad\forall j,k\in D
$$
$$
\iff 0 = e_{kk}K_{jj}e_{hh} =0
\quad,\quad\forall j,k\in D  \mbox{ such that }  K_{kk}\ne 0
$$
$$
\overset{\eqref{Kkk-ne-0}}{\iff} 0 = e_{kk}K_{jj}e_{hh} =0
\quad,\quad\forall j,k\in D \ , \ \forall h\ne k
$$
and this is equivalent to
\begin{equation}\label{Kjj-NDS}
K_{jj} =\sum_{k\in D}  e_{kk}K_{jj}e_{kk}
=\sum_{k\in D} \langle e_{k}, K_{jj}e_{k}\rangle e_{kk}
=:\sum_{k\in D} r_{jk} e_{kk}
\quad,\quad\forall j\in D
\end{equation}
Replacing \eqref{Kjj-NDS} in \eqref{K-NDS}, one obtains
$$
K= \sum_{j\in D}(e_{jj}\otimes K_{jj})
= \sum_{j,k\in D}r_{jk}  e_{jj}\otimes e_{kk}
$$
The fact that $K$ is a conditional density amplitude implies that
$$
1 = \overline{\hbox{Tr}}_{2}(K^*K)
= \overline{\hbox{Tr}}_{2}((\sum_{j,k\in D}r_{jk}  e_{jj}\otimes e_{kk})^*(\sum_{j,k\in D}r_{jk}  e_{jj}\otimes e_{kk}))
$$
$$
= \overline{\hbox{Tr}}_{2}((\sum_{j,k\in D}\overline{r_{jk}}r_{jk}    e_{jj}\otimes e_{kk})
= \sum_{j\in D}\left(\sum_{k\in D} |r_{jk}|^2\right) e_{jj}
$$
\begin{equation}\label{normaliz-diag-K-NDS}
\iff \sum_{k\in D} |r_{jk}|^2   \quad,\quad\forall j\in D
\end{equation}
Therefore the matrix with entries
$$
p_{jk} := |r_{jk}|^2 \quad,\quad\forall j,k\in D
$$
is a stochastic matrix and each $r_{jk}$ is a square root of the corresponding $p_{jk}$.
This proves the first statement.
The proof of the converse statement is included in the calculation \eqref{normaliz-diag-K-NDS}.
$\qquad\square$\\
zzz
\noindent One proves that, for any normal state $P_{0}$ on $\mathcal{B}$ with density operator
$w_{0}$, the pair $(P_{0},\mathcal{E})$ uniquely defines a state $P_{H} $ on
$\mathcal{A}_{H}\equiv\bigotimes_{\mathbb{N}} \mathcal{B} $ which is a quantum Markov chain
and its joint expectations have the form \eqref{joint-exp-QMC-PH}.
The (forward) \textbf{time--shift} $u$ is the $*$--endomorphism of $\mathcal{A}$ characterized
by the property
$$
u\circ j_{H_{n}} := j_{H_{n+1}} \quad,\quad\forall n\in\mathbb{N}
$$
For any conditional density amplitude $K$ and any $n$, one defines
\begin{equation}\label{CCDA-shifted-K}
K_{[n,n+1]} := u^{n}((j_{H_{0}}\otimes j_{H_{1}})(K))
= (j_{H_{n}}\otimes j_{H_{n+1}})(K)
\end{equation}
and one proves that, denoting
$$
W_{[0,n+1]} :=
$$
\begin{equation}\label{CCDA-dens_matr_QM1}
K_{[n,n+1]} K_{[n-1,n]}\dots  K_{[0,1]} w_0 K_{[0,1]}^* K_{[1,2]}^*\dots K_{[n,n+1]}^*
\in \bigl( \otimes\mathcal{B} \bigr)^{n+2}
\end{equation}
one has
\begin{equation}\label{CCDA-dens_matr_QM2}
\overline{\hbox{Tr}}_{n+1}\left(
K_{[n,n+1]} K_{[n-1,n]}\dots  K_{[0,1]} w_0 K_{[0,1]}^* K_{[1,2]}^*\dots K_{[n,n+1]}^*\right)
\in \bigl( \otimes\mathcal{B} \bigr)^{n+1}
\end{equation}
$$
=\hbox{density matrix of } \ P_{H,[0,n]}
$$
Conversely, given a conditional density amplitude (CDA) $K$
and an initial density operator $w_0 $ on $\mathcal{B} \cong \mathcal{A}_0 $ and defining
the $W_{[0,n+1]}$ by \eqref{CCDA-dens_matr_QM1}, the limit
$$
P_{H} (a) := \lim_{n\to\infty} \hbox{Tr}(W_{[0,n]}a)
\quad,\quad \forall a\in\mathcal{A}
$$
exists in the strongly finite sense on $\mathcal{A}$ and satisfies \eqref{CCDA-dens_matr_QM2}.
Equivalently:
\begin{equation}\label{CCDA-varphi}
P_{H}(a)
= Tr(W_{[0,n+1]}a) \qquad,\quad \forall a\in \mathcal{A}_{[0,n]}
\end{equation}
and $P_{H}$ is the quantum Markov chain
generated by by the pair $(\varphi_{0},\mathcal{E})$ with $\mathcal{E}$ given by
\eqref{purely-gen-calE} (see \cite{Ac81-Topics-QP} for more information on this construction).\\

\subsubsection{Diagonalizable conditional density amplitudes and HMP}\label{sec:Diag-CDA-HMP}

\noindent The simplest examples of quantum Markov chains can be constructed by fixing
an ortho--normal basis $e_{H}\equiv (e_{H;j})_{j\in D}$ of $\mathcal{H}$ (where $D$ is a set
whose cardinality is equal to the dimension of $\mathcal{H}$) and choosing $K$ satisfying the
conditions of Theorem \ref{thm:char-diagonaliz-CDA-NDS} with respect to this basis.
Then we know that $K$ is uniquely determined by a stochastic matrix $P_{H}\equiv (p_{H;i,j})$
through the identity
\begin{equation}\label{Comm-CDA-HM}
K_{H} := \sum_{i,j\in D} \sqrt{p_{H;i,j}} e_{H;ii}\otimes e_{H;jj}
\end{equation}
where $e_{H;ij}:=e_{H;i}e_{H;j}^*$ is the system of matrix units associated to the ONB $e_{H}$ and,
for each $i,j\in D$, $e_{H;i}e_{H;j}^*$ is defined by
\begin{equation}\label{df-eiej^*}
e_{H;ij}\xi=e_{H;i}e_{H;j}^*(\xi) := \langle e_{H;j}, \xi\rangle e_{H;i}
\qquad,\quad\forall \xi\in \mathcal{H}
\end{equation}
CDA of the form \eqref{Comm-CDA-HM} (and the associated QMC) are called
rank--$1$--\textbf{diagonalizable}. For such QMC, defining
\begin{equation}\label{CCDA-26.4}
P_0=K_{[0,1]}^*K_{[0,1]}\qquad \;\qquad P_n=u^n(P_0)\qquad ;\qquad n\in \mathbb{N}
\end{equation}
one verifies that
\begin{equation}\label{CCDA-26.5}
W_{[0,n]}=w_0\cdot P_0\cdot P_1\cdots  P_{n}
\end{equation}
$$
=\sum_{j_1,\dots, j_n\in D} p_{H;j_0}^{0} p_{H;j_0j_1}\dots p_{H;j_{n-1}j_n}
e_{H;j_0j_0}\otimes e_{H;j_1j_1}\otimes\dots \otimes e_{H;j_nj_n}
$$
Thus diagonalizable QMC are \textit{diagonal liftings} of classical Markov chains.

\subsubsection{ Classical hidden Markov processes associated to diagonalizable quantum Markov chains:
maximal observables }\label{sec:QMC-HMP-max-obs}

In this section we prove that \textbf{the restriction of any diagonalizable quantum Markov state
on any diagonal algebra is a classical hidden Markov process}.
\begin{definition}\label{df-calD-e(on)}{\rm
In the notations of section \ref{Diag-QMC-Cl-HMP}, for each $n\in\mathbb{N}$, fix an ortho--normal
basis $e_{O}\equiv (e_{O_{n};j})_{j\in D}$ of $\mathcal{H}$ and define the $e_{O_{n}}$--diagonal
algebra by
$$
\mathcal{D}_{e_{O_{n}}}
:= \hbox{lin. span}\left\{e_{O_{n};jj} \colon j\in D\right\}
\subset j_{O_{n}}(\mathcal{B})\subset \mathcal{A}_{O_{n}}
$$
where $e_{O_{n};jj}:= e_{O_{n};j}e_{O_{n};j}^*$ (see \eqref{df-eiej^*}).
The algebra
\begin{equation}\label{eO-diag-alg}
\mathcal{D}_{(e_{O})} :=\bigvee_{n\in\mathbb{N}}\mathcal{D}_{e_{O_{n}}}
\left(\sim \bigotimes_{n\in\mathbb{N}}\mathcal{D}_{e_{O_{n}}}\right)
\subset \mathcal{A}_{O}
\end{equation}
is called the diagonal algebra with respect to the sequence of bases  $e_{O_{n}}$ or simply the
$e_{O}$--diagonal algebra.
}\end{definition}
\begin{theorem}\label{diag-restr-diag-QMC}{\rm
In the notations of Section \ref{Diag-QMC-Cl-HMP}, let $\varphi$ be the quantum Markov chain
on $\mathcal{A}=\bigotimes_{\mathbb{N}} \mathcal{B} $ determined by the pair
$(\varphi_{0},\mathcal{E})$ where $\varphi_{0}$, $\mathcal{E}$ satisfy the conditions
of Theorem \ref{thm:char-diagonaliz-CDA-NDS}.
Then the restriction of $\varphi$ to the $e_{O}$--diagonal
algebra \eqref{eO-diag-alg} is characterized by the joint probabilities
\begin{equation}\label{joint-exp-O-diag-proc}
P_{O} (j_{O_{0}}(e_{O_0;k_0,k_0}) j_{O_{1}}(e_{O_1;k_1,k_1})\cdots j_{O_{n}}(e_{O_n;k_n,k_n})
\end{equation}
$$
=\sum_{j_0,\dots, j_n\in D}p_{H;j_0}^{(0)} p_{H_{0};j_0j_1}\cdots  p_{H_{n};j_{n-1}j_n}\cdot
p(O_0=k_0|H=j_0)\cdots p(O_n=k_n|H=j_n)
$$
where $k_0,\dots, k_n\in D$,  and
\begin{equation}\label{trans-prob-OH-diag-proc-prob}
p(O_m=k_m | H_{m}=j) =  |\langle e_{H_{m};j},e_{O_m;k}\rangle|^2
\end{equation}
while the stochastic matrix $P_{H}:=(p_{H;jk})$ is given by \eqref{struct-K-NDS}.
}\end{theorem}

\noindent\textbf{Remark}.
Due to the identification \eqref{identif-calA(H,O)},
in \eqref{trans-prob-OH-diag-proc-prob}, both hidden and observable random variables take values
in the state space $D$. Moreover the notation $H_{m}=j_m$ means that the hidden random variable 
$H_{m}$ at time $m$ is in the state $j_m$, $p(O_m=k_m | H_{m}=j_{m})$ denotes the conditional
probability that, at time $m$, $O_m=k_m $ given that, at time $m$, $H_{m}=j_{m}$ and
$p_{H;j_{m-1}j_m}$ is the conditional probability that, at time $m$, $H_{m}=j_m$ given that, 
at time $m$, $H=j_{m-1}$.
These identifications allow to simplify the notations but, for the physical interpretation of 
the conditional probabilities it is better to avoid them.\\

\noindent\textbf{Remark}.
Comparing the right hand side of \eqref{joint-exp-O-diag-proc} with that of
\eqref{joint-exp-HO-discr-proc}, one immediately recognizes that the former gives the joint
probabilities of the observable process with emission operators
$(p(O_m=k_m|H=j_m))$ and Markov operators $(p_{H;j_{m-1}j_m})$.
\\

\noindent\textbf{Proof}.
Recalling the structure \eqref{CCDA-26.5} of the density operators associated to diagonalizable
conditional density amplitude $K$, one finds
\begin{equation}\label{joint-prob-diag-QMC1}
\varphi (j_0(e_{O_0;k_0,k_0}) j_1(e_{O_1;k_1,k_1})\cdots j_n(e_{O_n;k_n,k_n})
\end{equation}
$$
=\sum_{j_1,\dots, j_n\in D} p_{H;j_0}^{0} p_{H;j_0j_1}\dots p_{H;j_{n-1}j_n}
\hbox{Tr}_{[0,n]}\left((e_{H;j_0j_0}\otimes e_{H;j_1j_1}\otimes\dots \otimes e_{H;j_nj_n})\right.
$$
$$
\left.(j_0(e_{O_0;k_0,k_0}) j_1(e_{O_1;k_1,k_1})\cdots j_n(e_{O_n;k_n,k_n}))\right)
$$
$$
=\sum_{j_1,\dots, j_n\in D} p_{H;j_0}^{0} p_{H;j_0j_1}\dots p_{H;j_{n-1}j_n}
\hbox{Tr}_{[0,n]}\left(j_0(e_{H;j_0j_0})j_1(e_{H;j_1j_1})\cdots j_{n-1}(e_{H;j_{n-1}j_{n-1}})\right.
$$
$$
\left.(j_0(e_{O_0;k_0,k_0}) j_1(e_{O_1;k_1,k_1})\cdots j_{n-1}(e_{O_{n-1};k_{n-1},k_{n-1}}))
j_n(e_{H;j_nj_n}e_{O_n;k_n,k_n})\right)
$$
$$
=\sum_{j_1,\dots, j_n\in D} p_{H;j_0}^{0} p_{H;j_0j_1}\dots p_{H;j_{n-1}j_n}
\hbox{Tr}_{[0,n-1]}\left(j_0(e_{H;j_0j_0})j_1(e_{H;j_1j_1})\cdots j_{n-1}(e_{H;j_{n-1}j_{n-1}})\right.
$$
$$
\left.(j_0(e_{O_0;k_0,k_0}) j_1(e_{O_1;k_1,k_1})\cdots j_{n-1}(e_{O_{n-1};k_{n-1},k_{n-1}}))
\right) \overline{\hbox{Tr}}_{n}\left(j_n(e_{H;j_nj_n}e_{O_n;k_n,k_n}) \right)
$$
$$
=\sum_{j_1,\dots, j_n\in D} p_{H;j_0}^{0} p_{H;j_0j_1}\dots p_{H;j_{n-1}j_n}
\hbox{Tr}_{[0,n-1]}\left(j_0(e_{H;j_0j_0})j_1(e_{H;j_1j_1})\cdots j_{n-1}(e_{H;j_{n-1}j_{n-1}})\right.
$$
$$
\left.(j_0(e_{O_0;k_0,k_0}) j_1(e_{O_1;k_1,k_1})\cdots j_{n-1}(e_{O_{n-1};k_{n-1},k_{n-1}}))\right)
| \langle e_{H;j_n},e_{O_n;k_n}\rangle|^2
$$
$$
\overset{\eqref{trans-prob-OH-diag-proc-prob}}{=}
\sum_{j_1,\dots, j_n\in D} p_{H;j_0}^{0} p_{H;j_0j_1}\dots p_{H;j_{n-1}j_n}p(H=j|O_n=k_n)
$$
$$
\hbox{Tr}_{[0,n-1]}\left(j_0(e_{H;j_0j_0})j_1(e_{H;j_1j_1})\cdots j_{n-1}(e_{H;j_{n-1}j_{n-1}})\right.
$$
$$
\left.(j_0(e_{O_0;k_0,k_0}) j_1(e_{O_1;k_1,k_1})\cdots j_{n-1}(e_{O_{n-1};k_{n-1},k_{n-1}}))
\right)
$$
From this, \eqref{joint-exp-O-diag-proc} follows by induction.
$\qquad\square$

\subsection{Diagonal restrictions of quantum hidden Markov processes}\label{sec:diag-restr-QHMP}

In the notations and assumptions of Theorem \ref{th:struct-QHMP}, let
\begin{equation}\label{df-QSP-calA(H,O)-t}
(\mathcal{A}_{H,O},P_{H,O}, (\mathcal{B}_{O_{n}}\otimes \mathcal{B}_{H_{n}})_{n\in \mathbb{N}}, (j_{O_{n}}\otimes j_{H_{n}})_{n\in \mathbb{N}})
\end{equation}
be the quantum hidden process defined in Theorem \ref{th:joint-exps-(O,H)-q} with the
special choice of the emission operators
\begin{equation}\label{diag-emiss-op}
\mathcal{E}_{O,H,n} \colon \mathcal{B}\otimes \mathcal{B} \
\overset{\eqref{df-jHn-calA-Hn}, \eqref{df-jOn-calA-On}}{\equiv} \
\mathcal{B}_{H_{n}}\otimes \mathcal{B}_{O_{n}}\to\mathcal{B}\equiv \mathcal{B}_{H_{n}}
\quad,\quad n\in \mathbb{N}
\end{equation}
given by
\begin{equation}\label{df-calE-I(O,H,n)}
\mathcal{E}_{O,H,n}(x)
:= \overline{\hbox{Tr}}_{2}(K_{O',H',n}^*x K_{O',H',n})
\quad,\quad x\in \mathcal{B}_{H_{n}}\otimes \mathcal{B}_{O_{n}}
\end{equation}
where $\overline{\hbox{Tr}}_{2}$ given by \eqref{df-partial-Tr} and $K_{O',H',n}$ is defined
as follows.
For each $n\in\mathbb{N}$, fix two ortho--normal bases $e_{O'_{n}}\equiv (e_{O'_{n};j})_{j\in D}$
(resp. $e_{H'_{n}}\equiv (e_{H'_{n};j})_{j\in D}$)
of $\mathcal{H}$ as in Definition \ref{df-calD-e(on)}, and
denote, in the notation \eqref{df-eiej^*}, $e_{O'_{n};jk}:=e_{O'_{n};j}e^*_{O'_{n};k}$
(resp. $e_{H'_{n};jk} := e_{H'_{n};j}e_{H'_{n};k}^*$)
(see \eqref{df-eiej^*}) be the associated system of matrix units.
Fix a \textbf{conditional density amplitude} (with respect to the trace)
\begin{equation}\label{Comm-CDA-KOH-diag}
K_{O',H',n} := \sum_{i,j\in D} K_{O',H',n;i,j} e_{H'_{n};ii}\otimes e_{O'_{n};jj}
\in \mathcal{B}_{H_{n}}\otimes \mathcal{B}_{O_{n}}
\end{equation}
Notice the difference between \eqref{Comm-CDA-KOH-diag}, where the projections
$e_{H'_{n};ii}$, $e_{O'_{n};jj}$ belong to \textbf{different} diagonal algebras and
\eqref{struct-K-NDS}, where the projections $e_{jj}$, $e_{kk}$ belong to \textbf{the same}
diagonal algebra.\\

\noindent Notice that, for any $g\in\mathcal{B} \equiv\mathcal{B}_{H_{n}} $ and
$f\in\mathcal{B} \equiv\mathcal{B}_{O_{n}} $, one has
\begin{equation}\label{range-calE(O,H,n)}
\mathcal{E}_{O,H,n}(g\otimes f)
=\overline{\hbox{Tr}}_{2}(K_{O',H',n}^*(g\otimes f)K_{O',H',n})
\end{equation}
$$
=\overline{\hbox{Tr}}_{2}\Big(\Big(\sum_{i_{1},j_{1}\in D} K_{O',H',n;i_{1},j_{1}} e_{H'_{n};i_{1}i_{1}}\otimes e_{O'_{n};j_{1}j_{1}}\Big)^* (g\otimes f)
\Big(\sum_{i_{2},j_{2}\in D} K_{O',H',n;i_{2},j_{2}} e_{H'_{n};i_{2}i_{2}}\otimes e_{O'_{n};j_{2}j_{2}}\Big)\Big)
$$
$$
=\sum_{i_{1},j_{1},i_{2},j_{2}\in D}\overline{K_{O',H',n;i_{1},j_{1}}}K_{O',H',n;i_{2},j_{2}}
\overline{\hbox{Tr}}_{2}\left(
e_{H'_{n};i_{1}i_{1}}ge_{H'_{n};i_{2}i_{2}}\otimes
e_{O'_{n};j_{1}j_{1}}fe_{O'_{n};j_{2}j_{2}}\right)
$$
$$
=\sum_{i_{1},j_{1},i_{2},j_{2}\in D}\overline{K_{O',H',n;i_{1},j_{1}}}K_{O',H',n;i_{2},j_{2}}
\langle e_{H'_{n};i_{1}},ge_{H'_{n};i_{2}}\rangle
\langle e_{O'_{n};j_{1}}, fe_{O'_{n};j_{2}}\rangle
\overline{\hbox{Tr}}_{2}\left(e_{H'_{n};i_{1}i_{2}}\otimes e_{O'_{n};j_{1}j_{2}}\right)
$$
$$
=\sum_{i_{1},j_{1},i_{2},j_{2}\in D}\overline{K_{O',H',n;i_{1},j_{1}}}K_{O',H',n;i_{2},j_{2}}
\langle e_{H'_{n};i_{1}},ge_{H'_{n};i_{2}}\rangle
\langle e_{O'_{n};j_{1}}, fe_{O'_{n};j_{2}}\rangle \delta_{j_{1},j_{2}}
e_{H'_{n};i_{1}i_{2}}
$$
$$
=\sum_{i_{1},j_{1},i_{2}\in D}\overline{K_{O',H',n;i_{1},j_{1}}}K_{O',H',n;i_{2},j_{1}}
\langle e_{H'_{n};i_{1}},ge_{H'_{n};i_{2}}\rangle
\langle e_{O'_{n};j_{1}}, fe_{O'_{n};j_{1}}\rangle e_{H'_{n};i_{1}i_{2}}
$$
For $n\in\mathbb{N}$, denote $\hat{\mathcal{D}}_{O'_{n}}$ (resp. $\hat{\mathcal{D}}_{H'_{n}}$) the diagonal algebra generated by $(e_{O'_{n};jj})_{j\in D}$ (resp. $(e_{H'_{n};ii})_{i\in D}$),
then \eqref{df-jHn-calA-Hn}, \eqref{df-jOn-calA-On} imply that
$$
j_{O_{n}}(\hat{\mathcal{D}}_{O'_{n}})=:\mathcal{D}_{O'_{n}}\subset \mathcal{B}_{O_{n}}\ ;\quad
j_{H_{n}}(\hat{\mathcal{D}}_{H'_{n}})=: \mathcal{D}_{H'_{n}}\subset \mathcal{B}_{H_{n}}
$$
Define respectively the $H'$--diagonal algebra and the $O'$--diagonal algebra by
\begin{equation}\label{df-H'-O'-diag-alg}
\mathcal{D}_{H'} := \bigvee \mathcal{D}_{H'_{n}}
\quad;\quad \mathcal{D}_{O'} := \bigvee \mathcal{D}_{O'_{n}}
\end{equation}
An important consequence of \eqref{range-calE(O,H,n)} is that $\mathcal{E}_{O,H,n}$ maps the diagonal algebra
$\hat{\mathcal{D}}_{H'_{n}}\otimes\hat{\mathcal{D}}_{O'_{n}}\subset \mathcal{B}\otimes\mathcal{B}$ into the diagonal algebra $\hat{\mathcal{D}}_{H'_{n}}\subset \mathcal{B}$.
In particular, putting $f=g=1_{\mathcal{B}}$ in \eqref{range-calE(O,H,n)}, one finds
\begin{align*}
1_{\mathcal{B}}\overset{\text{Le. \ref{|K(jk)|2-stoch-matr}}}{=}&
\mathcal{E}_{O,H,n}(1_{\mathcal{B}}\otimes 1_{\mathcal{B}})\\
\overset{\eqref{range-calE(O,H,n)}}{=}&
\sum_{i_{1},j_{1},i_{2}\in D}\overline{K_{O',H',n;i_{1}, j_{1}}}K_{O',H',n;i_{2},j_{1}}
\delta_{i_{1},i_{2}} e_{H'_{n};i_{1}i_{2}}\\
=&\sum_{i_{1}\in D}\Big(\sum_{j_{1}\in D} |K_{O',H',n;i_{1},j_{1}}|^2\Big)e_{H'_{n};i_{1}i_{1}}
\end{align*}
In other words
\begin{equation}\label{df-P(On=j|H'n=i)}
P(O_{n}=j | H'_{n}=i):= p_{O,H',n;i,j} := |K_{O',H',n;i,j}|^2 \ ,\quad\forall i, j \in D
\end{equation}
is a stochastic matrix.
\begin{theorem}\label{diag-restr-diag-QMC-t}{\rm
In the notations introduced above and in the assumptions of Theorem \ref{diag-restr-diag-QMC},
the joint probabilities of the quantum hidden Markov\\ $(H,O)$--process with emission operators
given by \eqref{df-calE-I(O,H,n)} and \eqref{Comm-CDA-KOH-diag}, restricted to the diagonal algebra
$\mathcal{D}_{H'} \otimes \mathcal{D}_{O'} $, are given by
\begin{align}\label{joint-prob-H'O'-diag-proc}
&P_{H',O'}\Big(\prod_{m=0}^{n}j_{H_{m}}(e_{H';j_m,j_m}) \prod_{m=0}^{n}j_{O'_{m}}(e_{O'_m;k_m,k_m})\Big)\\
=& p_{O',H,0;h_{0},k_{0}} p_{O',H,1;h_{1},k_{1}} \cdots p_{O',H,n-1;h_{n-1},k_{n-1}}p_{O',H,n;h_{n},k_{n}}\notag\\
&\sum_{j_0,\dots, j_n\in D}p_{H;j_0}^{(0)} p_{H;j_0j_1}\cdots  p_{H;j_{n-1}j_n}\cdot|\langle e_{H_0;j_0},e_{H_0';k_0}\rangle|^2 \cdots |\langle e_{H_n;j_n},e_{H_n';k_n}\rangle|^2\notag
\end{align}
where
$p_{H;jj'}$ and $p_{O,H',n;i,j} $ are given by \eqref{Comm-CDA-HM} and \eqref{df-P(On=j|H'n=i)} respectively.
}\end{theorem}
\textbf{Remark}. Notice that these joint probabilities depend on the transition probabilities
$|\langle e_{H_r;j_r},e_{H_r';k_r}\rangle|^2$, but not on the
$|\langle e_{O_r;j_r},e_{O_r';k_r}\rangle|^2$.\\

\noindent\textbf{Proof}.
From Theorem \eqref{th:struct-QHMP} we know that, under the above assumptions, the process \eqref{df-QSP-calA(H,O)-t} is the hidden quantum Markov process whose joint
expectations are given by
\begin{align}\label{joint-exp-qHOP-diag}
&P_{H,O}\left(\prod_{m=0}^{n}(j_{H_{m}}(f_{m})j_{O_{m}}(g_{m}) \right)\\
= &P_{H_{0}}\left(\mathcal{E}_{H_{0}}\left(\mathcal{E}_{O,H;0} (g_{0}\otimes f_{0})\otimes\mathcal{E}_{H;1} \mathcal{E}_{O,H;1} \left((g_{1}\otimes f_{1})\otimes \cdots\right.\right.\right. \notag\\
& \quad\ \left.\left.\left.\otimes \mathcal{E}_{H_{n-1}} \left(\mathcal{E}_{O,H;n}(g_{n-1}\otimes f_{n-1})\otimes
\mathcal{E}_{H_{n}}\left(\mathcal{E}_{O,H;n}\left(g_{n}\otimes f_{n}\right)\otimes 1_{H_{n+1}}
\right)\right)\right)\right)\right)\notag
\end{align}
for all $n\in\mathbb{N}$, $m\in \{0, \dots, n\}$, $g_{m}\in \mathcal{B}_{O_{n}}$ and
$f_{m}\in \mathcal{B}_{H_{n}}$.
First we calculate the expression
\begin{equation}\label{joint-exp-qHOP-diag1}
\mathcal{E}_{O,H;n}(g_{n}\otimes f_{n})
\overset{\eqref{df-calE-I(O,H,n)}}{=}
\overline{\hbox{Tr}}_{2}(K_{O',H',n}^*(g_{n}\otimes f_{n}) K_{O',H',n})
\end{equation}
$$
\overset{\eqref{range-calE(O,H,n)}}{=}
\overline{\hbox{Tr}}_{2}\Big(\Big(\sum_{i,j\in D} K_{O',H',n;i,j} e_{H'_{n};ii}\otimes e_{O'_{n};jj}\Big)^*
(g_{n}\otimes f_{n})\Big(\sum_{i,j\in D} K_{O',H',n;i,j} e_{H'_{n};ii}\otimes e_{O'_{n};jj}\Big)\Big)
$$
$$
=\sum_{i_{1},j_{1},i_{2},j_{2}\in D}   \overline{K_{O',H',n;i_{1},j_{1}}}K_{O',H',n;i_{2},j_{2}}
\overline{\hbox{Tr}}_{2}\left((e_{H'_{n};i_{1}i_{1}}\otimes e_{O'_{n};j_{1}j_{1}})(g_{n}\otimes f_{n}\right)
(e_{H'_{n};i_{2}i_{2}}\otimes e_{O'_{n};j_{2}j_{2}}))
$$
$$
=\sum_{i_{1},j_{1},i_{2},j_{2}\in D}   \overline{K_{O',H',n;i_{1},j_{1}}}K_{O',H',n;i_{2},j_{2}}
\overline{\hbox{Tr}}_{2}
\left(e_{H'_{n};i_{1}i_{1}}g_{n}e_{H'_{n};i_{2}i_{2}}
\otimes e_{O'_{n};j_{1}j_{1}}f_{n}e_{O'_{n};j_{2}j_{2}}\right)
$$
$$
=\sum_{i_{1},j_{1},i_{2},j_{2}\in D}   \overline{K_{O',H',n;i_{1},j_{1}}}K_{O',H',n;i_{2},j_{2}}
\langle e_{H'_{n};i_{1}},g_{n}e_{H'_{n};i_{2}} \rangle
\langle e_{O'_{n};j_{1}},f_{n}e_{O'_{n};j_{2}} \rangle
\overline{\hbox{Tr}}_{2}
\left(e_{H'_{n};i_{1}i_{2}}
\otimes e_{O'_{n};j_{1}j_{2}}\right)
$$
$$
=\sum_{i_{1},j_{1},i_{2},j_{2}\in D}   \overline{K_{O',H',n;i_{1},j_{1}}}K_{O',H',n;i_{2},j_{2}}
\langle e_{H'_{n};i_{1}},g_{n}e_{H'_{n};i_{2}} \rangle
\langle e_{O'_{n};j_{1}},f_{n}e_{O'_{n};j_{2}} \rangle
e_{H'_{n};i_{1}i_{2}}\delta_{j_{1},j_{2}}
$$
$$
=\sum_{i_{1},j_{1},i_{2}\in D}   \overline{K_{O',H',n;i_{1},j_{1}}}K_{O',H',n;i_{2},j_{1}}
\langle e_{H'_{n};i_{1}},g_{n}e_{H'_{n};i_{2}} \rangle
\langle e_{O'_{n};j_{1}},f_{n}e_{O'_{n};j_{1}} \rangle
e_{H'_{n};i_{1}i_{2}}
$$
In order to calculate the expectation value \eqref{joint-exp-qHOP-diag}, one has to replace
$f_{n}$ by a generic $e_{O_n;a,b}$ and $g_{n}$ by a generic $e_{H';a',b'}$ ($a,b,a',b'\in D$) in \eqref{joint-exp-qHOP-diag1}. This leads to
\begin{align}\label{joint-exp-qHOP-diag2}
&\mathcal{E}_{O,H;n}(e_{H_n';a',b'}\otimes e_{O_n';a,b}) \\
=&\sum_{i_{1},j_{1},i_{2}\in D}\overline{K_{O',H',n;i_{1},j_{1}}} K_{O',H',n;i_{2},j_{1}}\notag\\
&\qquad\qquad \langle e_{H'_{n};i_{1}},e_{H';a',b'} e_{H'_{n};i_{2}} \rangle\langle e_{O'_{n};j_{1}}, e_{O_n;a,b} e_{O'_{n};j_{1}} \rangle e_{H'_{n};i_{1}i_{2}}\notag\\
=&\sum_{i_{1},j_{1},i_{2}\in D}   \overline{K_{O',H',n;i_{1},j_{1}}}K_{O',H',n;i_{2},j_{1}}
\delta_{i_{1},a'}\delta_{b',i_{2}}\delta_{j_1,b}\delta_{a,j_{1}}
e_{H'_{n};i_{1}i_{2}}\notag\\
=&\delta_{a,b} \overline{K_{O',H',n;a',a}}K_{O',H',n;b',a} e_{H'_{n};a',b'}\notag
\end{align}
\eqref{joint-exp-qHOP-diag2} implies in particular that $\mathcal{E}_{O,H;n}$ maps
the diagonal algebra $\mathcal{D}_{e_{H_{n}'}}\otimes\mathcal{D}_{e_{O_{n}'}}$ into
$\mathcal{D}_{e_{H_{n}'}}\subset \mathcal{D}_{H}$ and the elements not in the algebra
$\mathcal{B}_{H_{n}'}\otimes\mathcal{D}_{e_{O_{n}'}}$ to zero.
Moreover, putting in \eqref{joint-exp-qHOP-diag2} $a'=b'=h_{m}$ and $a=b=k_{m}$, one finds
\begin{align}\label{joint-exp-qHOP-diag3}
&\mathcal{E}_{O,H;n}(e_{H_{n}';h_{n},h_{n}}\otimes e_{O_n';k_{n},k_{n}}) \\
=&\sum_{i_{1},j_{1},i_{2}\in D}  \overline{K_{O',H',n;i_{1}, j_{1}}}K_{O',H',n;i_{2},j_{1}}\notag\\
&\qquad\qquad\langle e_{H'_{n};i_{1}}, e_{H';h_{n},h_{n}}e_{H'_{n};i_{2}} \rangle\langle e_{O'_{n}; j_{1}}, e_{O_n';k_{n},k_{n}}e_{O'_{n};j_{1}} \rangle
e_{H'_{n};i_{1}i_{2}}\notag\\
=&\sum_{i_{1},j_{1},i_{2}\in D}   \overline{K_{O',H',n;i_{1}, j_{1}}}K_{O',H',n;i_{2},j_{1}}\delta_{i_{1},h_{n}}\delta_{h_{n},i_{2}}\delta_{k_{n},j_{1}} e_{H'_{n};i_{1}i_{2}}\notag\\
=&\overline{K_{O',H',n;h_{n},k_{n}}}K_{O',H',n;h_{n},k_{n}} e_{H'_{n};h_{n},h_{n}} \notag\\
=&|K_{O',H',n;h_{n},k_{n}}|^2 e_{H'_{n};h_{n},h_{n}} =: p_{O,H,n;h_{n},k_{n}}  e_{H'_{n};h_{n},h_{n}}\notag\\
\end{align}
In conclusion:
\begin{equation}\label{joint-exp-qHOP-diag4}
\mathcal{E}_{O,H;n}(e_{H_{n}';h_{n},h_{n}}\otimes e_{O_n';k_{n},k_{n}})
\overset{\eqref{df-P(On=j|H'n=i)}}{=} p_{O,H,n;h_{n},k_{n}}  e_{H'_{n};h_{n},h_{n}}
\end{equation}
From the above discussion it follows that the expectation value \eqref{joint-prob-H'O'-diag-proc}
is equal to
\begin{equation}\label{joint-exp-HO-diag-proc-t1}
P_{H,O}\left(\prod_{m=0}^{n}j_{H,m}(e_{H';j_m,j_m}) j_{O,m}(e_{O_m';k_m,k_m})\right)
\end{equation}
$$
= P_{H_{0}}\left(\mathcal{E}_{H_{0}}\left(\mathcal{E}_{O,H;0} (e_{H';j_0,j_0}\otimes e_{O_0';k_0k_0})\otimes
\mathcal{E}_{H_1}(\mathcal{E}_{O,H;1}\left((e_{H';j_1,j_1}\otimes e_{O_1';k_1,k_1})\otimes \cdots\right.\right.\right.
$$
$$
\left.\left.\left.\otimes \mathcal{E}_{H_{n-1}} \left(\mathcal{E}_{O,H;n}(e_{H';j_{n-1},j_{n-1}}
\otimes e_{O_{n-1}';k_{n-1},k_{n-1}})\otimes
\mathcal{E}_{H_{n}}\left(\mathcal{E}_{O,H;n}
\left(e_{H';j_{n},j_{n}}\otimes
e_{O_{n}';k_{n},k_{n}}\right)\otimes 1_{H_{n+1}}
\right)\right)\right)\right)\right)
$$
$$
= P_{H_{0}}\left(\mathcal{E}_{H_{0}}
\left((p_{O,H,0;h_{0},k_{0}}  e_{H'_{0};h_{0},h_{0}})\otimes
\mathcal{E}_{H_1}((p_{O,H,1;h_{1},k_{1}}  e_{H'_{1};h_{1},h_{1}})\otimes \cdots \right.\right.
$$
$$
\left.\left.\left.\otimes \mathcal{E}_{H_{n-1}}
\left((p_{O,H,n-1;h_{n-1},k_{n-1}}  e_{H'_{n-1};h_{n-1},h_{n-1}})\otimes
\mathcal{E}_{H_{n}}\left((p_{O,H,n;h_{n},k_{n}}  e_{H'_{n};h_{n},h_{n}})\right)
\otimes 1_{H_{n+1}}
\right)\right)\right)\right)
$$
$$
= p_{O,H,0;h_{0},k_{0}} p_{O,H,1;h_{1},k_{1}} \cdots p_{O,H,n-1;h_{n-1},k_{n-1}}p_{O,H,n;h_{n},k_{n}}
$$
$$
P_{H_{0}}\left(\mathcal{E}_{H_{0}}\left(e_{H'_{0};h_{0},h_{0}} \otimes \mathcal{E}_{H_1}(e_{H'_{1};h_{1},h_{1}}\otimes \cdots
\otimes \, \mathcal{E}_{H_{n-1}}\left(e_{H'_{n-1}; h_{n-1},h_{n-1}}\otimes\mathcal{E}_{H_{n}}\left(e_{H'_{n};h_{n},h_{n}}\otimes 1_{H_{n+1}}\right)\right)\right)\right)
$$
$$
= p_{O,H,0;h_{0},k_{0}} p_{O,H,1;h_{1},k_{1}} \cdots p_{O,H,n-1;h_{n-1},k_{n-1}}p_{O,H,n;h_{n},k_{n}}
P_{H}\left(\prod_{m=0}^{n} j_{m}(e_{H'_{0};h_{m},h_{m}})\right)
$$
where the probabilities $p_{O,H,m;h_{m},k_{m}}$ are given by \eqref{joint-exp-qHOP-diag4} and
$P_{H}\left(\prod_{m=0}^{n} j_{m}(e_{H'_{0};h_{m},h_{m}})\right)$ are the joint probabilities
of the restriction of the diagonalizable quantum Markov chain on the $H'$--diagonal algebra.\\
These joint probabilities have been calculated in Theorem \ref{diag-restr-diag-QMC}, where
the\\ $O$--process plays the role of the $H'$--process, and are given by
\begin{equation}\label{joint-exp-O-diag-proc-t}
P_{H}\left(\prod_{m=0}^{n} j_{m}(e_{H'_{0};k_{m},k_{m}})\right)
\end{equation}
$$
=\sum_{j_0,\dots, j_n\in D}p_{H;j_0}^{(0)} p_{H;j_0j_1}\cdots  p_{H;j_{n-1}j_n}\cdot
p(H'_0=k_0|H=j_0)\cdots p(H'_n=k_n|H=j_n)
$$
where the probabilities $p_{H;j}^{(0)}$ and $p_{H;jj'}$ are given by \eqref{CCDA-26.5} and
\begin{equation}\label{trans-prob-OH-diag-proc-t}
p(H'_m=k_m | H=j) =  |\langle e_{H;j},e_{H'_m;k}\rangle|^2
\end{equation}
Combining \eqref{joint-exp-HO-diag-proc-t1}, \eqref{joint-exp-O-diag-proc-t} and \eqref{trans-prob-OH-diag-proc-t}, one obtains
\begin{equation}\label{joint-exp-HO-diag-proc-t2}
P_{H,O} \left(\prod_{m=0}^{n}j_{H,m}(e_{H';j_m,j_m})j_{O,m}(e_{O_m;k_m,k_m})\right)
\end{equation}
$$
= p_{O,H,0;h_{0},k_{0}} p_{O,H,1;h_{1},k_{1}} \cdots p_{O,H,n-1;h_{n-1},k_{n-1}}p_{O,H,n;h_{n},k_{n}}
$$
$$
\sum_{j_0,\dots, j_n\in D}p_{H;j_0}^{(0)} p_{H;j_0j_1}\cdots  p_{H;j_{n-1}j_n}\cdot
|\langle e_{H;j_0},e_{H'_m;k_0}\rangle|^2 \cdots |\langle e_{H;j_n},e_{H'_m;k_n}\rangle|^2
$$
which is \eqref{joint-prob-H'O'-diag-proc}.
$\qquad\square$

\subsubsection{A family of hidden, but not hidden Markov, classical processes}
\label{sec:fam-hidd-not-HMP-class}

Recalling the notations: $\forall h_{n}, k_{n} \in D$,
\begin{equation}\label{df-P(On=kn|H'n=hn)}
P(O_{n}=k_{n} | H'_{n}=h_{n})
:= p_{O,H',n;h_{n},k_{n}} := |K_{O,H',n;h_{n},k_{n}}|^2
\end{equation}
\begin{equation}\label{df-P(H'n=kn|Hn=jn)}
P(H'_{n}=k_{n} | H_{n}=j_{n}) = |\langle e_{H;j_n},e_{H'_m;k_n}\rangle|^2
= P(H_{n}=j_{n} | H'_{n}=k_{n})
\end{equation}
one has, from Theorem \eqref{diag-restr-diag-QMC-t} that the joint probabilities of the $(H,O)$--process, restricted to the diagonal algebra $\mathcal{D}_{H'} \otimes \mathcal{D}_{O'}$,
are given by
\begin{equation}\label{joint-exp-HO-diag-proc-tt}
P_{H,O} \left(\prod_{m=0}^{n}j_{H_{m}}(e_{H';j_m,j_m})\prod_{m=0}^{n}j_{O_{m}}(e_{O_m';k_m,k_m})\right)
\end{equation}
$$
= p_{O,H,0;h_{0},k_{0}} p_{O,H,1;h_{1},k_{1}} \cdots p_{O,H,n-1;h_{n-1},k_{n-1}}p_{O,H,n;h_{n},k_{n}}
$$
$$
\sum_{j_0,\dots, j_n\in D}p_{H;j_0}^{(0)} p_{H;j_0j_1}\cdots  p_{H;j_{n-1}j_n}\cdot
|\langle e_{H;j_0},e_{H'_m;k_0}\rangle|^2 \cdots |\langle e_{H;j_n},e_{H_m';k_n}\rangle|^2
$$
$$
\overset{\eqref{df-P(On=kn|H'n=hn)}, \eqref{df-P(H'n=kn|Hn=jn)}}{=} \
\prod_{m=0}^{n}P(O_{m}=k_{m} | H'_{m}=h_{m})
$$
$$
\sum_{j_0,\dots, j_n\in D}
p_{H;j_0}^{(0)} p_{H;j_0j_1}\cdots  p_{H;j_{n-1}j_n}\cdot
\prod_{m=0}^{n} P(H'_{m}=k_{m} | H_{m}=j_{m})
$$
We know (see \eqref{cond-exps-O|H}) that the condition
$$
P_{H,H'} \left(H'_{0}=k_{0};H'_{1}=k_{1};\dots;H'_{n}=k_{n} \Big|
H_{0}=j_{0};H_{1}=j_{1};\dots;H_{n}=j_{n} \right)
$$
\begin{equation}\label{H'-obs-H-hidd}
= \prod_{m=0}^{n} P(H'_{m}=k_{m} | H_{m}=j_{m})
\end{equation}
characterizes the process $(H'_{m})$ as observable process of the underlying Markov process
$(H_{m})$ and that, in this case, one has
$$
\sum_{j_0,\dots, j_n\in D}
p_{H;j_0}^{(0)} p_{H;j_0j_1}\cdots  p_{H;j_{n-1}j_n}\cdot
\prod_{m=0}^{n} P(H'_{m}=k_{m} | H_{m}=j_{m})
$$
$$
=P_{H'} \left(H'_{0}=k_{0};H'_{1}=k_{1};\dots;H'_{n}=k_{n}\right)
$$
where the right hand side represents the joint probabilities of the $H'$--process.
In view of this, \eqref{joint-exp-HO-diag-proc-tt} can be written in the form
\begin{equation}\label{joint-exp-OH'-diag-proc}
P_{H,O} \left(\prod_{m=0}^{n}j_{H_{m}}(e_{H';j_m,j_m})\prod_{m=0}^{n}j_{O_{m}}(e_{O_m;k_m,k_m})\right)
\end{equation}
$$
=\prod_{m=0}^{n}P(O_{m}=k_{m} | H'_{m}=h_{m})
P_{H'} \left(H'_{0}=k_{0};H'_{1}=k_{1};\dots;H'_{n}=k_{n}\right)
$$
In analogy with \eqref{H'-obs-H-hidd}, we introduce the condition
\begin{equation}\label{O-obs-H'-hidd}
\prod_{m=0}^{n}P(O_{m}=k_{m} | H'_{m}=h_{m})
\end{equation}
$$
= P_{O,H'} \left(O_{0}=k_{0};O_{1}=k_{1};\dots O'_{n}=k_{n} \Big|
H'_{0}=h_{0};H'_{1}=h_{1};\dots;H'_{n}=h_{n} \right)
$$
that characterizes the process $(O_{m})$ as observable process of the underlying, not Markov
but hidden Markov, process $(H'_{m})$. Assuming that \eqref{O-obs-H'-hidd} holds,
in view of \eqref{O-obs-H'-hidd}, the right hand side of \eqref{joint-exp-OH'-diag-proc} becomes
$$
\prod_{m=0}^{n}P(O_{m}=k_{m} | H'_{m}=h_{m})
P_{H'} \left(H'_{0}=k_{0};H'_{1}=k_{1};\dots;H'_{n}=k_{n}\right)
$$
$$
=P_{O,H'} \left(O_{0}=k_{0};O_{1}=k_{1};\dots O'_{n}=k_{n} \Big|
H'_{0}=h_{0};H'_{1}=h_{1};\dots;H'_{n}=h_{n} \right)
$$
$$
P_{H'} \left(H'_{0}=k_{0};H'_{1}=k_{1};\dots;H'_{n}=k_{n}\right)
$$
$$
=P_{O,H'} \left(O_{0}=k_{0};O_{1}=k_{1};\dots O'_{n}=k_{n} \ ; \
H'_{0}=h_{0};H'_{1}=h_{1};\dots;H'_{n}=h_{n} \right)
$$
which is precisely the classical probabilistic interpretation of the left hand side of
\eqref{joint-exp-HO-diag-proc-tt}, i.e.
$$
P_{H,O} \left(\prod_{m=0}^{n}j_{H_{m}}(e_{H';j_m,j_m})\prod_{m=0}^{n}j_{O_{m}}(e_{O_m;k_m,k_m})\right)
$$
In conclusion: the restriction, on the $(H',O)$--diagonal algebra, of the quantum hidden
Markov $(H,O)$--process with emission operators given by \eqref{diag-emiss-op} and \eqref{Comm-CDA-KOH-diag}, produce a new family of classical stochastic processes.
In a process of this family the pair $(H,O)$ of hidden and observable process, in a usual hidden
Markov process, is replaced by a \textbf{triple} of processes $(H,H',O)$ such that:\\
(i) $H'$ is the observable process of the Markov process $H$;\\
(ii) $O$ is the observable process of the hidden Markov process $H'$.\
For this reason, a process in this family could be called a $3$--\textbf{tier hidden Markov
process}.\\

\noindent In particular, processes in this family produce a first class of examples of the
processes abstractly defined in Definition \ref{df:class-HP}, namely: hidden processes
whose underlying process is not Markov (hidden Markov in the family in question).

\subsection{ $D_{e}$--\textbf{preserving} quantum Markov chains and associated classical
hidden processes}\label{sec:De-pres-QMC-HP}

The family of $e$--\textbf{diagonalizable} quantum Markov chains (see Definition \ref{df:diagonalizable-QMC} the end of Section \ref{sec:QMC}) is a particular class of
homogeneous backward QMC whose transition expectation $\mathcal{E}$ enjoys the following propery:
\begin{equation}\label{calE(DexDe)-subseteq-De}
\mathcal{E}(D_{e}\otimes D_{e})\subseteq D_{e}\otimes 1_{e} \subseteq D_{e}\otimes D_{e}
\end{equation}
The following problem:\\
\textit{ Given a $C^*$--algebra $\mathcal{B}$, characterize the transition expectations  $\mathcal{E}\colon \mathcal{B}\otimes \mathcal{B} \to \mathcal{B}
\equiv \mathcal{B}\otimes 1_{\mathcal{B}}$ such that there exists an abelian sub--algebra
$\mathcal{D}\subset\mathcal{B}\otimes \mathcal{B}$ satisfying
\begin{equation}\label{calD-pres-trans-exp}
\mathcal{E}(\mathcal{D})\subseteq \mathcal{D}\cap (\mathcal{B}\otimes 1_{\mathcal{B}})
\end{equation}
has accompanied the development of quantum Markov chain since the its beginning}.\\
In its general formulation, this problem is still open even if
$\mathcal{B}=\mathcal{B}(\mathcal{H})$ for some Hilbert space $\mathcal{H}$.
However, in \cite{LuYG95-QMC}, the author solved this problem for a particular but important
family of transition expectations $\mathcal{E}\colon M_{d}(\mathbb{C})\otimes M_{d}(\mathbb{C})$,
namely those of the form
\begin{equation}\label{a(1.2)}
\mathcal{E}(x)=\overline{\hbox{Tr}}_2(\sum_{j\in D_{\mathcal{E}}}K^*_{j}\ x\ K_{j})
\end{equation}
and such that there exists a maximal abelian sub--algebra $D_{e}\subset M_{d}(\mathbb{C})$
satisfying \eqref{calE(DexDe)-subseteq-De}.
\begin{definition}\label{2.5}{\rm
A transition expectation $\mathcal{E}\colon M_{d}(\mathbb{C})\otimes M_{d}(\mathbb{C})$ satisfying condition \eqref{calE(DexDe)-subseteq-De} will be called $D_{e}$--\textbf{preserving}.
}\end{definition}
In the following section we give an improved version of main result of \cite{LuYG95-QMC} and,
in Section \ref{sec:Diag-restr-De-pres-QMC} we prove that their restrictions to the $e$--diagonal
sub--algebra lead to a new class of classical hidden processes that is not included in Definition \ref{df:class-HP}, but is included in Definition \ref{df:class-gen-hidd-proc}.

\subsection{ Characterization of the transition expectations on $M_d(\mathbb{C})$ satisfying
$\mathcal{E}(D_{e}\otimes D_{e}) \subseteq D_{e}$}

In this section \textbf{we restrict ourselves} to the case where $\hbox{dim}(\mathcal{H})=d<+\infty$.
The assignment of the system of matrix units $(e_{h,k})_{h,k\in D}$ allows to identify
$\mathcal{B}=\mathcal{B}(H)$ with the algebra $M_d(\mathbb{C})$ of all $d\times d$ complex matrices.
In the following, using this identification, we write $M_d(\mathbb{C})$ instead of $\mathcal{B}$.
\begin{proposition}\label{2.1}{\rm
Let $\mathcal{E}$ be a transition expectation on $M_d(\mathbb{C})$ of the form \eqref{a(1.2)}.
For \textbf{any} system $e\equiv (e_{hh'})_{h,h'\in D}$ of matrix units $M_d(\mathbb{C})$,
it is always possible to write the operators $K_{r}$ in \eqref{a(1.2)} in the form
\begin{equation}\label{a(1.4)}
K_{r}= \ \sum_{h,h'\in D}e_{hh'}\otimes\ K_{r,h,h'}   \quad,\quad\forall r\in D_{\mathcal{E}}
\end{equation}
Having fixed the decomposition \eqref{a(1.4)}, $\mathcal{E}$ acts on the system
of $e$--matrix units as follows:
\begin{equation}\label{act-calE-M-M}
\mathcal{E}(e_{i,m}\otimes e_{i'm'})
=\sum_{h,j'\in D}e_{hj'} \
\sum_{r\in D_{\mathcal{E}}}\hbox{Tr}\left(K_{r;hi}^*\ e_{i'm'} \ K_{r;mj'}\right)
\end{equation}
$\mathcal{E}$ is $D_{e}$--\textbf{preserving} if and only if, for any $ m, m', j, j'\in D$ such that
$ j\ne j' $
\begin{equation}\label{(3.5)1}
\sum_{r\in D_{\mathcal{E}}}\hbox{Tr}\left(K_{r;mj}^*\ e_{m'm'} \ K_{r;mj'}\right)
=\langle e_{m'}, \Big(\sum_{r\in D_{\mathcal{E}}} K_{r,m,j'}K^*_{r,m,j}\Big) e_{m'}\rangle
= 0
\end{equation}
In this case
\begin{equation}\label{(3.4)gen2}
\mathcal{E}(e_{mm}\otimes e_{m'm'})
=\sum_{j\in D}e_{j,j} Tr(\sum_{r\in D_{\mathcal{E}}}K^*_{r,m,j} e_{m'm'} K_{r,m,j})
\end{equation}
$\mathcal{E}$ is identity preserving if and only if
\begin{equation}\label{calE(1)=1-diag}
1 = Tr\left(\sum_{r\in D_{\mathcal{E}}}\sum_{m\in D} |K_{r,m,j}|^2\right)
\quad,\quad\forall j\in D
\end{equation}
i.e. if and only if, for each $r\in D$, the matrix  $P\equiv (P_{j,m})$
$$
P_{j,m} :=  Tr\left(\sum_{r\in D_{\mathcal{E}}} |K_{r,m,j}|^2\right)  \quad,\quad m,j\in D
$$
is a \textbf{stochastic matrix}
($P_{j,m}\ge 0 \, ; \, \sum_{m\in D}P_{j,m}=1$).
}\end{proposition}
$$
\iff Tr(\sum_{r\in D_{\mathcal{E}}}\sum_{m\in D} K^*_{r,m,j}K_{r,m,j}) = 1
\quad,\quad\forall j\in D
$$
\textbf{Proof}.
For $i, m\in D$ and $b\in M_d(\mathbb{C})$, one has
$$
\mathcal{E}(e_{i,m}\otimes b)
\overset{\eqref{df-cal-E-spec1}}{=}
\sum_{r\in D_{\mathcal{E}}}\sum_{h,j,i',j'\in D}
e_{jh}e_{i,m}e_{i'j'} \ \hbox{Tr}\left(K_{r;hj}^*\ b \ K_{r;i'j'}\right)
$$
$$
=\sum_{r\in D_{\mathcal{E}}}\sum_{h,j,i',j'\in D}
\delta_{hi}\delta_{mi'}e_{hj'} \ \hbox{Tr}\left(K_{r;hj}^*\ b \ K_{r;i'j'}\right)
$$
\begin{equation}\label{calE(e(i,m)-otimes-b)}
=\sum_{r\in D_{\mathcal{E}}}\sum_{j,j'\in D}
e_{jj'} \ \hbox{Tr}\left(K_{r;ij}^*\ b \ K_{r;mj'}\right)
\end{equation}
\eqref{calE(e(i,m)-otimes-b)} shows that, for any $e_{i,m},\ e_{i'm'}$,
$$
\mathcal{E}(e_{i,m}\otimes e_{i'm'})
=\sum_{j,j'\in D}e_{jj'} \
\sum_{r\in D_{\mathcal{E}}}\hbox{Tr}\left(K_{r;ij}^*\ e_{i'm'} \ K_{r;mj'}\right)
$$
and this proves \eqref{act-calE-M-M}.
Taking $i=m$ and $i'=m'$, one gets
$$
\mathcal{E}(e_{mm}\otimes e_{m'm'})
=\sum_{j,j'\in D}e_{jj'} \
\sum_{r\in D_{\mathcal{E}}}\hbox{Tr}\left(K_{r;mj}^*\ e_{m'm'} \ K_{r;mj'}\right)
$$
The linear independence of the $(e_{hj'})$ implies that  $\mathcal{E}$ maps $D_{e} \otimes D_{e}$ into $D_{e}$, i.e. is $D_{e}$--\textbf{preserving}, if and only if
$$
j\ne j' \ \Rightarrow \
\sum_{r\in D_{\mathcal{E}}}\hbox{Tr}\left(K_{r;mj}^*\ e_{m'm'} \ K_{r;mj'}\right) = 0
\quad, \quad \forall m, m'\in D
$$
which is equivalent to \eqref{(3.5)1}.
In this case $\mathcal{E}$ is identity preserving iff
$$
1_{\mathcal{B}}
=\sum_{m,m'\in D}\mathcal{E}(e_{mm}\otimes e_{m'm'})
=\sum_{m,m'\in D}\sum_{j\in D}e_{j,j} Tr(\sum_{r\in D_{\mathcal{E}}}K^*_{r,m,j} e_{m'm'} K_{r,m,j})
$$
$$
=\sum_{j\in D}e_{j,j} Tr\left(\sum_{r\in D_{\mathcal{E}}}\sum_{m\in D} |K_{r,m,j}|^2\right)
$$
$$
\iff 1 = \sum_{m\in D} Tr\left(\sum_{r\in D_{\mathcal{E}}}|K_{r,m,j}|^2\right)
\quad,\quad\forall j\in D
$$
which is \eqref{calE(1)=1-diag}.
$\qquad\square$\\

\noindent The following theorem gives the structure of the restriction on the $e$--diagonal algebra
of a quantum Markov chain generated by a $D_{e}$--preserving transition expectation.
\begin{theorem}\label{2.3}{\rm
Let $\mathcal{E}$ be a transition expectation from $\mathcal{B} \otimes \mathcal{B}$ into
$\mathcal{B}$ mapping $D_{e} \otimes D_{e}$ into $D_{e}$.
Denote, for all $r,m,i\in \{1,\cdots ,d\}$,
\begin{equation}\label{(3.11b)}
P_{r;m,i}:= \langle e_i, \Big(\sum_{j\in D_{\mathcal{E}}} K_{j,m,r}K^*_{j,m,r}\Big)e_i\rangle
\end{equation}
\begin{equation}\label{(3.11b)2}
P_{r;m}:=\sum_{i\in D}P_{r;m,i}
\end{equation}
Then $P\equiv (P_{r;m,i})$ is a $(d\times d^2)$--stochastic matrix, namely:
\begin{equation}\label{(2.2a)}
P_{r;m,i}\ge 0\quad,\quad \forall m,r,i\in \{1,\cdots, d\}
\end{equation}
\begin{equation}\label{(2.2b)}
\sum_{m,i}P_{r;m,i}=\ 1\quad,\quad \forall r\in \{1,\cdots, d\}
\end{equation}
Moreover, if $\varphi \equiv (\varphi_0, \mathcal{E})$ is the homogeneous quantum Markov chain
on $\bigotimes_{\mathbb{N}} M_d(\mathbb{C})$ determined by the pair $(\varphi_0, \mathcal{E})$
where $\varphi_0$ (initial state) is given by
\begin{equation}\label{e-diag-w(0)}
\varphi_0 := Tr(\sum_{i\in D}p^{(0)}_{i}e_{i,i} \, \cdot \, )
:= Tr(p^{(0)} \, \cdot \, )
\quad;\quad p^{(0)}_{i}\ge 0 \ , \  \sum_{i\in D} p^{(0)}_{i} = 1
\end{equation}
Then one has
\begin{align}\label{cal-E0]-cal-E-De-pres}
&\mathcal{E}(e_{i_0,i_0}\otimes\mathcal{E}(e_{i_1,i_1}\otimes
\cdots \otimes \mathcal{E}(e_{i_{n-1}i_{n-1}}\otimes
\mathcal{E}(e_{i_{n},i_{n}}\otimes 1))\cdots ))\\
=&\sum_{r_{n}, r_{n-1},\dots, r_{1}, r_0}
P_{r_0;i_0,r_{1}}P_{r_{1};i_{1},r_{2}} \cdots P_{r_{n-1};i_{n-1},r_{n}}P_{r_{n};i_{n}}e_{r_0r_0}\notag
\end{align}
In particular
\begin{align}\label{e-diag-proc}
&\varphi (e_{i_0,i_0}\otimes e_{i_1,i_1}\otimes \cdots \otimes e_{i_{n},i_{n}})\\
=&\sum_{r_{n}, r_{n-1},\dots, r_{1}, r_0}
p^{(0)}_{i_0}P_{r_0;i_0,r_{1}}P_{r_{1};i_{1},r_{2}} \cdots P_{r_{n-1};i_{n-1},r_{n}}P_{r_{n};i_{n}}\notag
\end{align}
Conversely, given a family of numbers $(P_{r;m,i})_{r,m,i\in \{1,\cdots ,d\}}$ satisfying conditions \eqref{(2.2a)} and \eqref{(2.2b)} and a probability measure
$(p^{(0)}_{i})_{i\in \{1,\cdots, d\}}$ on $\{1,\cdots, d\}$, there exists a unique classical stochastic process $(X_{n})_{n\in \mathbb{N}}$ with state space $\{1,\cdots, d\}$, whose joint
probabilities are given by the right hand side of \eqref{e-diag-proc} i.e., for all
$i_0, i_1,\dots ,i_n\in \{1,\cdots, d\}$,
\begin{align}\label{(2.2d)}
&\hbox{Prob}\bigl(X_0=i_0,X_1=i_1,\cdots ,X_n=i_n\bigr)\\
=&\sum_{r_{n}, r_{n-1},\dots, r_{1}, r_0}
p^{(0)}_{i_0}P_{r_0;i_0,r_{1}}P_{r_{1};i_{1},r_{2}} \cdots P_{r_{n-1};i_{n-1},r_{n}}P_{r_{n};i_{n}}\notag
\end{align}
with the $P_{r_{n};i_{n}}$ given by \eqref{(3.11b)2}.
}\end{theorem}
\textbf{Proof}.
Given \eqref{(3.11b)}, \eqref{(2.2a)} is clear. \eqref{(2.2b)} follows from the fact that, for all
$r\in \{1,\cdots, d\}$, one has
$$
\sum_{m,i}P_{r;m,i}= \sum_{m,i}\Big(\sum_{j\in D_{\mathcal{E}}}|K^*_{j,m,r}|^2\Big)(i,i)
 \ \overset{\eqref{calE(1)=1-diag}}{=} \ 1
$$
For each $m,i\in\{1,\cdots ,d\}$, one has,
\begin{equation}\label{(3.5)T}
\mathcal{E}(e_{mm}\otimes e_{ii})
\ \overset{\eqref{(3.4)gen2}, \eqref{(3.11b)}}{=} \
\sum_{r\in D}e_{rr} P_{r;m,i}
\end{equation}
Therefore
\begin{equation}\label{(3.13a)}
\mathcal{E}(e_{i_1,i_1}\otimes
\cdots \otimes \mathcal{E}(e_{i_{n-1}i_{n-1}}\otimes
\mathcal{E}(e_{i_{n},i_{n}}\otimes 1))\cdots )
\end{equation}
$$
\ \overset{\eqref{(3.5)T}, \eqref{(3.11b)2}}{=} \
\sum_{i\in D}\sum_{r_{n}} \mathcal{E}(e_{i_1,i_1}\otimes
\cdots \otimes \mathcal{E}(e_{i_{n-1}i_{n-1}}\otimes e_{r_{n}r_{n}})\dots )P_{r_{n};i_{n},i}
$$
$$
\ \overset{\eqref{(3.11b)}}{=} \
\sum_{r_{n}} \mathcal{E}(e_{i_1,i_1}\otimes
\cdots e_{i_{n-2}i_{n-2}}\otimes
\mathcal{E}(e_{i_{n-1}i_{n-1}}\otimes e_{r_{n}r_{n}} )\dots )P_{r_{n};i_{n}}
$$
$$
\ \overset{\eqref{(3.5)T}}{=} \
\sum_{r_{n}, r_{n-1}} \mathcal{E}(e_{i_1,i_1}\otimes
\cdots \mathcal{E}(e_{i_{n-2}i_{n-2}}\otimes e_{r_{n-1}r_{n-1}}) \dots )P_{r_{n-1};i_{n-1},r_{n}}P_{r_{n};i_{n}}
$$
Suppose by induction that
\begin{align}\label{(3.13b)}
&\mathcal{E}(e_{i_1,i_1}\otimes
\cdots \otimes \mathcal{E}(e_{i_{n-1}i_{n-1}}\otimes
\mathcal{E}(e_{i_{n},i_{n}}\otimes 1))\cdots )\\
= &\sum_{r_{n}, r_{n-1},\dots, , r_{1}} e_{r_{1}r_{1}}
P_{r_{1};i_{1},r_{2}} \cdots P_{r_{n-1};i_{n-1},r_{n}}P_{r_{n};i_{n}}\notag
\end{align}
Then
\begin{align}\label{(3.13c)}
&\mathcal{E}(e_{i_0,i_0}\otimes\mathcal{E}(e_{i_1,i_1}\otimes
\cdots \otimes \mathcal{E}(e_{i_{n-1}i_{n-1}}\otimes
\mathcal{E}(e_{i_{n},i_{n}}\otimes 1))\cdots ))\\
\overset{\eqref{(3.13b)}}{=}& \sum_{r_{n}, r_{n-1},\dots, r_{1}}
\mathcal{E}(e_{i_0,i_0}\otimes e_{r_{1}r_{1}})
P_{r_{1};i_{1},r_{2}} \cdots P_{r_{n-1};i_{n-1},r_{n}}P_{r_{n};i_{n}}\notag\\
\overset{\eqref{(3.5)T}}{=} &
\sum_{r_{n}, r_{n-1},\dots, r_{1}, r_0}
e_{r_0r_0} P_{r_0;i_0,r_{1}}P_{r_{1};i_{1},r_{2}} \cdots P_{r_{n-1};i_{n-1},r_{n}}P_{r_{n};i_{n}}\notag
\end{align}
It follows by induction that \eqref{cal-E0]-cal-E-De-pres} holds.
Therefore, in the notations of the statement,
\begin{align}\label{(3.13d)}
&\varphi (e_{i_0,i_0}\otimes e_{i_1,i_1}\otimes \cdots \otimes e_{i_{n},i_{n}})\\
=& \overline{\hbox{Tr}}_{0}\left(p^{(0)}\mathcal{E} (e_{i_0,i_0}\otimes\mathcal{E}(e_{i_1,i_1}\otimes\cdots \otimes \mathcal{E}(e_{i_{n-1}}e_{i_{n-1}}^*\otimes\mathcal{E} (e_{i_{n},i_{n}}\otimes 1))\cdots ))\right)\notag\\
=&\sum_{r_{n}, r_{n-1},\dots, r_{1}, r_0}
p^{(0)}_{i_0}P_{r_0;i_0,r_{1}}P_{r_{1};i_{1},r_{2}} \cdots P_{r_{n-1};i_{n-1},r_{n}}P_{r_{n};i_{n}}\notag
\end{align}
which is \eqref{e-diag-proc}.
To prove the converse statement, we have to show that, defining the left hand side of \eqref{(2.2d)}
by its right hand side, Kolmogorov compatibility conditions are satisfied. This follows from
\begin{align*}
&\sum_{i_{n}\in D}\hbox{Prob}\bigl(X_0=i_0,X_1=i_1,\cdots ,X_n=i_n\bigr)\\
=&\sum_{r_{n}, r_{n-1},\dots, r_{1}, r_0}p^{(0)}_{i_0} P_{r_0;i_0,r_{1}}P_{r_{1};i_{1},r_{2}} \cdots P_{r_{n-1};i_{n-1}, r_{n}}\sum_{i_{n}\in D}P_{r_{n};i_{n}}\\
\overset{\eqref{(3.11b)2}}{=}&\sum_{r_{n}, r_{n-1},\dots, r_{1}, r_0}p^{(0)}_{i_0}P_{r_0;i_0,r_{1}}P_{r_{1};i_{1},r_{2}} \cdots P_{r_{n-1};i_{n-1},r_{n}}
\sum_{i_{n}, i\in D}P_{r_{n};i_{n},i}
\end{align*}
\begin{align*}
\overset{\eqref{(2.2b)}}{=}&\sum_{r_{n}, r_{n-1},\dots, r_{1}, r_0} p^{(0)}_{i_0}P_{r_0;i_0,r_{1}}P_{r_{1};i_{1},r_{2}} \cdots P_{r_{n-1};i_{n-1},r_{n}}\\
\overset{\eqref{(3.11b)2}}{=}&\sum_{r_{n-1},\dots, r_{1}, r_0}
p^{(0)}_{i_0}P_{r_0;i_0,r_{1}}P_{r_{1};i_{1},r_{2}} \cdots P_{r_{n-1};i_{n-1}}\\
=&\sum_{i_{n}\in D}\hbox{Prob}\bigl(X_0=i_0,X_1=i_1,\cdots ,X_n=i_{n-1}\bigr) \tag*{$\square$}
\end{align*}

\begin{theorem}\label{2.6}{\rm
The classical stochastic process $(X_{n})_{n\in \mathbb{N}}$  with state space $\{1,\cdots, d\}$,
obtained by restriction to the $e$--diagonal algebra of a $D_{e}$--preserving quantum Markov
chain, is a classical Markov chain if and only if the operators $\{K_r\}$ (see \eqref{df-cal-E-Tr-CDA})
have the form:
\begin{equation}\label{(2.6)}
K_r=\sum_{h\in D}e_{hh}\otimes K_{r,h,h}\ ,\qquad\qquad
\forall\ j\in D_{\mathcal{E}}
\end{equation}
}\end{theorem}
\textbf{Proof}.
The classical process $(X_{n})_{n\in \mathbb{N}}$ in the statement is a classical Markov
chain if and only if
\begin{equation}\label{(2.5)}
\mathcal{E}(a\otimes b)=a\mathcal{E}(1\otimes b) = aP(b)
\quad,\quad \forall a,b\in D_{e}
\end{equation}
where $P$ is a Markov operator on $D_{e}$.
This is equivalent to say
\begin{align*}
&\sum_{r\in D}e_{r,r} Tr( \sum_{j\in D_{\mathcal{E}}}|K_{j,m,r}|^2 e_{ii})=\sum_{r\in D}e_{r,r} Tr( \sum_{j\in D_{\mathcal{E}}}K^*_{j,m,r} e_{ii} K_{j,m,r})\\
=&\mathcal{E}(e_{mm}\otimes e_{ii})=e_{mm}\mathcal{E}(1\otimes e_{ii})=e_{mm} \sum_{h\in D}\mathcal{E}(e_{hh}\otimes e_{ii})\\
=&e_{mm} \sum_{h\in D}\sum_{r\in D}e_{r,r} Tr( \sum_{j\in D_{\mathcal{E}}}|K_{j,h,r}|^2e_{ii})
=e_{mm} \sum_{h\in D} Tr( \sum_{j\in D_{\mathcal{E}}}|K_{j,h,m}|^2e_{ii})
\end{align*}
and so
\begin{equation}\label{nec-cond-e-proc-Mark}
\sum_{r\in D}e_{r,r} Tr( \sum_{j\in D_{\mathcal{E}}}|K_{j,m,r}|^2 e_{ii}) =e_{mm} \sum_{h\in D} Tr( \sum_{j\in D_{\mathcal{E}}}|K_{j,h,m}|^2e_{ii})
\end{equation}
\eqref{nec-cond-e-proc-Mark} implies that
$$
Tr( \sum_{j\in D_{\mathcal{E}}}|K_{j,m,r}|^2 e_{ii}) = 0  \ ,\quad\forall i\in D\,,\ m\ne r
$$
equivalently
$$
Tr( \sum_{j\in D_{\mathcal{E}}}|K_{j,m,r}|^2) = 0\ ,\quad \forall m\ne r
$$
This equality and the positivity of the operator $\vert K_{j,m,r}\vert ^2$ guarantee that $\vert K_{j,m,r}\vert ^2=0$ whenever $m\ne r$.
Conversely, if this condition is satisfied, one has trivially
\eqref{(2.6)}.
$\qquad\square$
\begin{corollary}\label{charact-diag-QMC-|DE|=1}{\rm
In the assumptions of Theorem \ref{2.6}, suppose in addition that
in the notation \eqref{(2.6)}, $|D_{\mathcal{E}}|=1$ (i.e. there is only one conditional density amplitude $K$). Then $\varphi$ is diagonalizable in the sense of Definition \ref{df:diagonalizable-QMC}, i.e.
$K$ has the form \eqref{K-e-diag}.
}\end{corollary}
\textbf{Proof}.
By Theorem \ref{2.6} and the diagonalizability econdition, $K$ has the form
$$
K=\sum_{hh\in D}e_{hh}\otimes K_{h}
$$
and satisfies
$$
0=[K\otimes 1, 1\otimes K]
= [\sum_{h}e_{hh}\otimes K_{h}\otimes 1, 1\otimes \sum_{h'}e_{h'h'}\otimes K_{h'}]
$$
$$
\iff \sum_{h,h'}e_{hh}\otimes K_{h}e_{h'h'}\otimes K_{h'}
= \sum_{h,h'}e_{hh}\otimes e_{h'h'}K_{h}\otimes K_{h'}
$$
Multiplying both sides by a fixed $e_{hh}\otimes 1\otimes 1$, one finds
$$
e_{hh}\otimes \sum_{h'}K_{h}e_{h'h'}\otimes K_{h'}
= e_{hh}\otimes \sum_{h'}e_{h'h'}K_{h}\otimes K_{h'}
\quad,\quad\forall h\in D
$$
$$
\iff \sum_{h'}K_{h}e_{h'h'}\otimes K_{h'}
= \sum_{h'}e_{h'h'}K_{h}\otimes K_{h'}
\quad,\quad\forall h\in D
$$
Multiplying on the left both sides by a fixed $e_{kk}\otimes 1$, one finds
$$
\sum_{h'}e_{kk}K_{h}e_{h'h'}\otimes K_{h'}
= e_{kk}K_{h}\otimes K_{k,k}
\quad,\quad\forall h, k\in D
$$
Multiplying on the right both sides by a fixed $e_{k'k'}\otimes 1$, one finds
$$
e_{kk}K_{h}e_{k'k'}\otimes K_{k',k'}
= e_{kk}K_{h}e_{k'k'}\otimes K_{k,k}
\quad,\quad\forall h, k, k'\in D
$$
$$
\iff e_{kk}K_{h}e_{k'k'}\otimes (K_{k',k'}-K_{k,k}) = 0
\quad,\quad\forall h, k, k'\in D
$$
Since the $K_{m,n}$ can always be supposed to be independent up to a change of the index set
$D_{\mathcal{E}}$, this is equivalent to say that
$$
e_{kk}K_{h}e_{k'k'} = 0
\quad,\quad\forall h, k, k'\in D \ , \ k\ne k'
$$
i.e.,
$$
K_{h} = \sum_{j\in D} K_{hj}e_{jj}
\quad,\quad\forall h \in D \ , \ \hbox{for some } K_{hj}\in\mathbb{C}
$$
We know that the normalization condition implies that $(|K_{hj}|^2)_{h,j\in D}$ is a stochastic
matrix. Therefore $K$ has the form \eqref{K-e-diag} and $\varphi$ is diagonalizable.
$\qquad\square$

\subsection{Diagonal restrictions of $D_{e}$--preserving quantum Markov chains}
\label{sec:Diag-restr-De-pres-QMC}

Theorem \ref{2.3} shows that, if a transition expectation $\mathcal{E}$ is $D_{e}$--preserving
for some ortho--normal basis $e\equiv (e_{j})_{j\in D}$ of $H\equiv\mathbb{C}^{d}$, then for
any $e$--diagonal state $\varphi_0$ on $M$ (see \eqref{e-diag-w(0)}), the restriction of the
quantum Markov chain $\varphi$ on $\bigotimes_{n\in\mathbb{N}} M$, defined by the pair
$(\varphi_0, \mathcal{E})$, to the diagonal algebra\\
$\mathcal{D}_{e} = \bigotimes_{n\in\mathbb{N}} D_{e}$ gives rise to a classical stochastic process
$(X_{n})_{n\in \mathbb{N}}$, with state space $\{1,\cdots, d\}$, whose joint probabilities are
uniquely determined by $\varphi_0$ and a family of real numbers
$(P_{r;m,i})_{r,m,i\in \{1,\cdots ,d\}}$ satisfying conditions \eqref{(2.2a)} and \eqref{(2.2b)}
through the identity \eqref{e-diag-proc}. Conversely, any such a triple determines a unique
classical process. The following theorem (see Theorem 2.4 of \cite{LuYG95-QMC}) characterizes the classical processes obtained in this way as restrictions of a special class of Markov processes.
\begin{theorem}\label{hmm-th12}{\rm
Let $P\equiv (P_{r;m,i})$ be a $(d\times d^2)$--stochastic matrix as in Theorem \ref{2.3} and let
$X \equiv (X_n)_{n\in\mathbb{N}}$ be a classical stochastic process with state space
$\{1,\cdots, d\}$ and joint probabilities given by the right hand side of \eqref{e-diag-proc}.
Define the probability measure $P^{(0)}\equiv p^{(0)}(\cdot,\cdot)$ on $\{1,\cdots ,d\}^2$ by
\begin{equation}\label{(3.19)}
p^{(0)}(j,i):=\sum_r p_0(r)P_{r;i,j}
\end{equation}
and the matrix $P_{D^2}\equiv\bigl(p_{(j,i),(j',i')}\bigr)$ by
\begin{equation}\label{df-PD2}
p_{(j,i),(j',i')}:=P_{j;i',j'}\ ,\quad  \forall\ (j,i), (j',i')\in \{1,\cdots ,d\}^2
\end{equation}
and let $(Z,Y)\equiv\{(Z_n,Y_n)\}_{n=0}^\infty$ be the classical Markov chain with state space
$\{1,\cdots, d\}^2$ determined by the pair $(P^{(0)},P_{D^2})$.\\

\noindent Then, denoting $P_{X}$ (resp. $P_{Z,Y}$) the probability distribution of $X$
(resp. $(Z,Y)$), one has
\begin{align}\label{(2.4)}
&P_{X}\bigl(X_0=i_0,X_1=i_1,\cdots ,X_n=i_n\bigr)\\
=&P_{Z,Y}\bigl(Y_0=i_0,Y_1=i_1,\cdots ,Y_n=i_n\bigr)\notag
\end{align}
in other words, the process $X$ is stochastically equivalent to the sub--process $Y$, obtained as
restriction of the process $(Z,Y)$ to the $\sigma$--algebra generated by the $Y_n$.\\
Conversely, if $(Z,Y)\equiv\{(Z_n,Y_n)\}_{n=0}^\infty$ is a classical Markov chain with state
space $\{1,\cdots, d\}^2$ determined by the pair $(P^{(0)},P)$ such that:\\
(i) $P\equiv\bigl(p_{(j,i),(j',i')}\bigr)$ with
\begin{equation}\label{p((j,i),(j',i'))-indep-i}
p_{(j,i),(j',i')}  \ \hbox{independent of } \ i
\end{equation}
(ii) $P^{(0)}\equiv p^{(0)}(\cdot,\cdot)$ is the probability measure on $\{1,\cdots ,d\}^2$
defined by
\begin{equation}\label{df-p(0)(j,i)}
p^{(0)}(j,i):=\sum_{r\in D} p_0(r)p_{(r,i),(j,i)}
\end{equation}
where $(p_0(r))_{r\in D}$ is an arbitrary probability measure on$\{1,\cdots ,d\}$,\\
then there exists a quantum Markov chain $\varphi$ on $\bigotimes_{n\in\mathbb{N}}M$ and a diagonal
algebra $\mathcal{D}_{e}:= \bigotimes_{n\in\mathbb{N}}D_{e}\subseteq \bigotimes_{n\in\mathbb{N}}M$
such that the restriction of $\varphi$ on $\mathcal{D}_{e}$ is a classical random process
stochastically equivalent to the sub--process $(Y_n)_{n\in\mathbb{N}}$ of $(Z,Y)$.
Moreover the choice of $D_{e}$ is arbitrary.
}\end{theorem}
\textbf{Proof}.
Let the classical stochastic process $X$ be as in the first statement of the theorem.
Then the matrix $\bigl(p_{(j,i),(j',i')}\bigr)$ defined by \eqref{df-PD2} is a transition matrix
on $\{1,\cdots ,d\}^2$ because of \eqref{(2.2b)}.
Then, if $P^{(0)}$ is defined by \eqref{(3.19)}, the pair $(P^{(0)},P_{D^2})$ determines
a, unique up to stochastic equivalence, classical Markov chain
$(Z,Y)\equiv\{(Z_n,Y_n)\}_{n=0}^\infty$ with state space $\{1,\cdots, d\}^2$.
By construction, the joint distribution $P_{Y}$ of the sub--process $Y$ is given by
$$
P_{Y}(Y_0=i_0,\cdots ,Y_1=i_1)
=\sum_{r_{n}, r_{n-1},\dots, r_{1}, r_0}
p^{(0)}_{i_0}P_{r_0;i_0,r_{1}}P_{r_{1};i_{1},r_{2}} \cdots P_{r_{n-1};i_{n-1},r_{n}}P_{r_{n};i_{n}}
$$
where the $P_{r;m}$ are given by \eqref{(3.11b)2}. Comparing this with the right hand side
of \eqref{(2.2d)} one sees that \eqref{(2.4)} holds.\\
Conversely, given a classical Markov chain $\{(Z_n,Y_n)\}_{n=0}^\infty$ with state space
$\{1,\cdots ,d\}^2$ and determined by a pair $(P^{(0)},P)$ satisfying conditions (i) and (ii)
above, then, defining
$$
P_{j;i',j'} := p_{(j,i),(j',i')}\quad,\quad  \forall\ (j,i), (j',i')\in \{1,\cdots ,d\}^2
$$
the $(P_{j;i,k})$ satisfy \eqref{(2.2a)}, \eqref{(2.2b)}.
Therefore, by Theorem \ref{2.3}, for any choice of a probability measure $p_0$ on $\{1,\cdots ,d\}$,
and of a diagonal algebra $D_{e}\subset M$, the pair $(p_0, (P_{j;i,k})_{j,i,k\in D}$ defines
a quantum Markov state $\varphi$ whose joint probabilities on
$\mathcal{D}_{e}:=\bigotimes_{\mathbb{N}} D_{e}$ are given by \eqref{(2.2d)}.
Since, by construction, these coincide with the joint probabilities of the sub--process
$(Y_n)$, the thesis follows.
$\qquad\square$

\subsection{Hidden processes defined by diagonal restrictions of $D_{e}$--preserving quantum Markov chains}\label{sec:HP-Diag-restr-De-pres-QMC}

In this section we prove that, from diagonal restrictions of $D_{e}$--preserving
(non--diagonal) quantum Markov chains, one can also obtain time--consecutive Hidden Markov process
in the sense of Definition \eqref{df:time-cons-HMP}.\\

\noindent Let $(Z,Y)$ be a Markov process as in the second statement of Theorem \ref{hmm-th12}.
If $I$ is any sub--set of the set $\{Z_{k},Y_{k},\dots,Z_{0},Y_{0}\},$ we denote
$$
\mathcal{F}_{I}
:= \sigma\hbox{--algebra generated by the random variables in } I
$$
thus, for example, for $k\in\mathbb{N}$,
$$
\mathcal{F}_{(Z_{h},Y_{h})_{h=0}^{k}}
:= \sigma\hbox{--algebra generated by } \ \{(Z_{m},Y_{m})\}_{m\in \{0,\dots, k\} }
$$
$$
\mathcal{F}_{Z_{k}}
:= \sigma\hbox{--algebra generated by } \  Z_{k}
$$
We denote $E_{\mathcal{F}_{I}}$ the conditional expectation of the  $(Z,Y)$ process given the
$\sigma$--algebra $\mathcal{F}_{I}$. It is known that, if $I\subseteq J$, then
\begin{equation}\label{proj-cond-exp}
E_{\mathcal{F}_{I}}E_{\mathcal{F}_{J}} = E_{\mathcal{F}_{J}}E_{\mathcal{F}_{I}} = E_{\mathcal{F}_{I}}
\end{equation}
The Markovianity of the process $(Z,Y)$ means that
\begin{equation}\label{Mark-prop-L}
E_{\mathcal{F}_{(Z_{h},Y_{h})_{h=0}^{n}}}(G_{\mathcal{F}_{(Z_{h},Y_{h})_{h=n+1}^{\infty}}})
= E_{\mathcal{F}_{(Z_{n},Y_{n})}}(G_{\mathcal{F}_{(Z_{h},Y_{h})_{h=n+1}^{\infty}}})
\end{equation}
for any $\mathcal{F}_{(Z_{h},Y_{h})_{h=n+1}^{\infty}}$--measurable function
$G_{\mathcal{F}_{(Z_{h},Y_{h})_{h=n+1}^{\infty}}}$.
By assumption, the process $(Z,Y)$ satisfies condition \eqref{p((j,i),(j',i'))-indep-i} which is equivalent to
$$
P_{Y,Z}\bigl((Z_n,Y_n)\bigr| \mathcal{F}_{(Z_{h},Y_{h})_{h=0}^{n-1}})\bigr)
= P_{Y,Z}\bigl((Z_n,Y_n)\bigr| \ Z_{n-1}\bigr)
$$
i.e.,
\begin{equation}\label{HM2}
E_{\mathcal{F}_{(Z_{h},Y_{h})_{h=0}^{n-1}}}\bigl(f_{n}(Z_n,Y_n) \bigr)=  E_{\mathcal{F}_{Z_{n-1}}}\bigl(f_{n}(Z_n,Y_n)\bigr)
\end{equation}
for any $\mathcal{F}_{(Z_{n},Y_{n})}$--measurable function $f_{n}$. From \eqref{HM2}, one deduces that, since $\mathcal{F}_{Z_{n-1},(Z_{h},Y_{h})_{h=0}^{n-2}}
\subseteq\mathcal{F}_{(Z_{h},Y_{h})_{h=0}^{n-1}}$, one has
\begin{align}\label{HM3}
&E_{\mathcal{F}_{Z_{n-1},(Z_{h},Y_{h})_{h=0}^{n-2}}} \bigl(f_{n}(Z_n,Y_n)\bigr)\\
=&E_{\mathcal{F}_{Z_{n-1},(Z_{h},Y_{h})_{h=0}^{n-2}}}
E_{\mathcal{F}_{(Z_{h},Y_{h})_{h=0}^{n-1}}}\bigl(f_{n}(Z_n,Y_n) \bigr)\notag\\
\overset{\eqref{HM2}}{=}& E_{\mathcal{F}_{Z_{n-1},(Z_{h},Y_{h})_{h=0}^{n-2}}}
E_{\mathcal{F}_{Z_{n-1}}}\bigl(f_{n}(Z_n,Y_n)\bigr)
=E_{\mathcal{F}_{Z_{n-1}}}\bigl(f_{n}(Z_n,Y_n)\bigr)\notag
\end{align}
and this implies that
$$
E\bigl(f_{n}(Z_n,Y_n)\cdots  f_{0}(Z_0,Y_0)\bigr)
$$
$$
=E\bigl(E_{\mathcal{F}_{(Z_{h},Y_{h})_{h=0}^{n-1}}}\bigl(f_{n}(Z_n,Y_n)\bigr)
f_{n-1}(Z_{n-1},Y_{n-1})\cdots  f_{0}(Z_0,Y_0)\bigr)
$$
$$
\overset{\eqref{HM2}}{=}
E\bigl(E_{\mathcal{F}_{Z_{n-1}}}\bigl(f_{n}(Z_n,Y_n)\bigr)
f_{n-1}(Z_{n-1},Y_{n-1})\cdots  f_{0}(Z_0,Y_0)\bigr)
$$
$$
=E\bigl(E_{\mathcal{F}_{Z_{n-1}}}\bigl(f_{n}(Z_n,Y_n)\bigr)\cdot
$$
$$
\cdot \ E_{\mathcal{F}_{Z_{n-1},(Z_{h},Y_{h})_{h=0}^{n-2}}}
\bigl(f_{n-1}(Z_{n-1},Y_{n-1})\bigr)f_{n-2}(Z_{n-2},Y_{n-2})\cdots  f_{0}(Z_0,Y_0)\bigr)\bigr)
$$
$$
\overset{\eqref{HM2}}{=}
E\bigl(E_{\mathcal{F}_{Z_{n-1}}}\bigl(f_{n}(Z_n,Y_n)\bigr)
E_{\mathcal{F}_{Z_{n-1}}}\bigl(
\bigl(f_{n-1}(Z_{n-1},Y_{n-1})\bigr)f_{n-2}(Z_{n-2},Y_{n-2})\cdots  f_{0}(Z_0,Y_0)\bigr)\bigr)
$$
Suppose by induction that
\begin{align}\label{HM-ind-ass}
&E\bigl(f_{n}(Z_n,Y_n)\cdots  f_{0}(Z_0,Y_0)\bigr)\\
=&E\bigl(\prod_{h=n-(k-1)}^{n}E_{\mathcal{F}_{Z_{n-1}}} \bigl(f_{h}(Z_{h},Y_{h})\bigr)
\prod_{m=0}^{n-k}  f_{m}(Z_{m},Y_{m})\bigr)\notag
\end{align}
Then
$$
E\bigl(f_{n}(Z_n,Y_n)\cdots  f_{0}(Z_0,Y_0)\bigr)
$$
$$
=E\bigl(\prod_{h=n-(k-1)}^{n}E_{\mathcal{F}_{Z_{n-1}}}\bigl(f_{h}(Z_{h},Y_{h})\bigr)
f_{n-k}(Z_{n-k},Y_{n-k})\prod_{m=0}^{n-k-1}  f_{m}(Z_{m},Y_{m})\bigr)
$$
$$
=E\bigl(\prod_{h=n-(k-1)}^{n}E_{\mathcal{F}_{Z_{n-1}}}\bigl(f_{h}(Z_{h},Y_{h})\bigr)
E_{\mathcal{F}_{Z_{n-1},(Z_{h},Y_{h})_{h=0}^{n-k-1}}}
\bigl(f_{n-k}(Z_{n-k},Y_{n-k})\bigr)\cdot
$$
\begin{equation}\label{HM4}
\cdot \ \prod_{m=0}^{n-k-1}  f_{m}(Z_{m},Y_{m})\bigr)
\end{equation}
Since $\mathcal{F}_{Z_{n-1},(Z_{h},Y_{h})_{h=0}^{n-k-1}}
\subseteq \mathcal{F}_{Z_{n-1},(Z_{h},Y_{h})_{h=0}^{n-2}}$,
it follows that
\begin{align}\label{HM5}
&E_{\mathcal{F}_{Z_{n-1},(Z_{h},Y_{h})_{h=0}^{n-k-1}}}
\bigl(f_{n-k}(Z_{n-k},Y_{n-k})\bigr)\\
=& E_{\mathcal{F}_{Z_{n-1},(Z_{h},Y_{h})_{h=0}^{n-k-1}}}
E_{\mathcal{F}_{Z_{n-1},(Z_{h},Y_{h})_{h=0}^{n-2}}}
\bigl(f_{n-k}(Z_{n-k},Y_{n-k})\bigr)\notag\\
\overset{\eqref{HM3}}{=}& E_{\mathcal{F}_{Z_{n-1},(Z_{h},Y_{h})_{h=0}^{n-k-1}}}
E_{\mathcal{F}_{Z_{n-1}}}\bigl(f_{n-k}(Z_{n-k},Y_{n-k}\bigr)
\notag\\
=&E_{\mathcal{F}_{Z_{n-1}}}\bigl(f_{n-k}(Z_{n-k},Y_{n-k}\bigr) \notag
\end{align}
In conclusion
$$
E\bigl(f_{n}(Z_n,Y_n)\cdots  f_{0}(Z_0,Y_0)\bigr)
$$
$$
\overset{\eqref{HM4}}{=}E\bigl(\prod_{h=n-(k-1)}^{n}E_{\mathcal{F}_{Z_{n-1}}}\bigl(f_{h}(Z_{h},Y_{h})\bigr)
E_{\mathcal{F}_{Z_{n-1},(Z_{h},Y_{h})_{h=0}^{n-k-1}}}
\bigl(f_{n-k}(Z_{n-k},Y_{n-k})\bigr)
$$
$$
\prod_{m=0}^{n-k-1}  f_{m}(Z_{m},Y_{m})\bigr)
$$
$$
\overset{\eqref{HM5}}{=}
E\bigl(\prod_{h=n-(k-1)}^{n}E_{\mathcal{F}_{Z_{n-1}}}\bigl(f_{h}(Z_{h},Y_{h})\bigr)
E_{\mathcal{F}_{Z_{n-1}}}\bigl(f_{n-k}(Z_{n-k},Y_{n-k}\bigr)
\prod_{m=0}^{n-k-1}  f_{m}(Z_{m},Y_{m})\bigr)
$$
$$
=E\bigl(\prod_{h=n-k}^{n}E_{\mathcal{F}_{Z_{n-1}}}\bigl(f_{h}(Z_{h},Y_{h})\bigr)
\prod_{m=0}^{n-k-1}  f_{m}(Z_{m},Y_{m})\bigr)
$$
Thus \eqref{HM-ind-ass} holds for $k+1$, hence, by induction, we conclude that \eqref{HM-ind-ass} holds for all $k\le n$. Putting $k = n-1$ in \eqref{HM-ind-ass}, one finds
\begin{align}\label{HM3b}
&E\bigl(f_{n}(Z_n,Y_n)\cdots  f_{0}(Z_0,Y_0)\bigr)\\
=&E\bigl(\prod_{h=1}^{n}E_{\mathcal{F}_{Z_{n-1}}}\bigl(f_{h}(Z_{h},Y_{h})\bigr)f_{0}(Z_{0},Y_{0})\bigr)\notag\\
=&E\bigl(\prod_{h=1}^{n}E_{\mathcal{F}_{Z_{n-1}}} \bigl(f_{h}(Z_{h},Y_{h})\bigr)E_{\mathcal{F}_{Z_{n-1}}} \bigl(f_{0}(Z_{0},Y_{0})\bigr)\bigr)\notag\\
=&E\Big(\prod_{h=0}^{n}E_{\mathcal{F}_{Z_{n-1}}} \bigl(f_{h}(Z_{h},Y_{h})\bigr)\Big)\notag
\end{align}
Choosing
$$
f_{h}(Z_{h},Y_{h}) := F_{h}(Z_{h})G(Y_{h})
$$
one has
\begin{equation}\label{(Z,Y)-HP1}
E_{Z,Y}\bigl(\prod_{h=0}^{n}F_{h}(Z_{h})G(Y_{h})\bigr)
=E_{Z}\bigl(\prod_{h=0}^{n}E_{\mathcal{F}_{Z_{n-1}}}
\bigl(F_{h}(Z_{h})G(Y_{h})\bigr)\bigr)
\end{equation}
Now notice that, since
$\mathcal{F}_{(Z_{h})_{h=0}^{n-1}}\subseteq\mathcal
{F}_{(Z_{h},Y_{h})_{h=0}^{n-1}}$,
one has
\begin{align}\label{Z-Mark1}
&E_{\mathcal{F}_{(Z_{h})_{h=0}^{n-1}}}\bigl(f_{n}(Z_n,Y_n)\bigr)
=E_{\mathcal{F}_{(Z_{h})_{h=0}^{n-1}}}E_{\mathcal{F}_ {(Z_{h},Y_{h}) _{h=0}^{n-1}}}\bigl(f_{n}(Z_n,Y_n)\bigr)\notag\\
=&E_{\mathcal{F}_{(Z_{h})_{h=0}^{n-1}}}E_{\mathcal{F}_{Z_{n-1}}} \bigl(f_{n}(Z_n,Y_n)\bigr)= E_{\mathcal{F}_{Z_{n-1}}} \bigl(f_{n}(Z_n,Y_n)\bigr)\notag
\end{align}
In particular, if $f_{n}$ does not depend on the second variable, one has
\begin{equation}\label{Z-Mark2}
E_{\mathcal{F}_{(Z_{h})_{h=0}^{n-1}}}\bigl(f_{n}(Z_n)\bigr)
= E_{\mathcal{F}_{Z_{n-1}}}\bigl(f_{n}(Z_n)\bigr)
\end{equation}
i.e. the sub--process $(Z_n)$ is Markov.
Thus, denoting $(P_{Z;n,n+1})$ the sequence of its transition operators, one has
\begin{align}\label{(Z,Y)-HP2}
&E_{Z,Y}\bigl(\prod_{h=0}^{n}f_{h}(Z_{h},Y_{h})\bigr)
\overset{\eqref{(Z,Y)-HP1}}{=}E_{Z}\bigl(\prod_{h=0}^{n} E_{\mathcal{F}_{Z_{n-1}}}\bigl(f_{h}(Z_{h},Y_{h})
\bigr)\bigr)\\
=&E_{Z_{0}}\Big(P_{Z;0,1}P_{Z;1,2} \cdots P_{Z;n-2,n-1}
\Big(\prod_{h=0}^{n}E_{\mathcal{F}_{Z_{n-1}}} \left(f_{h}(Z_{h},Y_{h}\right)\Big)\Big)\notag
\end{align}
and one recognizes that, up to notations, the right hand side of \eqref{(Z,Y)-HP2}
is \textbf{exactly of the form} that characterizes time consecutive hidden Markov process
in the sense of Definition \eqref{df:time-cons-HMP}.\\

\noindent Notice that, if $(Z,Y)$ where a HMP with hidden process $Z$, one should have
\begin{equation}\label{joint-probs-Y-Z-Mark}
P_{Z,Y}\Big(\prod_{m=0}^{n}f_m(Y_m)g_m(Z_m)\Big)
\end{equation}
\[=p_{Z_{0}}\left( B_{Y,Z,0}(f_{0})g_{0}P_{Z_{0}} \left(B_{Y,Z,1}(f_{1})g_{1}\right)\left(\cdots P_{Z_{n-1}} \left(B_{Y,Z,n}(f_{n})g_{n}\right)\right)\right)
\]
Comparing \eqref{HM3b} and \eqref{joint-probs-Y-Z-Mark} one sees that $(Z,Y)$ is a
time--consecutive hidden Markov process (with hidden process $Z$) in the sense of Definition \ref{df:class-gen-hidd-proc}, \textbf{but not a hidden Markov process} in the sense
of Definition \ref{df:class-HP}.


\begin{thebibliography}{99}



\bibitem[Ac73-past-dep-PE]{Ac73-past-dep-PE}
L. Accardi: {\it On a class of measures connected with past dependent probability evolutions},
Avtomatika i Telemechanika,(in Russian),  p.50--61 (1974),
Inner report, Laboratorio di cibernetica (1973)

\bibitem[Ac74d-Camerino]{Ac74-Camerino}
L. Accardi: {\it Non--commutative Markov chains},
Proceedings International School of Mathematical Physics,
Universit\`a di Camerino 30 Sept., 12 Oct. p.268-295 (1974)

\bibitem[Ac74-FAA]{Ac74-FAA}
L. Accardi: {\it The noncommutative markovian property},
Func. Anal. Appl., (in russian) 9 (1), p.1-8 (1975),
submitted 31-1-1974, translated in: Funct. Anal. and its Appl. 9 (1) p.1-8  (1975)

\bibitem[Ac81-Topics-QP]{Ac81-Topics-QP}
L. Accardi: {\it Topics in quantum probability},
Phys. Rep. 77 (3), p.169-192 (1981). Review section of Physics Letters (eds. C. DeWitt-Morette, K.D. Elworthy)

\bibitem[Ac82c-Some trends and problems]{Ac82-Some trends and problems}
L. Accardi: {\it Some trends and problems in quantum probability},
 Quantum  probability and applications to the
quantum theory of irreversible processes,
L. Accardi, A. Frigerio and V. Gorini (eds.),
Proc. 2-d Conference: Quantum Probability and
applications to the quantum theory of irreversible
processes, 6-11, 9 (1982) Villa  Mondragone (Rome),
Springer LNM N. 1055, p.1-19 (1984)\\

\bibitem[Ac91-Q-Kalman-filters]{Ac91-Q-Kalman-filters}
Accardi L.:\\
Quantum Kalman filters,
in: Mathematical system theory,\\
The influence of R.E. Kalman, A.C. Antoulas (ed.)
Springer (1991) 135-143\\
Invited contribution to the memorial volume for
the 60-th birthday of R.E.Kalman\\

\bibitem[AcFi01b]{AcFi01b}
 L. Accardi, F. Fidaleo: {\it On the structure of quantum Markov fields}, Proceedings Burg Conference 15--20 March 2001, W. Freudenberg (ed.), World Scientific, QP--PQ Series 15 p.1--20 (2003)

\bibitem[AcFi03-EMC]{AcFi03-EMC}
L. Accardi, F. Fidaleo: {\it Entangled Markov Chains}
Annali di Matematica Pura e Applicata, 184 (3), p.327--346 (2005).
DOI: 10.1007/s10231-004-0118-4. Preprint Volterra N.556 (2003)
http://www.springerlink.com/index/10.1007/s10231-004-0118-4

\bibitem[AcFiMu07]{AcFiMu07}
L. Accardi, F. Fidaleo, F. Mukhamedov:
 {\it Markov states and chains on the CAR algebra},
Infin. Dimen. Anal. Quantum Probab. Relat. Top. (IDAQP) 10, p.165--183 (2007)

\bibitem[AcSouElG20]{AcSouElG20}
L. Accardi, A. Souissi, El. Soueidy: {\it Quantum Markov chains, A unification approach}, Infin. Dimens. Anal. Quantum Probab. Relat. Top. (IDAQP)  23 (2) (2020)

\bibitem[AcWa87]{AcWa87}
L. Accardi, G.S. Watson:
{\it Markov states of the quantum electromagnetic field}, Phys. Rev. A 35, p.1275--1283 (1987)

\bibitem[Algh2016]{Algh2016}
R. Alghamdi:
{\it Hidden Markov models (HMMs) and security applications}.
Int. J. Adv. Comput. Sci. Appl. 7 (2), p.39--47 (2016)

\bibitem[BaumPetr66]{BaumPetr66}
L.E. Baum, T. Petrie:
{\it Statistical inference for probabilistic functions of finite state Markov chains}.
The Annals of Mathematical Statistics, 37, p.1554-1563 (1966)

\bibitem[CGGK2017]{CGGK2017}
M. Cholewa, P. Gawron, P. Glomb, D. Kurzyk:
{\it Quantum hidden Markov models based on transition operation matrices}.
Quantum Information Processing, 16 (4), p.1--19 (2017)

\bibitem[Eddy98]{Eddy98}
S.R. Eddy: {\it Profile hidden Markov models}. Bioinformatics, 14 (9),  p.755-763 (1998)

\bibitem[JasKell12]{JasKell12} J. Ernst, M. Kellis:
{\it ChromHMM: automating chromatin-state discovery and characterization}. Nature methods, 28, 9, p.215-216 (2012)

\bibitem[FaNaWe92]{FaNaWe92}
M. Fannes, B. Nachtergaele B., R.F. Werner: {\it Finitely correlated states on quantum spin chains},
Commun. Math. Phys., v.144, p.443--490 (1992)

\bibitem[FelsChur92]{FelsChur92}
J. Felsenstein, G.A. Churchill:
{\it A Hidden Markov Model approach to variation among sites in rate of evolution}.
 Molecular biology and evolution, 13 (1), p.93-104 (1996)

\bibitem[GhahrJord97]{GhahrJord97}
Z. Ghahramani, M.I. Jordan:
{\it Factorial Hidden Markov Models}. Machine Learning, 29,  p. 245-273 (1997)

\bibitem[HuaYasMerv90]{HuaYasMerv90}
X.D. Huang, A. Yasuo, M. Jack:
{\bf Hidden Markov Models for Speech Recognition},  Columbia University Press (1990) ISBN:978-0-7486-0162-2

\bibitem[JelBahMer75]{JelBahMer75}
F. Jelinek, L. Bahl, R. Mercer:
 {\it Design of a linguistic statistical decoder for the recognition of continuous speech}.
IEEE Transactions on Information Theory, 21 (3), p.250-256 (1975)

\bibitem[HasNat05]{HasNat05}
M.R. Hassan, B. Nath,
{\it Stock market forecasting using hidden Markov model: a new approach},
5th International Conference on Intelligent Systems Design and Applications (ISDA'05), p.192-196 (2005)

\bibitem[LiSte03]{LiSte03}
Na Li, M. Stephens:
{\it Modeling linkage disequilibrium and identifying recombination hotspots using single-nucleotide polymorphism data}. Genetics, 165 (4),  p.2213-33 (2003), doi: 10.1093/genetics/165.4.2213

\bibitem[LuYG95-QMC]{LuYG95-QMC}
Y.G. Lu:
{\it Quantum Markov chain and classical random sequences}.
Nagoya Math. J., 139, p.173--183 (1995)

\bibitem[MGK2021]{MGK2021}
B. Mor, S. Garhwal, A. Kumar:
{\it A systematic review of hidden markov models and their applications}. Archives of computational methods in engineering, 28 (3), p.1429-1448 (2021)

\bibitem[MonrWiesn11]{MonrWiesn11} A. Monras, A. Beige, and K. Wiesner: {\it Hidden Quantum Markov Models and non-adaptive read-out of many-body states},
App. Math. Comput. Sci. 3, 93 (2011)

\bibitem[Nguyet18]{Nguyet18} N. Nguyen: {\it Hidden Markov Model for Stock Trading}.
International Journal of Financial Studies, 6 (2), 36 (2018) https://doi.org/10.3390/ijfs6020036

\bibitem[OP]{OP} M. Ohya, D. Petz: {\bf Quantum entropy and its use}, Springer, Berlin-Heidelberg-New York (1993)

\bibitem[RosPent98]{RosPent98} N.M. Oliver, B. Rosario, A. Pentland: {\it Graphical models for recognizing human interactions}. Advances in Neural Information Processing Systems, 11, p.24-30  (1998)

\bibitem[PardBirm05]{PardBirm05}B. Pardo, W. Birmingham.
{\it Modeling Form for On-line Following of Musical Performances}. AAAI'05: Proceedings of the 20th national conference on Artificial intelligence, v.2, p.1018-1023 (2005)

\bibitem[Rab88]{Rab88} L.R. Rabiner:
 {\it Mathematical foundations of hidden Markov models}.
 In: Recent advances in speech understanding and dialog systems,
 Springer, Berlin, Heidelberg, p.183-205 (1988)

\bibitem[Rab89]{Rab89} L.R. Rabiner,
{\it A tutorial on hidden Markov models and selected applications in speech recognition},
Proceedings of the IEEE 77 (2), p.257--286 (1989)


\bibitem[RabJua86]{RabJua86}L.R. Rabiner, B.H. Juang:
{\it An introduction to hidden Markov models}.
IEEE ASSP magazine, 3 (1), p.4-16 (1986)

\bibitem[RabLeSo83]{RabLeSo83}
L.R. Rabiner, S.E. Levinsion, M.M. Sondhi,
{\it On the Application of Vector Quantization and Hidden Markov Models to Speaker-Independent Isolated Word Recognition},
Bell System Tech. J., v.62, n.4, p.1075--1105 (1983)

\bibitem[RaoSaz93-proj]{RaoSaz93-proj}
M.M. Rao, V.V. Sazonov: {\it A projective limit theorem for probability spaces and applications},
Theor. Probab. Appl. 38 (2),  p.307--315 (1993) Translated from Russian, https://doi.org/10.1137/1138027

\bibitem[RebSa17]{RebSa17}
Sara Rebagliati, Emanuela Sasso\\
Pattern recognition using hidden Markov models in financial time series\\
ACTA ET COMMENTATIONES UNIVERSITATIS TARTUENSIS DE MATHEMATICA, (21) 1 (2017) 1--17\\
http://acutm.math.ut.ee

\bibitem[SatGuru93]{SatGuru93}L. Satish, B.I. Gururaj,
 {\it Use of hidden Markov models for partial discharge pattern classification}, IEEE Transactions on Electrical Insulation, 28 (2), p.172-182 (1993)

\bibitem[SSGGBB18]{SSGGBB18}S. Srinivasan, G. Gordon, B. Boots: {\it Learning hidden quantum Markov models}.
Proceedings of the 21st International Conference on Artificial Intelligence and Statistics (AISTATS) 2018, PMLR: v.84, p.1979-1987

\bibitem[WiesnCrutc08]{WiesnCrutc08} K. Wiesner, C.P. Crutchfield:
 {\it Computation in finitary stochastic and quantum processes},
Physica D, v.237, iss.9, p.1173-1195 (2008), https://doi.org/10.1016/j.physd.2008.01.021

\bibitem[YamOhyIsh92]{YamOhyIsh92}J.
Yamato, J. Ohya, K. Ishii:
 {\it Recognizing human action in time-sequential images using hidden Markov model}.
 Proceedings 1992 IEEE Computer Society Conference on Computer Vision and Pattern Recognition (CVPR), p. 379-385, doi:10.1109/CVPR.1992


\end{thebibliography}
\end{document}